\input{epsf}

\newcounter{letter}     
\newenvironment{alphalist}{\begin{list}
{{\normalshape(\alph{letter})}}{\usecounter{letter}}}{\end{list}}                                                  
        \newcommand{\be}{\begin{equation}}
        \newcommand{\ee}{\end{equation}}
        \newcommand{\ba}{\begin{eqnarray}}
        \newcommand{\ea}{\end{eqnarray}}
        \newcommand{\ban}{\begin{eqnarray*}}
        \newcommand{\ean}{\end{eqnarray*}}
        \newcommand{\barr}{\begin{array}}
        \newcommand{\earr}{\end{array}}

\def\=>{\Rightarrow}

\def\mapright#1{\smash{\mathop{\longrightarrow}\limits^{#1}}}
\renewcommand{\fam}{{\rm fam}}
\newcommand{\svf}{{\rm svf}}
\newcommand{\prof}{{\rm prof}}
\newcommand{\type}{{\rm type}}

\newcommand{\sig}{{\rm sig}}
\newcommand{\et}{\hspace{-0.08in}{\bf .}\hspace{0.1in}}
\newcommand{\BOX}{\hbox {$\sqcap$ \kern -1em $\sqcup$}}
\newcommand{\qed}{\hskip 3em \hbox{\BOX} \vskip 2ex}
\renewcommand{\hom}{{\rm hom}}

\newcommand{\Set}{{\rm Set}}

\renewcommand{\to}{\rightarrow}
\newcommand{\tensor}{\otimes}
\newcommand{\maps}{\colon}
\newcommand{\op}{{\rm op}}

\newcommand{\iso}{\cong}
\newcommand{\End}{{\rm end}}
\newcommand{\elt}{{\rm elt}}

\newcommand{\N}{{\Bbb N}}
\newcommand{\G}{{\cal G}}

\newtheorem{thm}{Theorem}    

\newtheorem{prop}[thm]{Proposition}
\newtheorem{defn}[thm]{Definition}
\newtheorem{ex}[thm]{Example}

	\documentstyle[12pt,amsfonts,diagram]{article}
\begin{document}

\hfuzz=9pt      

      \begin{center}
      {\bf Higher-Dimensional Algebra III: \\
      $n$-Categories and the Algebra of Opetopes\\}
       \vspace{0.3cm}
      {\em John C.\ Baez and James Dolan\\}
      \vspace{0.3cm}
      {\small Department of Mathematics,  University of California\\ 
      Riverside, California 92521 \\
      USA\\ }
      \vspace{0.3cm}
      {\small email: baez@math.ucr.edu, jdolan@math.ucr.edu\\}
      \vspace{0.3cm}
      {\small February 5, 1997 \\ }
      \end{center}

\begin{abstract} We give a definition of weak $n$-categories
based on the theory of operads.  We work with operads having an
arbitrary set $S$ of types, or `$S$-operads', and given such an operad
$O$, we denote its set of operations by $\elt(O)$.  Then for any
$S$-operad $O$ there is an $\elt(O)$-operad $O^+$ whose algebras are
$S$-operads over $O$.  Letting $I$ be the initial operad with a
one-element set of types, and defining $I^{0+} = I$, $I^{(i+1)+} =
(I^{i+})^+$, we call the operations of $I^{(n-1)+}$ the
`$n$-dimensional opetopes'.  Opetopes form a category, and presheaves
on this category are called `opetopic sets'.  A weak $n$-category is
defined as an opetopic set with certain properties, in a manner
reminiscent of Street's simplicial approach to weak
$\omega$-categories.  In a similar manner, starting from an arbitrary
operad $O$ instead of $I$, we define `$n$-coherent $O$-algebras',
which are $n$ times categorified analogs of algebras of $O$.  Examples
include `monoidal $n$-categories', `stable $n$-categories', `virtual
$n$-functors' and `representable $n$-prestacks'.  We also describe how
$n$-coherent $O$-algebra objects may be defined in any
$(n+1)$-coherent $O$-algebra.
\end{abstract}

\section{Introduction}

A fundamental problem in higher-dimensional algebra is to set up a
convenient theory of weak $n$-categories.   Since there seems to be
quite a bit of freedom in what such a theory could look like, we begin
with a rough sketch of what is called for, and then summarize the
ideas behind our approach.
  
As traditionally conceived, an $n$-category should be some sort of
algebraic structure having objects or 0-morphisms, 1-morphisms between
0-morphisms, 2-morphisms between 1-morphisms, and so on up to
$n$-morphisms.  There should be various ways of composing
$j$-morphisms, and these composition operations should satisfy various
laws, such as associativity laws.  In the so-called `strict'
$n$-categories, these laws are equations.  While well-understood and
tractable, strict $n$-categories are insufficiently general for many
applications: what one usually encounters in nature are `weak'
$n$-categories, in which composition operations satisfy the
appropriate laws {\it only up to equivalence}.  Here the idea is that
$n$-morphisms are equivalent precisely when they are equal, while for
$j < n$ an equivalence between $j$-morphisms is recursively defined as
a $(j+1)$-morphism from one to the other that is invertible up to
equivalence.

What makes it difficult to define weak $n$-categories is that laws
formulated as equivalences should satisfy laws of their own ---
so-called `coherence laws' --- so that one can manipulate them with some
of the same facility as equations.  Moreover, these coherence laws
should also be equivalences satisfying their own coherence laws, again
up to equivalence, and so on.  

For example, a weak 1-category is just an ordinary category, defined
by Eilenberg and MacLane \cite{EM} in their 1945 paper.   In a
category, composition of 1-morphisms is associative `on the nose':
\[               (fg)h = f(gh)  .\]
Weak 2-categories first appeared in the work of B\'enabou
\cite{Benabou} in 1967, under the name of `bicategories'.  In a
bicategory, composition of 1-morphisms is associative only up to an
invertible 2-morphism, the `associator': 
\[              A_{f,g,h} \maps (fg)h \to f(gh) .\]
The associator allows one to rebracket parenthesized composites of
arbitrarily many 1-morphisms, but there may be many ways to use it to
go from one parenthesization to another.  For all these to be equal,
the associator must satisfy a coherence law, the pentagon identity, 
which says that the following diagram commutes:
\[
\begin{diagram}[(fg)(hi)]
\node{((fg)h)i} \arrow{e} \arrow{s} 
\node{(fg)(hi)} \arrow{e}
\node{f(g(hi))} \\
\node{(f(gh))i} \arrow[2]{e} \node[2]{f((gh)i)} \arrow{n} 
\end{diagram}
\]
where all the arrows are 2-morphisms built using the associator.
Weak 3-categories or `tricategories' were defined by Gordon, Power
and Street \cite{GPS} in a paper that appeared in 1995.   In a 
tricategory, the pentagon identity holds only up to an invertible
3-morphism, which satisfies a further coherence law of its own.

When one explicitly lists the coherence laws this way, the definition
of weak $n$-category tends to grow ever more complicated with
increasing $n$.   To get around this, one must carefully study the
origin of these coherence laws.   So far, most of our insight into
coherence laws has been won through homotopy theory, where it is
common to impose equations {\it only up to homotopy}, with these
homotopies satisfying coherence laws, again up to homotopy, and so on.
For example, the pentagon identity and higher coherence laws for
associativity first appeared in Stasheff's \cite{Stasheff} work on the
structure inherited by a space equipped with a homotopy equivalence
to a space with an associative product.  Subsequent work by Boardman
and Vogt, May, Segal and others led to a systematic treatment of
coherence laws in homotopy theory through the formalism of topological
operads \cite{Adams}.   

Underlying the connection between homotopy theory and $n$-category
theory is a hypothesis made quite explicit by Grothendieck \cite{Gro}: to any
topological space one should be able to associate an $n$-category
having points as objects, paths between points as 1-morphisms, certain
paths of paths as 2-morphisms, and so on, with certain
homotopy classes of $n$-fold paths as $n$-morphisms.  This should be a
special sort of weak $n$-category called a `weak $n$-groupoid', in which all
$j$-morphisms ($0 < j \le n$) are equivalences.  Moreover, the
process of assigning to each space its `fundamental $n$-groupoid', as
Grothendieck called it, should set up a complete correspondence
between the theory of homotopy $n$-types (spaces whose homotopy
groups vanish above the $n$th) and the theory of weak $n$-groupoids.  
This hypothesis explains why all the coherence laws for weak
$n$-groupoids should be deducible from homotopy theory.  It also
suggests that weak $n$-categories will have features not found in
homotopy theory, owing to the presence of $j$-morphisms that are not
equivalences.  

In addition, this hypothesis makes it clear in which
contexts the laws governing composition of $j$-morphisms should hold
only up to equivalence: namely, in those where {\it there is no preferred
composite of $j$-morphisms; instead, the composite is best
regarded as only unique up to equivalence}.   In homotopy theory this
arises from the arbitrary choice involved in parametrizing the
composite of two paths.   Because of this arbitrariness, 
composition of paths fails to be associative `on the nose'.  Instead,
it is associative up to a homotopy, the associator, with this homotopy
satisfying a coherence law, the pentagon identity, but again only
up to homotopy, and so on.

While many ways around this problem have been explored, here we prefer
to accept it as a fact of nature and develop a theory of weak
$n$-categories in which composition of $j$-morphisms is not an
operation in the traditional sense, but something a bit more subtle. 
Indeed, many forms of `composition' in mathematics are of this sort,
such as the disjoint union of sets or the tensor product of vector
spaces.  While one can artificially treat them as operations in the
traditional sense, it is better to define them by universal
properties.   Uniqueness up to equivalence then follows automatically.  
Taking this as a hint, we shall define the composite of $j$-morphisms
by a universal property.  

Homotopy theory also makes it clear that when setting up a theory of
$n$-categories, there is some choice involved in the shapes of ones
$j$-morphisms --- or in the language of topology, `$j$-cells'.  The
traditional approach to $n$-categories is `globular'.  This means that
for $j > 0$, each $j$-cell $f\maps x \to y$ has two $(j-1)$-cells
called its `source', $sf = x$, and `target', $tf = y$, which
for $j > 1$ satisfy
\[  s(sf) = s(tf), \qquad t(sf) = t(tf)).\]
Thus a $j$-cell can be visualized as a `globe', a $j$-dimensional ball whose
boundary is divided into two $(j-1)$-dimensional hemispheres corresponding
to its source and target.
In homotopy theory, however, the simplicial approach is much more
popular.  In a `simplicial set', each $j$-cell $f$ is
shaped like a $j$-dimensional simplex, and has $j+1$ faces, certain
$(j-1)$-cells $d_0f, \dots, d_nf$.  In addition to these there are
$(j+1)$-cells $i_0f, \dots, i_{n+1}f$ called `degeneracies', and 
the face and degeneracy maps satisfy certain well-known relations.

In the simplicial approach, weak $n$-groupoids are described using
`Kan complexes'.  It is worth recalling these here, because they begin
to illustrate how composite $j$-morphisms can be defined by a
universal property.  A `$j$-dimensional horn' in a simplicial set is,
roughly speaking, a configuration in which all but one of the faces of
a $j$-simplex have been filled in by $(j-1)$-cells in a consistent
way.  A simplicial set for which any horn can be extended to a
$j$-cell is called a `Kan complex'.  Kan complexes serve to describe
arbitrary homotopy types.  Algebraically, we may think of them as a
simplicial version of `weak $\omega$-groupoids', since they can have
nontrivial $j$-cells for arbitrarily large $j$.  A Kan complex
represents a homotopy $n$-type, or in other words a weak $n$-groupoid,
if for $j > n + 1$ any configuration in which all the faces of a
$j$-simplex have been filled in by $(j-1)$-cells in a consistent way
can be uniquely extended to a $j$-cell.   

Consider for example the
case $j = 2$.   Suppose, as shown in Figure 1, that two faces of a
2-simplex have been filled in by 1-cells $f$ and $g$ such that $d_1 f
= d_0 g = y$.   Then in a Kan complex we can extend this horn to a 2-cell
$F$, which has as its third face a 1-cell $h$.  
\medskip

\begin{center}
\setlength{\unitlength}{0.0100in}%

\begingroup\makeatletter\ifx\SetFigFont\undefined
\def\x#1#2#3#4#5#6#7\relax{\def\x{#1#2#3#4#5#6}}%
\expandafter\x\fmtname xxxxxx\relax \def\y{splain}%
\ifx\x\y   
\gdef\SetFigFont#1#2#3{%
  \ifnum #1<17\tiny\else \ifnum #1<20\small\else
  \ifnum #1<24\normalsize\else \ifnum #1<29\large\else
  \ifnum #1<34\Large\else \ifnum #1<41\LARGE\else
     \huge\fi\fi\fi\fi\fi\fi
  \csname #3\endcsname}%
\else
\gdef\SetFigFont#1#2#3{\begingroup
  \count@#1\relax \ifnum 25<\count@\count@25\fi
  \def\x{\endgroup\@setsize\SetFigFont{#2pt}}%
  \expandafter\x
    \csname \romannumeral\the\count@ pt\expandafter\endcsname
    \csname @\romannumeral\the\count@ pt\endcsname
  \csname #3\endcsname}%
\fi
\fi\endgroup
\begin{picture}(141,120)(228,422)
\thicklines

\put(241,441){\circle*{6}}
\put(360,440){\circle*{6}}
\put(240,440){\line( 3, 4){ 60}}
\put(300,520){\line( 3,-4){ 60}}
\multiput(260,475)(0.50000,0.25000){21}{\makebox(0.4444,0.6667){\SetFigFont{7}{8.4}{rm}.}}
\put(270,480){\line( 0,-1){ 10}}
\put(336,482){\line( 0,-1){ 10}}
\multiput(336,472)(-0.50000,0.25000){21}{\makebox(0.4444,0.6667){\SetFigFont{7}{8.4}{rm}.}}
\multiput(240,440)(11.42857,0.00000){11}{\line( 1, 0){  5.714}}
\multiput(298,445)(0.50000,-0.25000){21}{\makebox(0.4444,0.6667){\SetFigFont{7}{8.4}{rm}.}}
\multiput(308,440)(-0.50000,-0.25000){21}{\makebox(0.4444,0.6667){\SetFigFont{7}{8.4}{rm}.}}
\multiput(300,460)(-0.25000,0.50000){21}{\makebox(0.4444,0.6667){\SetFigFont{7}{8.4}{rm}.}}
\put(300,520){\circle*{6}}
\multiput(300,460)(0.25000,0.50000){21}{\makebox(0.4444,0.6667){\SetFigFont{7}{8.4}{rm}.}}
\multiput(302,479)(0.00000,-10.00000){2}{\line( 0,-1){  5.000}}
\multiput(298,479)(0.00000,-10.00000){2}{\line( 0,-1){  5.000}}

\put(280,465){\makebox(0,0)[lb]{\raisebox{0pt}[0pt][0pt]{$F$}}}
\put(228,431){\makebox(0,0)[lb]{\raisebox{0pt}[0pt][0pt]{$x$}}}
\put(299,527){\makebox(0,0)[lb]{\raisebox{0pt}[0pt][0pt]{$y$}}}
\put(369,433){\makebox(0,0)[lb]{\raisebox{0pt}[0pt][0pt]{$z$}}}
\put(261,481){\makebox(0,0)[lb]{\raisebox{0pt}[0pt][0pt]{$f$}}}
\put(338,484){\makebox(0,0)[lb]{\raisebox{0pt}[0pt][0pt]{$g$}}}
\put(301,422){\makebox(0,0)[lb]{\raisebox{0pt}[0pt][0pt]{$h$}}}
\end{picture}
\end{center}

\centerline{1.  Extending a horn to a cell}
\medskip

\noindent In this situation, we may think of $h$ as `a composite' of
$f$ and $g$, and $F$ as `a process of composing' $f$ and $g$.  There
is not a unique preferred composite.  However, it automatically
follows from the definition of Kan complex that any two composites are
equivalent.  Here two $j$-cells with all the same faces are said to be
`equivalent' if there is a $(j+1)$-cell having these $j$-cells as two
of its faces, the rest being degenerate.

Kan complexes serve as a highly efficient formalism in which to do
homotopy theory \cite{May2}.  In particular, there is no need to
explicitly list coherence laws!  They are all implicit in the fact that
every horn can be extended to a cell, and they all become explicit ---
in their simplicial forms --- if one makes composition into
an operation of the traditional sort by arbitrarily choosing an
extension of every horn.  It is tempting, therefore, to develop a
simplicial approach to weak $n$-categories.   

This was done by Street \cite{Street}, who actually dealt with weak
$\omega$-categories.  Like Kan complexes, these are simplicial sets.
However, only certain `admissible' horns, having the correct sort of
orientation, are required to have extensions.  For example, we do not
require the horn shown in Figure 2 to have an extension, since the
missing face would correspond to a composite of $f$ and an inverse of
$g$, which we expect to exist in a weak $n$-groupoid, but not in a
weak $n$-category.
\medskip

\begin{center}
\setlength{\unitlength}{0.0100in}%
\begingroup\makeatletter\ifx\SetFigFont\undefined
\def\x#1#2#3#4#5#6#7\relax{\def\x{#1#2#3#4#5#6}}%
\expandafter\x\fmtname xxxxxx\relax \def\y{splain}%
\ifx\x\y   
\gdef\SetFigFont#1#2#3{%
  \ifnum #1<17\tiny\else \ifnum #1<20\small\else
  \ifnum #1<24\normalsize\else \ifnum #1<29\large\else
  \ifnum #1<34\Large\else \ifnum #1<41\LARGE\else
     \huge\fi\fi\fi\fi\fi\fi
  \csname #3\endcsname}%
\else
\gdef\SetFigFont#1#2#3{\begingroup
  \count@#1\relax \ifnum 25<\count@\count@25\fi
  \def\x{\endgroup\@setsize\SetFigFont{#2pt}}%
  \expandafter\x
    \csname \romannumeral\the\count@ pt\expandafter\endcsname
    \csname @\romannumeral\the\count@ pt\endcsname
  \csname #3\endcsname}%
\fi
\fi\endgroup
\begin{picture}(141,120)(228,422)
\thicklines
\put(241,441){\circle*{6}}
\put(360,440){\circle*{6}}
\put(240,440){\line( 3, 4){ 60}}
\put(300,520){\line( 3,-4){ 60}}
\multiput(260,475)(0.50000,0.25000){21}{\makebox(0.4444,0.6667){\SetFigFont{7}{8.4}{rm}.}}
\put(270,480){\line( 0,-1){ 10}}
\multiput(240,440)(11.42857,0.00000){11}{\line( 1, 0){  5.714}}
\multiput(298,445)(0.50000,-0.25000){21}{\makebox(0.4444,0.6667){\SetFigFont{7}{8.4}{rm}.}}
\multiput(308,440)(-0.50000,-0.25000){21}{\makebox(0.4444,0.6667){\SetFigFont{7}{8.4}{rm}.}}
\multiput(300,460)(-0.25000,0.50000){21}{\makebox(0.4444,0.6667){\SetFigFont{7}{8.4}{rm}.}}
\multiput(300,460)(0.25000,0.50000){21}{\makebox(0.4444,0.6667){\SetFigFont{7}{8.4}{rm}.}}
\put(300,520){\circle*{6}}
\multiput(302,479)(0.00000,-10.00000){2}{\line( 0,-1){  5.000}}
\multiput(298,479)(0.00000,-10.00000){2}{\line( 0,-1){  5.000}}
\multiput(340,475)(-0.50000,0.25000){21}{\makebox(0.4444,0.6667){\SetFigFont{7}{8.4}{rm}.}}
\put(330,480){\line( 0,-1){ 10}}
\put(280,465){\makebox(0,0)[lb]{\raisebox{0pt}[0pt][0pt]{$F$}}}
\put(228,431){\makebox(0,0)[lb]{\raisebox{0pt}[0pt][0pt]{$x$}}}
\put(299,527){\makebox(0,0)[lb]{\raisebox{0pt}[0pt][0pt]{$y$}}}
\put(369,433){\makebox(0,0)[lb]{\raisebox{0pt}[0pt][0pt]{$z$}}}
\put(261,481){\makebox(0,0)[lb]{\raisebox{0pt}[0pt][0pt]{$f$}}}
\put(338,484){\makebox(0,0)[lb]{\raisebox{0pt}[0pt][0pt]{$g$}}}
\put(301,422){\makebox(0,0)[lb]{\raisebox{0pt}[0pt][0pt]{$h$}}}
\end{picture}
\end{center}

\centerline{2.  A horn that need not have an extension in a $n$-category}
\medskip

Second, for those horns that are required to have extensions, we
require the existence of a `universal' extension.  The point is that,
unlike $\omega$-groupoid case, we cannot think of every $(j+1)$-cell
as a `process of composing' all but one of its faces to obtain the
remaining face.  Instead, Street's weak $\omega$-categories are
equipped with a distinguished set of `universal' cells which we can
think of this way.  These satisfy some axioms: there are no universal
0-cells, all universal 1-cells are degenerate, and all degenerate
cells are universal.  Last, and most importantly, any composite of
universal cells is universal.

Our definition of weak $n$-categories resembles Street's, but with two
major differences.  First, while simplices are convenient in algebraic
topology, they are not well adapted to the `unidirectional' or
`noninvertible' character of the $j$-morphisms in $n$-category theory,
as is clear from the rather technical combinatorics involved in
orienting the faces of a simplex and defining admissible horns.  This
raises the possibility that a more convenient theory could be set up
with $j$-cells of some other shapes --- shapes motivated more by the
inner logic of $n$-category theory than by traditional concerns of
algebraic topology.  In our approach we use certain shapes called
`opetopes'.  ({\it Nota bene}: The first two syllables of `opetope'
are pronounced exactly as in the word `operation'.)

Opetopes arise naturally from the theory of operads.   Roughly
speaking, an `operad' is an algebraic gadget consisting of a
collection of abstract operations closed under composition. These
operations may have any finite number of arguments, and we work with
operads in which the arguments are `typed' or many-sorted.  Any such
operad is determined by: 1) its types, 2) its operations, and 3) its
`reduction laws', or equations stating that some composite of
operations equals a given operations.  This description of an operad
is like a presentation in terms of generators and relations.  An operad
also has `algebras' in which its operations are represented as actual
functions.   From the viewpoint of mathematical logic, an operad is a
kind of theory, and its algebras are models of that theory.  As
always, it is useful to study operads both syntactically, in terms
of their presentations, and semantically, in terms of their algebras.  

We define the `slice operad' $O^+$ of an operad $O$ in
such a way that an algebra of $O^+$ is precisely an operad over $O$, i.e.,
an operad with the same set of types as $O$, equipped with an operad 
homomorphism to $O$.   Syntactically, it turns out that:
\begin{enumerate}
\item The types of $O^+$ are the operations of $O$.
\item The operations of $O^+$ are the reduction laws of $O$.
\item The reduction laws of $O^+$ are the ways of combining reduction
laws of $O$ to give other reduction laws.  
\end{enumerate}
This gets at the heart of the process of `categorification', in which
laws are promoted to operations and these operations satisfy new
coherence laws of their own.  Here the coherence laws arise simply from
the ways of combining the the old laws.

The simplest operad of all is the initial operad, $I$.  Syntactically
speaking, this is the operad with only one type and only one
operation, the identity.  Semantically, $I$ is the operad whose
algebras are just sets, without any extra structure at all.  Starting
with $I$ and iterating the slice operad construction $(j-1)$ times, we
obtain an operad whose operations we call `$j$-dimensional opetopes.'
A $0$-dimensional opetope is just a point, and a $1$-dimensional
opetope is just an oriented interval.  For $j > 1$, a $j$-dimensional
opetope may have any number of `infaces' but only one
`outface'.  Thanks to the above syntactic description of the slice
operad construction, it turns out that a $j$-dimensional opetope
corresponds simply to a way of pasting together its
infaces --- certain $(j-1)$-dimensional opetopes --- to 
obtain its outface.  

A weak $n$-category will be an `opetopic set' with certain extra
properties, similar to those defining a Kan complex, but a bit more
complicated.  The analog of an admissible horn is a `niche', which
is a configuration in which all the infaces of an opetope have
been filled in with cells, but not the outface.  We require that
every niche can be extended to a universal cell, and regard the outface
of such a universal cell as `a composite' of its infaces.
We also require composites of universal cells to be universal.

Here we must note the second major difference between our approach and
Street's.  It turns out that if one works with $n$-categories instead
of $\omega$-categories, one need not (and should not) arbitrarily
designate certain cells as universal; instead, universality becomes a
property.  In our framework an $n$-category typically has $j$-cells
for arbitrarily large $j$, but they act like `equations' for $j > n$,
so every $j$-cell is defined to be universal for $j > n$.
Universality for $j$-cells of lower dimension is defined in a
recursive manner.  The basic idea is that a given cell occupying some
niche is universal if any other occupant factors through that one, up
to equivalence.  Here the notion of `equivalence' must also be
recursively defined.

A brief outline of our paper is as follows.  In Section \ref{operads}
we introduce some necessary material on operads.  In Section
\ref{opetopes} we describe the slice operad construction, opetopes and
opetopic sets.  In Section \ref{n-categories} we define weak
$n$-categories and begin to study them, along with the more general
`$n$-coherent $O$-algebras', which are the $n$-categorical analogs of
operad algebras.  In the Conclusions we compare other approaches to
weak $n$-categories and discuss the all-important question of when two
approaches can be considered equivalent.

Henceforth by `$n$-category' we always mean `weak $n$-category', as
defined in this paper.   For more background on $n$-category theory
and why it should be interesting, see our previous papers, which we
refer to as HDA0 \cite{BD}, HDA1 \cite{BN}, and HDA2 \cite{B}.  As in
those papers, we use the ordering in which the composite of morphisms
$f \maps x \to y$ and $g \maps y \to z$ is written as $fg$, but when
dealing with operads we write the composite of a $k$-ary operation $f$
with the operations $g_1,\dots,g_k$ as $f \cdot (g_1,\dots,g_k)$.    

\section{Operads} \label{operads}

It turns out to be convenient to describe weak $n$-categories using the
theory of operads.  Operads are a formalism for dealing with algebraic
structures having operations of arbitrary finite arity satisfying
arbitrary `reduction laws', that is, equational laws saying that some
composite of operations equals some operation.  For the benefit of the
reader unfamiliar with operads, we begin in Section
\ref{untyped.operads} by recalling the traditional sort of operad
\cite{May}.  We call these `untyped' operads because they are suited to
the case when the inputs and output of every operation are of the same
type.  Then, with the help of some generalities about monoid objects in
Section \ref{monoid.objects}, we introduce the more general `typed'
operads needed for this paper in Section \ref{typed.operads}.  While all
we really need are operads with an arbitrary {\it set} of types, we find
it somewhat illuminating to define operads with an arbitrary small {\it
category} of types.  In Section \ref{pullback.operads} we show how a
functor $F \maps C \to D$ gives a way to turn operads with type category
$D$ into operads with type category $C$, and in Section
\ref{slice.operads} we conclude with another basic operad construction,
the `slice operad of an operad algebra'.

\subsection{Untyped operads} \label{untyped.operads}

An {\it untyped operad} $O$ has, for each $k \ge 0$, a set $O_k$ of
{\it $k$-ary operations}.  We may visualize an element of $O_k$ as a
tree as in Figure 3.  This tree has one black dot or {\it node}
representing the operation itself, $k$ lines or {\it edges} coming in
from above representing the inputs of the operation, 
and one edge going out from below representing the output.  

\vbox{

\begin{center}
\setlength{\unitlength}{0.000500in}%
\begingroup\makeatletter\ifx\SetFigFont\undefined
\def\x#1#2#3#4#5#6#7\relax{\def\x{#1#2#3#4#5#6}}%
\expandafter\x\fmtname xxxxxx\relax \def\y{splain}%
\ifx\x\y   
\gdef\SetFigFont#1#2#3{%
  \ifnum #1<17\tiny\else \ifnum #1<20\small\else
  \ifnum #1<24\normalsize\else \ifnum #1<29\large\else
  \ifnum #1<34\Large\else \ifnum #1<41\LARGE\else
     \huge\fi\fi\fi\fi\fi\fi
  \csname #3\endcsname}%
\else
\gdef\SetFigFont#1#2#3{\begingroup
  \count@#1\relax \ifnum 25<\count@\count@25\fi
  \def\x{\endgroup\@setsize\SetFigFont{#2pt}}%
  \expandafter\x
    \csname \romannumeral\the\count@ pt\expandafter\endcsname
    \csname @\romannumeral\the\count@ pt\endcsname
  \csname #3\endcsname}%
\fi
\fi\endgroup
\begin{picture}(1844,2144)(3279,-3083)
\thicklines
\put(4201,-1861){\circle*{150}}
\put(4201,-3061){\line( 0, 1){1200}}
\put(4201,-1861){\line(-1, 3){300}}
\put(4201,-1861){\line( 1, 3){300}}
\put(4201,-1861){\line( 1, 1){900}}
\put(4201,-1861){\line(-1, 1){900}}
\end{picture}
\end{center}

\centerline{3.  An element of $O_k$ for $k = 4$}
\medskip
}

\noindent We may compose these trees by attaching the output edges of $k$
of them to the input edges of a tree with $k$ inputs, as shown in Figure
4.   (In the resulting tree some of the edges may be drawn as broken
lines for convenience.)  More precisely, for any integers $i_1, \dots,
i_k \ge 0$ there is a function  
\ban    O_k \times O_{i_1} \times \cdots \times
O_{i_k} &\to& O_{i_1 + \cdots + i_k}  . \\
(f,g_1, \dots, g_k) &\mapsto& f \cdot (g_1, \dots , g_k)  \ean

\vbox{

\begin{center}
\setlength{\unitlength}{0.0002500in}%
\begingroup\makeatletter\ifx\SetFigFont\undefined
\def\x#1#2#3#4#5#6#7\relax{\def\x{#1#2#3#4#5#6}}%
\expandafter\x\fmtname xxxxxx\relax \def\y{splain}%
\ifx\x\y   
\gdef\SetFigFont#1#2#3{%
  \ifnum #1<17\tiny\else \ifnum #1<20\small\else
  \ifnum #1<24\normalsize\else \ifnum #1<29\large\else
  \ifnum #1<34\Large\else \ifnum #1<41\LARGE\else
     \huge\fi\fi\fi\fi\fi\fi
  \csname #3\endcsname}%
\else
\gdef\SetFigFont#1#2#3{\begingroup
  \count@#1\relax \ifnum 25<\count@\count@25\fi
  \def\x{\endgroup\@setsize\SetFigFont{#2pt}}%
  \expandafter\x
    \csname \romannumeral\the\count@ pt\expandafter\endcsname
    \csname @\romannumeral\the\count@ pt\endcsname
  \csname #3\endcsname}%
\fi
\fi\endgroup
\begin{picture}(4544,4244)(3279,-5183)
\thicklines
\put(4201,-1861){\circle*{150}}
\put(5701,-1861){\circle*{150}}
\put(7201,-1861){\circle*{150}}
\put(5701,-3661){\circle*{150}}
\put(4201,-3061){\line( 0, 1){1200}}
\put(4201,-1861){\line(-1, 3){300}}
\put(4201,-1861){\line( 1, 3){300}}
\put(4201,-1861){\line( 1, 1){900}}
\put(4201,-1861){\line(-1, 1){900}}
\put(4201,-3061){\line( 5,-2){1500}}
\put(5701,-3661){\line( 0, 1){600}}
\put(5701,-3061){\line( 0, 1){1200}}
\put(5701,-3661){\line( 5, 2){1500}}
\put(7201,-3061){\line( 0, 1){1200}}
\put(7201,-1861){\line(-2, 3){600}}
\put(7201,-1861){\line( 0, 1){900}}
\put(7201,-1861){\line( 2, 3){600}}
\put(5701,-3661){\line( 0,-1){1500}}
\end{picture}
\end{center}

\centerline{4.  Composition in an operad}
\medskip
}

\noindent We require composition to be `associative', in the sense that
\[      f \cdot (g_1 \cdot (h_{11}, \dots, h_{1i_1}), \dots, 
g_k \cdot (h_{k1}, \dots, h_{ki_k})) = \]
\[ (f \cdot (g_1, \dots g_k)) \; \cdot \; (h_{11}, \dots, h_{1i_1},
\dots\dots , h_{k1}, \dots, h_{ki_k}) \]
whenever both sides are well-defined.  
This makes composites such as those shown in Figure 5 
unambiguous.  

\vbox{

\begin{center}

\setlength{\unitlength}{0.000250in}%
\begingroup\makeatletter\ifx\SetFigFont\undefined
\def\x#1#2#3#4#5#6#7\relax{\def\x{#1#2#3#4#5#6}}%
\expandafter\x\fmtname xxxxxx\relax \def\y{splain}%
\ifx\x\y   
\gdef\SetFigFont#1#2#3{%
  \ifnum #1<17\tiny\else \ifnum #1<20\small\else
  \ifnum #1<24\normalsize\else \ifnum #1<29\large\else
  \ifnum #1<34\Large\else \ifnum #1<41\LARGE\else
     \huge\fi\fi\fi\fi\fi\fi
  \csname #3\endcsname}%
\else
\gdef\SetFigFont#1#2#3{\begingroup
  \count@#1\relax \ifnum 25<\count@\count@25\fi
  \def\x{\endgroup\@setsize\SetFigFont{#2pt}}%
  \expandafter\x
    \csname \romannumeral\the\count@ pt\expandafter\endcsname
    \csname @\romannumeral\the\count@ pt\endcsname
  \csname #3\endcsname}%
\fi
\fi\endgroup
\begin{picture}(4980,5444)(3218,-5183)
\thicklines
\put(4201,-1861){\circle*{150}}
\put(5701,-1861){\circle*{150}}
\put(7201,-1861){\circle*{150}}
\put(5701,-3661){\circle*{150}}
\put(3901,-361){\circle*{150}}
\put(4501,-361){\circle*{150}}
\put(5101,-361){\circle*{150}}
\put(6601,-361){\circle*{150}}
\put(7801,-361){\circle*{150}}
\put(3301,-361){\circle*{150}}
\put(7201,-361){\circle*{150}}
\put(4201,-3061){\line( 0, 1){1200}}
\put(4201,-1861){\line(-1, 3){300}}
\put(4201,-1861){\line( 1, 3){300}}
\put(4201,-1861){\line( 1, 1){900}}
\put(4201,-1861){\line(-1, 1){900}}
\put(4201,-3061){\line( 5,-2){1500}}
\put(5701,-3661){\line( 0, 1){600}}
\put(5701,-3061){\line( 0, 1){1200}}
\put(5701,-3661){\line( 5, 2){1500}}
\put(7201,-3061){\line( 0, 1){1200}}
\put(7201,-1861){\line(-2, 3){600}}
\put(7201,-1861){\line( 0, 1){900}}
\put(7201,-1861){\line( 2, 3){600}}
\put(5701,-3661){\line( 0,-1){1500}}
\put(3301,-961){\line( 0, 1){600}}
\put(3901,-961){\line( 0, 1){600}}
\put(4501,-961){\line( 0, 1){600}}
\put(5101,-961){\line( 0, 1){600}}
\put(6601,-961){\line( 0, 1){600}}
\put(7201,-1036){\line( 0, 1){675}}
\put(7801,-961){\line( 0, 1){600}}
\put(5101,-361){\line( 1, 4){150}}
\put(5101,-361){\line(-1, 4){150}}
\put(4501,-361){\line( 0, 1){600}}
\put(3901,-361){\line( 1, 2){300}}
\put(3901,-361){\line( 0, 1){600}}
\put(3901,-361){\line(-1, 2){300}}
\put(6601,-361){\line( 1, 2){300}}
\put(6601,-361){\line( 0, 1){600}}
\put(6601,-436){\line(-2, 5){274.138}}
\put(7801,-361){\line(-3, 5){363.971}}
\put(7801,-361){\line(-1, 4){150}}
\put(7801,-361){\line( 1, 4){150}}
\put(7801,-361){\line( 3, 5){363.971}}
\end{picture}
\end{center}

\centerline{5.  Associativity for composition in an operad}
\medskip
}

\noindent We also require the existence of an `unit' $1 \in O_1$
such that
\[     1 \cdot (f) = f, \qquad f \cdot (1,\dots, 1) = f \]
for all $f \in O_k$.  

What we have so far is an {\it planar untyped operad}.   For a
full-fledged untyped operad, we also assume there is a right action of
the symmetric group $S_k$ on $O_k$ for all $k$, for which the
following compatibility condition holds: for any $f \in O_k$,
$\sigma \in S_k$, and $g_j \in O_{i_j}$ for $1 \le j \le k$, we
have
\[     (f\sigma) \cdot (g_1, \dots, g_k) = (f \cdot (g_1, \dots,
g_k))\, \rho(\sigma), \]
where
\[      \rho \maps S_k \to S_{i_1 + \cdots + i_k}  \]
is the obvious homomorphism.  We illustrate this condition in Figure 6.

\vbox{
\bigskip
\centerline{\epsfysize=1.5in\epsfbox{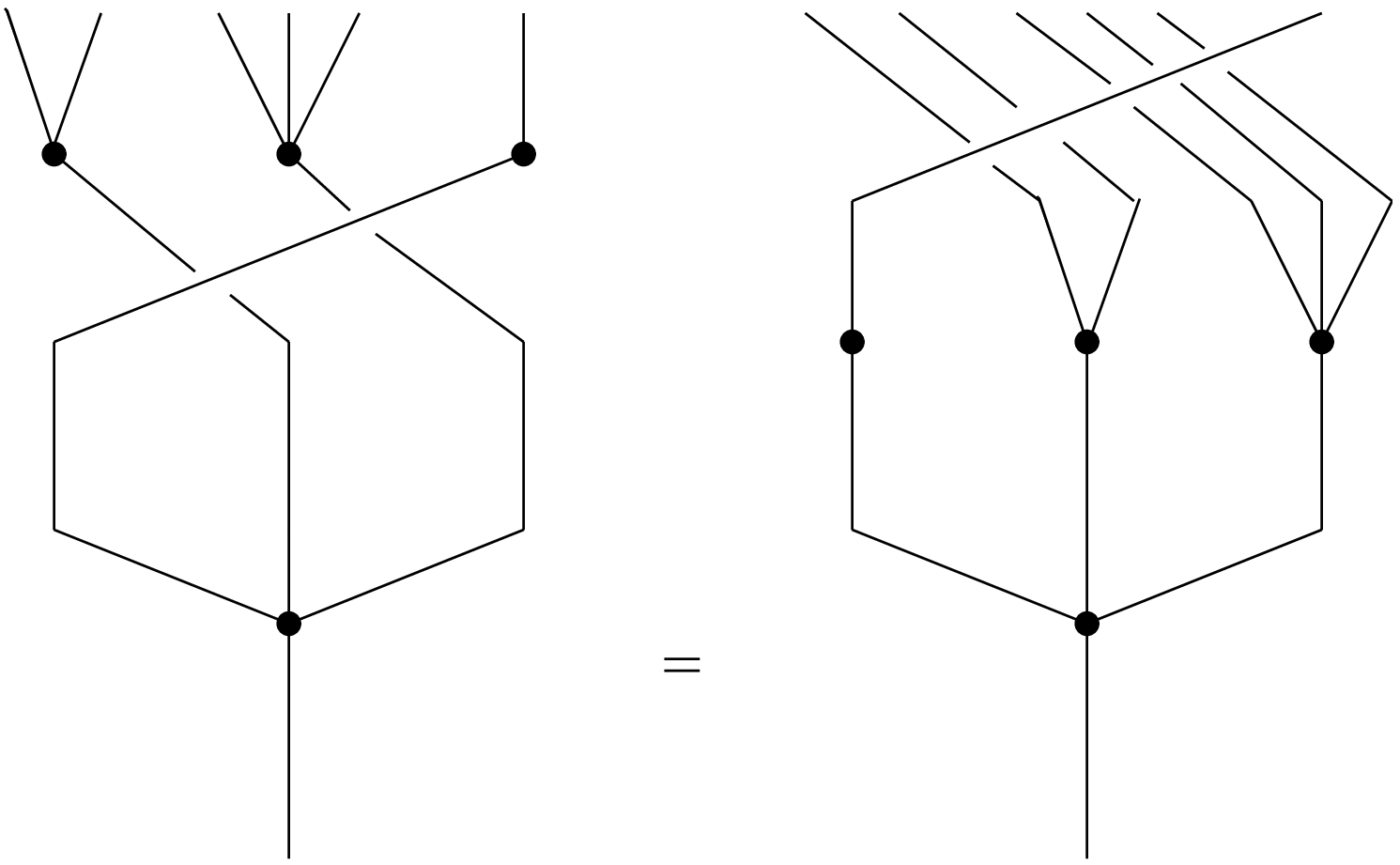}}
\medskip
\centerline{6.  Compatibility condition for symmetric group actions}
\medskip
}

Operads are mainly interesting for their algebras.   Given an untyped
operad $O$ as above, one defines an {\it $O$-algebra} to be a set $A$ on
which the operations of $O$ are rendered concrete.  In other words,
there are maps 
\[      \alpha \maps O_k \to \hom(A^k,A) \] 
sending the
identity operation $1 \in O_1$ to the identity function from $x$ to
itself, and sending composites to composites: 
\[      \alpha(f \cdot
(g_1, \dots, g_k)) = \alpha(f) \circ (\alpha(g_1) \times \cdots \times
\alpha(g_k)) .\] 
We also require that the maps $\alpha$ satisfy 
\[        \alpha(f\sigma) = \alpha(f)\sigma, \] 
where $f \in O_k$ and $\sigma \in
S_k$ acts on $\hom(A^k,A)$ on the right by permuting the factors in
$A^k$.  We omit this requirement if $O$ is merely planar.

\subsection{Monoid objects and their actions} \label{monoid.objects}

An untyped operad $O$ for which only $O_1$ is nonempty is just a
monoid, so we may think of a monoid as a kind of operad. 
Interestingly, however, there is a rather different way to think of
any operad as a kind of monoid, or more precisely, a `monoid object'. 
By the internalization principle discussed in HDA0, we can generalize
the definition of `monoid' from the category of sets to any
sufficiently similar category --- in fact, any monoidal category.  If
$M$ is a strict monoidal category, a {\it monoid object} in $M$ is an
object $m \in M$ equipped with a product $\mu \maps m \tensor m \to m$
and unit $\iota \maps 1 \to m$, such that the following diagrams
commute:
\[
\begin{diagram} [m \tensor m \tensor m] 
\node{m \tensor m \tensor m} 
\arrow{e,t}{\mu \tensor 1}  \arrow{s,l}{1 \tensor \mu} 
\node{m \tensor m}  \arrow{s,r}{\mu}  \\ 
\node{m \tensor m} \arrow{e,t}{\mu} \node{m}  
\end{diagram} 
\] 
\[ 
\begin{diagram} [m \tensor m] 
\node{1 \tensor m} \arrow[2]{e,t}{\iota \tensor 1} \arrow{se,b}{1} 
\node[2]{m \tensor m} \arrow{sw,r}{\mu}   \\
\node[2]{m}  
\end{diagram} 
\] 
\[ 
\begin{diagram} [m \tensor m]
\node{m \tensor 1} \arrow[2]{e,t}{1 \tensor \iota} \arrow{se,b}{1} 
\node[2]{m \tensor m}\arrow{sw,r}{\mu}   \\
\node[2]{m}  
\end{diagram} 
\] 
These represent associativity and the left and right unit laws,
respectively.  When the monoidal category $M$ is not strict, one simply
inserts the natural isomorphisms $(m \tensor m) \tensor m \iso m \tensor
(m \tensor m)$ and $1 \tensor m \iso m \iso m \tensor 1$ where needed. 

One may then define an action of the monoid object $m$ on any object in
$M$.  More generally, one may define an action of $m$ on any object in a
category on which $M$ acts.  Recall that an {\it action} of $M$ on a
category $C$ is a monoidal functor $A \maps M \to \End(C)$, where the
monoidal category $\End(C)$ has endofunctors on $C$ as objects and
natural transformations between these as morphisms.  Equivalently, we
may think of the action $A$ as a functor $A \maps M \times C \to C$
satisfying certain conditions.  Here it is convenient to write $A(m,c)$
simply as $m \tensor c$.

Suppose that $m \in M$ is an monoid object and $A \maps M \times C \to
C$ is an action.  If $A$ is a strict monoidal functor, we define an {\sl
action} of $m$ in $C$ {\it riding} the action $A$ to be a morphism
\[              \alpha \maps m \tensor c \to c \]
in $C$ making the following diagrams commute:
\[
\begin{diagram} [m \tensor m \tensor c]
\node{m \tensor m \tensor c} \arrow{e,t}{\mu \tensor 1} 
\arrow{s,l}{1 \tensor \alpha} \node{m \tensor c} 
\arrow{s,r}{\alpha}  \\
\node{m \tensor c} \arrow{e,t}{\alpha} \node{c} 
\end{diagram}
\]
\[
\begin{diagram} [m \tensor c]
\node{1 \tensor c} \arrow[2]{e,t}{i \tensor 1} \arrow{se,b}{1}
\node[2]{m \tensor c}\arrow{sw,r}{\alpha}   \\
\node[2]{c} 
\end{diagram}
\]
When $A$ is not strict, one simply inserts the natural isomorphisms
$(m \tensor m) \tensor c \iso m \tensor (m \tensor c)$ and 
$1 \tensor c \iso c$ where needed.

Given a monoid object $m \in M$ and an action $A$ of $M$ on $C$, we
define the {\it category of actions} of $m$ in $C$ {\it riding} the
action $A$ as follows.  The objects of this category are actions of $m$
in $C$ riding $A$, and given two such actions
\[             \alpha \maps m \tensor c \to c , \qquad
\alpha' \maps m \tensor c' \to c', \]
we define a morphism from $\alpha$ to $\alpha'$ to be a morphism
$f \maps c \to c'$ such that the following diagram commutes:
\[
\begin{diagram} [m \tensor c']
\node{m \tensor c} \arrow{e,t}{\alpha} \arrow{s,l}{1 \tensor f} 
\node{c} \arrow{s,r}{f}  \\
\node{m \tensor c'} \arrow{e,t}{\alpha'} \node{c'} 
\end{diagram}
\]

It is worth noting an interesting pattern.  We may generalize the
notion of `monoid' to that of a monoid object $m$ in any monoidal
category $M$, but the notion of `monoidal category' is itself a
categorification of the notion of `monoid'.  Similarly, we may define an
action of the monoid object $m \in M$ on $c \in C$ whenever the monoidal
category $M$ acts on $C$.

We see here two instances of the following principle: {\it certain
algebraic structures can be defined in any category equipped
with a categorified version of the same structure.}  Another instance
was mentioned in HDA0: we may define a commutative monoid object in
any symmetric monoidal category.  We name this principle the {\it
microcosm principle,} after the theory, common in pre-modern
correlative cosmologies, that every feature of the microcosm (e.g.\
the human soul) corresponds to some feature of the macrocosm.  Of
course, the above formulation of the microcosm principle is rather
vague; we give a precise version in Section \ref{microcosm}.

Even without a precise formulation, the microcosm principle can serve as
a useful guide when seeking the most general way to internalize certain
algebraic structures.  For example, we may apply the microcosm principle
to morphisms between monoid objects.  Suppose we are given a monoidal
functor $F \maps M \to M'$ between monoidal categories.  Then given
monoid objects $m \in M$ and $m' \in M'$, we define a morphism $f \maps
m \to m'$ {\it riding} $F$ to be a morphism $f \maps F(m) \to m'$ in
$M'$ making the following diagrams commute:
\[
\begin{diagram} [F(m \tensor m)]
\node{F(m \tensor m)} \arrow{s,l}{F(\mu)}
\node{F(m) \tensor F(m)} \arrow{e,t}{f\tensor f} \arrow{w}
\node{m' \tensor m'} \arrow{s,r}{\mu'} \\
\node{F(m)} \arrow[2]{e,t}{f} \node[2]{m'}  
\end{diagram}
\]
\[
\begin{diagram} [F(1_M)]
\node{F(1_M)} \arrow{s,l}{F(\iota)} 
\node{1_{M'}} \arrow{w} \arrow{s,r}{\iota'}\\
\node{F(m)} \arrow{e,t}{f} \node{m'} \\
\end{diagram}
\]

Here $\mu,\iota$ are the product and unit for $m$, while
$\mu',\iota'$ are the product and unit for $m'$.   If $F$ is a
strict monoidal functor, the unlabelled arrows $F(m) \tensor F(m)
\to F(m \tensor m)$ and $1_{M'} \to F(1_M)$ are identity
morphisms.  If $F$ is a weak monoidal functor, these arrows are
isomorphisms supplied by the definition of weak monoidal functor.	
However, in Sections \ref{pullback.operads} and
\ref{slice.operads} we will need the case where $F$ is merely a
lax monoidal functor (which Eilenberg and Kelly \cite{EK} call
simply a monoidal functor).   Then these arrows are morphisms,
not necessarily isomorphisms, supplied by the
definition of a lax monoidal functor.  

We call a morphism of monoid objects riding an identity functor
a `homomorphism'.  In other words, given
monoid objects $m,m'$ in a monoidal category $M$, we define a
{\it homomorphism} $f \maps m \to m'$ to be a morphism 
in $M$ for which the following diagrams commute:
\[
\begin{diagram} [m' \tensor m']
\node{m \tensor m} \arrow{e,t}{f \tensor f} \arrow{s,l}{\mu}
\node{m' \tensor m'} \arrow{s,r}{\mu'}  \\
\node{m} \arrow{e,t}{f} \node{m'} 
\end{diagram}
\]
\[
\begin{diagram} [m' \tensor m']
\node[2]{1} \arrow{sw,t}{\iota} \arrow{s,r}{\iota'} 
\\
\node{m} \arrow{e,t}{f} \node{m'} 
\end{diagram}
\]
Whenever $F\maps M \to M'$ is a lax monoidal functor and $m$ is a monoid
object in $M$, $F(m)$ becomes a monoid object in $M'$ in a natural way,
and a morphism of monoid objects $f \maps m \to m'$ riding $F$ can also
be thought of as a homomorphism from $F(m)$ to $m'$.

In the next section, we define an operad
to be a monoid object in a certain monoidal category of `signatures'.
An algebra of the operad will then be an action of this monoid
object, riding a certain action of the category of signatures.
We will use the concepts of morphisms and homomorphisms between
monoid objects to define morphisms and homomorphisms between operads.

\subsection{Typed operads} \label{typed.operads}

To define weak $n$-categories, we need operads for which the inputs and
output of each operation are `many-sorted', or `typed'.  In what
follows, we first define these operads and their algebras, and
then give a rather lengthy explanation of our definitions.  

\begin{defn}\et For a category $C$, let the category $\fam(C)$ of
{\rm $C$-families} be the category where an object is a finite list
of objects of $C$, and where a morphism from $(x_1,\dots,x_j)$ to 
$(y_1,\dots,y_k)$ is a bijection $b \maps \{1,\dots,j\} \to
\{1,\dots,k\}$ together with, for each $i$, a morphism from $x_i$ to
$y_{b(i)}$, with composition of morphisms given by the
obvious rule. \end{defn}

\begin{defn}\et For a category $C$, let $\svf(C)$ be the category of
set-valued functors on $C$, that is, the category whose objects are
functors from $C$ to $\Set$ and whose morphisms are natural
transformations between these.  \end{defn}

Notice that $\fam(C)$ is a symmetric monoidal category, where the tensor
product of the families $(x_1,\dots,x_j)$ and $(y_1,\dots,y_k)$ is
$(x_1,\dots,x_j,y_1,\dots,y_k)$.  In fact, $\fam(C)$ is the free
symmetric monoidal category on $C$ in an appropriate sense.  In a
similar sense, the category $\svf(\fam(C)^\op)$ of set-valued
contravariant functors on $\fam(C)$ is the `free symmetric 2-rig on
$C$', where a {\it 2-rig} is a symmetric monoidal cocomplete category
for which the monoidal structure preserves small colimits in each
argument.  By this universal property, the monoidal category
$\End(\svf(\fam(C)^\op))$ of endomorphisms preserving small colimits and
the symmetric monoidal structure is equivalent to $\svf(\fam(C)^\op
\times C)$.  

\begin{defn} \et Given a category $C$, we define the category $\prof(C)$
of {\rm $C$-profiles} to be $\fam(C)^\op \times C$. \end{defn}

\begin{defn} \et  Given a category $C$, we define the
category $\sig(C)$ of {\rm $C$-signatures} to be $\svf(\prof(C))$.  
\end{defn}

By the above remarks $\sig(C)$ is a monoidal category.  Note that
$\sig(C)$ has an action on $\svf(C)$, which we call the {\it tautologous
action}.   Thus we may make the following definitions:

\begin{defn}\et Given a small category $C$, we define a {\rm $C$-operad}
be a monoid object in $\sig(C)$, and define $\op(C)$, the category of
$C$-operads, to be the category of monoid objects in $\sig(C)$.
\end{defn}

\begin{defn} \et Given a $C$-operad $O$, we say $C$ is the category of
{\rm types} of $O$, and write $C = \type(O)$.  \end{defn}

\begin{defn} \et For a $C$-operad $O$, we define the category $\,O$-$\rm
alg\,$ of {\rm $O$-algebras} to be the category of actions of $O$ in
$\svf(C)$ riding the tautologous action of $\sig(C)$ on
$\svf(C)$. \end{defn}

To get a feel for these definitions, let us see how in a special
case they reduce to the definitions of untyped operads and their algebras.

\begin{ex}\label{1-ops} \et Untyped operads as $C$-operads with $C = 1$.  
{\rm Here we take $C$ to be the terminal category $1$, the category with one
object $x$ and one morphism $1_x$.  We denote the objects of $\fam(C)$ 
as $1, x, x^2 ,\dots$.  Note that $\hom(x^j,x^k)$ is the empty set
unless $j = k$, in which case it  is the symmetric group $S_k$.   An
object $A$ of $\svf(\fam(C)^\op)$ assigns to each object $x^k$ of
$\fam(C)$ a set $A_k$ equipped with an $S_k$-action.     It is
illuminating to write $A$ as a formal power series: 
\[     A = A_0 + A_1 x + A_2 x^2  + \cdots. \]
Then the coproduct in $\svf(\fam(C)^\op)$ corresponds to addition of
formal power series, where we add coefficients by taking their disjoint
union.  Similarly, the monoidal structure corresponds to multiplication of
formal power series, but where we multiply the coefficients as follows:
to multiply a set with $S_j$-action and a set with $S_k$-action, 
we take the Cartesian product with its natural $S_j \times
S_k$-action and then induce an action of $S_{j+k}$ along the 
obvious inclusion $S_j \times S_k \hookrightarrow S_{j+k}$.   
(Here by `inducing an action' we mean the left adjoint
of restricting an action to a subgroup.)

An endomorphism $P$ of $\svf(\fam(C)^\op)$ preserving small colimits
and the symmetric monoidal structure will be determined by its action on
the generating object $x$ (by which we mean the formal power series
$A$ with $A_1$ being a one-element set and $A_i$ empty for $i\ne 1$).
We have  
\[    P(x) =  P_0 + P_1 x + P_2 x^2  + \cdots,  \]
so $P$ is determined by the sets $P_k$.  Note that each set $P_k$
is equipped with an action of $S_k$, by the functoriality of $P$.
Conversely, any collection of sets $P_k$ with $S_k$-actions
determines such an endomorphism $P$.  We may also think of $P$ as
the object of $\sig(C)$ assigning to each $C$-profile $(x^k,x)$ a set
$P_k$.  We call the elements of $P_k$ the `$k$-ary
operations' of $P$.

A $C$-operad is a monoid object in $\sig(C)$.  To understand what this
amounts to, we must understand the monoidal structure in $\sig(C)$. 
This corresponds to the composition of endomorphisms of
$\svf(\fam(C)^\op)$, or in other words, composition of formal power
series.  Given $C$-signatures $P$ and $Q$, their composite is given by 
\[      (P \circ Q)(x) = P(Q(x)).  \]
Thus we have 
\ban (P \circ Q)_0 &=& P_0 + P_1 Q_0 + P_2 Q_0^2 +
\cdots \\ 
(P \circ Q)_1 &=& P_1 Q_1 + P_2 (Q_0 Q_1 + Q_1 Q_0) + \cdots \\
(P \circ Q)_2 &=& P_1 Q_2 + P_2 (Q_0 Q_2 + Q_1^2 + Q_2 Q_0) + \cdots
\ean
and so on, where we add and multiply the coefficients as before.    
Note that an element of $(P \circ Q)_j$ consists of
an element of $P_k$, for arbitrary $k \ge 0$, together with a
choice of elements of the sets $Q_{i_1}, \dots, Q_{i_k}$, where
$i_1 + \cdots + i_k = j$.   

Given a $C$-operad $O$, the product 
$\mu \maps O \circ O \to O$ gives a collection of functions from
$(O \circ O)_j$ to $O_j$.  This amounts to a collection of functions
\[ O_k \times O_{i_1} \times \cdots \times O_{i_k} \to O_{i_1 +
\cdots+ i_k}. \]
We leave it to the reader to check that in this special
case $C = 1$, $O$ being a monoid object in the
category of $C$-signatures is precisely equivalent to the conditions in
the definition of an operad given at the beginning of Section
\ref{operads}.   In particular, the associativity and unit laws there
correspond to the associativity and unit laws required of a monoid
object, while the conditions involving the symmetric groups correspond
to the fact that the product $\mu \maps O \circ O \to O$ is a symmetric
monoidal natural transformation between symmetric monoidal functors. 

In this case, the tautologous action of $\sig(C)$ on
$\svf(C)$ works as follows.  The category $\svf(C)$ is just
$\Set$, so suppose we are given a $C$-signature $P$ and a set
$A$.   Then $P$ acts on $A$ to give the set
\[       P(A) = P_0 + P_1 A + P_2 A^2 + \cdots \]
If $O$ is a $C$-operad, an $O$-algebra is a set $A$ together with an
action of $O$ on $A$, that is, a function $\alpha \maps O(A) \to A$
satisfying certain conditions.  Alteratively, as in the definition of
the algebra of an untyped operad, we can think of this action as a
collection of functions $\alpha \maps O_k \to \hom(A^k,A)$.  We leave it
for the reader to check that the conditions $\alpha$ must satisfy to be
an action are just the conditions given in Section \ref{operads}.
} \end{ex} 

Our use of formal power series above appears already in the generating
function approach to combinatorics \cite{DRS} and its categorical
interpretation in terms of `species' by Joyal \cite{Joyal}.  As shown in
Figure 7, what is at work here is the analogy between ordinary
set-theoretic linear algebra and categorified linear algebra.

\vskip 0.75em
\begin{center}
{\small
\begin{tabular}{|c|c|}  \hline
commutative rig $k$       &  symmetric 2-rig $\Set$       \\    \hline
set $S$                   &  category $C$                 \\    \hline
$k\langle S \rangle $ or $\hom(S,k)$ & $\svf(C^\op)$      \\    \hline
$FS$                      &  $\fam(C)$                    \\    \hline
$k[S]$ or $k[[S]]$        &  $\svf(\fam(C)^\op)$          \\   \hline
$\End(k[S])$ or $\End(k[[S]])$ &  $\End(\svf(\fam(C)^\op)) = \sig(C)$
\\  \hline
\end{tabular}} \vskip 1em
7.   Set-theoretic linear algebra versus categorified linear algebra
\end{center}
\vskip 0.5em 

\noindent Recall that a {\it rig} is a set with two monoid structures
$+$ and $\cdot$, where $+$ is commutative and $\cdot$ distributes over
$+$.  A 2-rig, as defined earlier, is a categorified analog of a rig.
In set-theoretic linear algebra we may work over any commutative rig
$k$, while in categorified linear algebra we may work over any symmetric
2-rig.  The free commutative rig on one element is $\N$, while the free
symmetric 2-rig on one object is $\Set$.  For simplicity, in Figure 7 we
only consider categorified linear algebra over $\Set$, although other
symmetric 2-rigs are also interesting.  It is most common in
set-theoretic linear algebra to work over a field or commutative ring,
but working over $\N$ is important in combinatorics, and heightens the
analogy to categorified linear algebra over $\Set$.

Given a set $S$ one may form the free $k$-module $k\langle S\rangle$ on
$S$.  Similarly, given a category $C$ one may form the free cocomplete
category $\svf(C^\op)$ on $C$; note that a cocomplete category is
automatically a $\Set$-module in the sense of Kapranov and Voevodsky
\cite{KV}.  One may also form the free commutative monoid $FS$ on the
set $S$.  The free commutative $k$-algebra on $S$ is then $k\langle FS
\rangle$, usually denoted by $k[S]$.  Similarly, one may form the free
symmetric monoidal category $\fam(C)$ on the category $C$.  The free
symmetric 2-rig on $C$ is then $\svf(\fam(C)^\op)$.  The monoidal
category $\End(\svf(\fam(C)^\op)) = \sig(C)$ is thus a categorified
version of the monoid $\End(k[S])$.

There are some rough spots in this analogy.  In particular, while we can
pull back $k$-valued functions along any function $f \maps S \to T$,
obtaining a $k$-linear map $f^\ast \maps k\langle T \rangle \to k
\langle S \rangle$, we cannot in general push them forwards.  In
contrast to this, not only can we pull back set-valued functors along
any functor $f \maps C \to D$, obtaining a functor $f^\ast \maps \svf(D)
\to \svf(C)$, we can also push them forward using the left adjoint
$f_\ast \maps \svf(C) \to \svf(D)$.  Both $f^\ast$ and $f_\ast$ preserve
small colimits.  In short, while the free $k$-module on a set transforms
only contravariantly under functions, the free cocomplete category on a
category transforms both covariantly and contravariantly under functors.
This plays an important role in Section \ref{pullback.operads}.

There is a kind of substitute for the free $k$-module on a set that
transforms covariantly: the $k$-module $\hom(S,k)$ of functions from $S$
to $k$.  In some ways $\svf(\fam(C^\op))$ resembles $\hom(FS,k) =
k[[S]]$ more than $k[S]$, which explains the importance of formal power
series in the generating function approach to combinatorics.  Of course,
$k\langle S \rangle$ and $\hom(S,k)$ are isomorphic when $S$ is finite;
the categorified situation works more smoothly because cocomplete
categories are closed under arbitrary colimits, while $k$-modules are
only closed under finite linear combinations.

To conclude this section, let us unpack our abstract definitions of
general $C$-operads and their algebras to obtain equivalent
`nuts-and-bolts' descriptions along more traditional lines.  First we
introduce some handy notation.  Given an object $(x_1,\dots,x_k) \in
\fam(C)$ and an object $x \in C$, we write the corresponding
$C$-profile as $(x_1,\dots,x_k,x')$.  A $C$-signature $P$ assigns to
this $C$-profile a set $P(x_1,\dots,x_k,x')$ which we call the set of
{\it operations} of $P$ with profile $(x_1,\dots,x_k,x')$.  As in
Figure 8, we may visualize such an operation as a gadget with $k$
inputs of types $x_1,\dots,x_k$ and one output of type $x'$.  Given an
operation with this profile, we call $x_1,\dots,x_k$ its {\it input
types} and $x'$ its {\it output type}, and the tuple $(x_1,\dots,x_k)$
its {\it arity}.  (In the untyped case we sometimes call the integer
$k$ the arity.)

\begin{center}
\setlength{\unitlength}{0.0007500in}%
\begingroup\makeatletter\ifx\SetFigFont\undefined
\def\x#1#2#3#4#5#6#7\relax{\def\x{#1#2#3#4#5#6}}%
\expandafter\x\fmtname xxxxxx\relax \def\y{splain}%
\ifx\x\y   
\gdef\SetFigFont#1#2#3{%
  \ifnum #1<17\tiny\else \ifnum #1<20\small\else
  \ifnum #1<24\normalsize\else \ifnum #1<29\large\else
  \ifnum #1<34\Large\else \ifnum #1<41\LARGE\else
     \huge\fi\fi\fi\fi\fi\fi
  \csname #3\endcsname}%
\else
\gdef\SetFigFont#1#2#3{\begingroup
  \count@#1\relax \ifnum 25<\count@\count@25\fi
  \def\x{\endgroup\@setsize\SetFigFont{#2pt}}%
  \expandafter\x
    \csname \romannumeral\the\count@ pt\expandafter\endcsname
    \csname @\romannumeral\the\count@ pt\endcsname
  \csname #3\endcsname}%
\fi
\fi\endgroup
\begin{picture}(1897,2565)(3226,-3304)
\thicklines
\put(4201,-1861){\circle*{150}}
\put(4201,-3061){\line( 0, 1){1200}}
\put(4201,-1861){\line(-1, 3){300}}
\put(4201,-1861){\line( 1, 3){300}}
\put(4201,-1861){\line( 1, 1){900}}
\put(4201,-1861){\line(-1, 1){900}}
\put(3226,-811){\makebox(0,0)[lb]{\raisebox{0pt}[0pt][0pt]{$x_1$}}}
\put(4126,-3286){\makebox(0,0)[lb]{\raisebox{0pt}[0pt][0pt]{$x'$}}}
\put(3901,-811){\makebox(0,0)[lb]{\raisebox{0pt}[0pt][0pt]{$x_2$}}}
\put(4501,-811){\makebox(0,0)[lb]{\raisebox{0pt}[0pt][0pt]{$x_3$}}}
\put(5101,-811){\makebox(0,0)[lb]{\raisebox{0pt}[0pt][0pt]{$x_4$}}}
\put(4426,-1966){\makebox(0,0)[lb]{\raisebox{0pt}[0pt][0pt]{$f$}}}
\end{picture}
\end{center}

\medskip
\centerline{8. An operation $f$ with profile $(x_1,x_2,x_3,x_4,x')$ }
\medskip

Since the tensor product of objects in $\sig(C)$ is given by composing
endomorphisms of $\svf(\fam(C)^\op)$, we may write the monoidal
structure in $\sig(C)$ as $\circ$.  One may check that given
$C$-signatures $P$ and $Q$, and an operation $f$ of $P$ and operations
$g_1,\dots,g_k$ of $Q$ for which the arity of $f$ is the product of
the arities of the $g_i$, we obtain an operation of $P \circ Q$.
We denote this operation by $f \circ (g_1,\dots,g_k)$. The output type
of $f \circ (g_1,\dots,g_k)$ is the output type of $f$, while its arity
type is the product of the arities of $g_1,\dots,g_k$.  We may
visualize $f \circ (g_1,\dots,g_k)$ as in Figure 9.

\begin{center}
\setlength{\unitlength}{0.000500in}%
\begingroup\makeatletter\ifx\SetFigFont\undefined
\def\x#1#2#3#4#5#6#7\relax{\def\x{#1#2#3#4#5#6}}%
\expandafter\x\fmtname xxxxxx\relax \def\y{splain}%
\ifx\x\y   
\gdef\SetFigFont#1#2#3{%
  \ifnum #1<17\tiny\else \ifnum #1<20\small\else
  \ifnum #1<24\normalsize\else \ifnum #1<29\large\else
  \ifnum #1<34\Large\else \ifnum #1<41\LARGE\else
     \huge\fi\fi\fi\fi\fi\fi
  \csname #3\endcsname}%
\else
\gdef\SetFigFont#1#2#3{\begingroup
  \count@#1\relax \ifnum 25<\count@\count@25\fi
  \def\x{\endgroup\@setsize\SetFigFont{#2pt}}%
  \expandafter\x
    \csname \romannumeral\the\count@ pt\expandafter\endcsname
    \csname @\romannumeral\the\count@ pt\endcsname
  \csname #3\endcsname}%
\fi
\fi\endgroup
\begin{picture}(3405,4252)(3879,-5183)
\thicklines
\put(4201,-1861){\circle*{150}}
\put(5701,-1861){\circle*{150}}
\put(7201,-1861){\circle*{150}}
\put(5701,-3661){\circle*{150}}
\put(4201,-3061){\line( 0, 1){1200}}
\put(4201,-3061){\line( 5,-2){1500}}
\put(5701,-3661){\line( 0, 1){600}}
\put(5701,-3061){\line( 0, 1){1200}}
\put(5701,-3661){\line( 5, 2){1500}}
\put(7201,-3061){\line( 0, 1){1200}}
\put(7201,-1861){\line( 0, 1){900}}
\put(5701,-3661){\line( 0,-1){1500}}
\put(5701,-1861){\line( 1, 1){900}}
\put(5678,-1853){\line( 1, 3){300}}
\put(5723,-1853){\line(-1, 3){300}}
\put(5701,-1861){\line(-1, 1){900}}
\put(4201,-1861){\line( 1, 3){300}}
\put(4201,-1936){\line(-1, 3){322.500}}
\put(3851,-1936){\makebox(0,0)[lb]{\raisebox{0pt}[0pt][0pt]{$g_1$}}}
\put(5351,-1936){\makebox(0,0)[lb]{\raisebox{0pt}[0pt][0pt]{$g_2$}}}
\put(6851,-1936){\makebox(0,0)[lb]{\raisebox{0pt}[0pt][0pt]{$g_3$}}}
\put(5326,-3886){\makebox(0,0)[lb]{\raisebox{0pt}[0pt][0pt]{$f$}}}
\end{picture}
\end{center}

\medskip
\centerline{9. An operation $f \circ (g_1,\dots,g_k)$ of $P \circ
Q$, where $k = 3$}
\medskip

Now suppose that $O$ is a $C$-operad.  Then it is a monoid object in
$\sig(C)$, and the product $\mu \maps O \circ O \to O$ sends each
operation $f \circ (g_1,\dots,g_k)$ in $O \circ O$ to an operation in
$O$ which we denote by $f \cdot (g_1,\dots,g_k)$.  One may check that
the associativity of the product $\mu$ implies an associativity law like
that for untyped operads. Also, the unit $\iota \maps 1 \to O$ gives
$O$ an operation $\iota_f$ of profile $(x,x')$ for every morphism $f
\maps x \to x'$ in $C$.  One may also check that the unit law and
compatibility with symmetric group actions hold as in an untyped
operad.  With a little more work, one can verify:

\begin{prop}\et \label{c.op} For any small category $C$, a
$C$-operad $O$ gives:
{\rm \begin{enumerate}  
\item {\it for any $C$-profile $(x_1,\dots,x_k,x')$, a set
$O(x_1,\dots,x_k,x')$}
\item  {\it for any $f \in O(x_1,\dots,x_k,x')$ and any
$g_1 \in O(x_{11},\dots,x_{1i_1},x_1), \dots,$ \hfill \break 
$g_k \in O(x_{k1},\dots,x_{ii_k},x_k)$, an element
\[ f \circ (g_1,\dots, g_k) \in O(x_{11},\dots,x_{1i_1}, \dots \dots,
x_{k1}, \dots , x_{ii_k},x') \] }
\item {\it for each morphism $f \maps x \to x'$ in $C$, an element
$\iota(f) \in O(x,x')$}
\item {\it for any permutation $\sigma \in S_k$, a map
\ban     \sigma \maps O(x_1,\dots,x_k,x') &\to&
          O(x_{\sigma(1)}, \dots, x_{\sigma(k)},x')  \\
f &\mapsto& f\sigma \ean  }
\end{enumerate}}
\noindent such that:
\begin{alphalist}
\item  whenever both sides make sense, 
\[      f \cdot (g_1 \cdot (h_{11}, \dots, h_{1i_1}), \dots, 
g_k \cdot (h_{k1}, \dots, h_{ki_k})) = \]
\[ (f \cdot (g_1, \dots g_k)) \; \cdot \; (h_{11}, \dots, h_{1i_1},
\dots\dots , h_{k1}, \dots, h_{ki_k}) \]
\item  for any $f \in O(x_1,\dots,x_k,x')$,
\[ f =  \iota(1_{x'}) \cdot f = f \cdot (1_{x_1},\dots,1_{x_k}) \]
\item  for any $f \in O(x_1,\dots,x_k,x')$ and 
$\sigma,\sigma' \in S_k$, 
\[      f(\sigma \sigma') = (f \sigma)\sigma' \]
\item  for any $f \in O(x_1,\dots,x_k,x')$, $\sigma \in S_k$, and 
$g_1 \in O(x_{11},\dots,x_{1i_1},x_1),\dots,$ \hfill \break
$g_k \in O(x_{k1},\dots,x_{ki_k},x_k)$, 
\[     (f\sigma) \cdot (g_\sigma(1), \dots, g_\sigma(k)) = (f \cdot
(g_1, \dots, g_k))\, \rho(\sigma), \]
where $\rho \maps S_k \to S_{i_1 + \cdots + i_k}$ is the obvious
homomorphism.
\end{alphalist}
\noindent Conversely, such data determine a unique $C$-operad.
\end{prop}

We can give a similar description of the algebras of a $C$-operad $O$.
An $O$-algebra is an action $\alpha \maps O(A) \to A$, but we usually
denote it simply as $A$.  Given an $O$-algebra $A$ and an object $x
\in C$, we call $A(x)$ the set of {\it elements of type $x$} of $A$.
For any $C$-profile $(x_1,\dots,x_k,x')$, the action $\alpha$ gives a
function 
\[ O(x_1,\dots,x_k,x') \times A(x_1) \times \cdots \times A(x_k) \to
A(x') \]
which we write as
\[  (f,a_1,\dots,a_k) \mapsto f(a_1,\dots,a_k)   .\]
Alternatively, we sometimes write this as a function
\[   O(x_1,\dots,x_k,x') \to \hom(A(x_1) \times \cdots
\times A(x_k),A(x')) \]
which by abuse of language we also call $\alpha$.  One may then
verify the following:

\begin{prop} \label{c.op.algebras}\et  For any $C$-operad
$O$, an $O$-algebra $A$ gives: 
{\rm \begin{enumerate}
\item {\it for any object $x \in C$, a set $A(x)$.}
\item {\it for any $C$-profile $(x_1,\dots,x_k,x')$, a function
\[   \alpha \maps O(x_1,\dots,x_k,x') \to \hom(A(x_1) \times \cdots
\times A(x_k),A(x')) \]  }
\end{enumerate}}
\noindent such that:
\begin{alphalist}
\item whenever both sides make sense,
\[      \alpha(f \cdot (g_1, \dots, g_k)) = \alpha(f) \circ (\alpha(g_1)
\times \cdots \times \alpha(g_k)) \]
\item  for any $x \in C$, $\alpha(\iota(1_x))$ acts as the
identity on $A(x)$
\item for any $f \in O(x_1,\dots,x_k,x')$ and $\sigma \in S_k$,
\[         \alpha(f\sigma) = \alpha(f)\sigma, \]
where $\sigma \in S_k$ acts on $\hom(A(x_1) \times \cdots \times
A(x_k),A)$ on the right by permuting the factors.  
\end{alphalist}
\noindent Conversely, such data determine a unique $O$-algebra.
\end{prop}

Starting in Section \ref{opetopes} we will restrict attention to operads
whose type category has only identity morphisms.  Such a category is
said to be discrete.  Since the category $\Set$ is isomorphic to the
category having small discrete categories as objects and functors as
morphisms, we need not worry much about the difference between small
discrete categories and sets.  Thus we may easily extend the terminology
above to define $S$-profiles, $S$-signatures, $S$-operads, and so on
when $S$ is a set.  For example, we define an {\it $S$-operad} to be an
operad whose type category is the discrete category with $S$ as its set
of objects.

\subsection{Pullback operads} \label{pullback.operads}

Given a functor $F \maps C \to D$ and a $D$-operad $O$, we now
construct a certain $C$-operad, the `pullback' $F^\ast O$.  
First recall that $D$-signatures can be regarded as set-valued
functors on $\prof(D) = \fam(D)^\op \times D$, and likewise for
$C$-signatures.   Thus we may pull back $D$-signatures to 
$C$-signatures along $F$, giving a functor
\[    F^\ast \maps \sig(D) \to \sig(C).\]
The proposition below makes $F^\ast$ into a lax monoidal functor.  
As in Section \ref{monoid.objects}, for any $D$-operad $O$,
the pullback $F^\ast O$ then becomes a $C$-operad.  

\begin{prop} \label{lax.monoidal.functor} \et  For any functor $F
\maps C \to D$,  $F^\ast \maps \sig(D) \to \sig(C)$ can be given
the structure of a lax monoidal functor.  \end{prop}

Proof - Note that $F \maps C \to D$ induces a pullback functor 
\[     F^\sharp \maps \svf(D^\op) \to \svf(C^\op) , \]
preserving small colimits, and also, because $\svf(C^\op)$ is the free
cocomplete category on $C$, a functor
\[     F_\sharp \maps \svf(C^\op) \to \svf(D^\op) \]
preserving small colimits.  In fact, $F^\sharp$ is right adjoint 
to $F_\sharp$.   By the universal property of $\svf(\fam(C^\op))$ and
$\svf(\fam(D^\op))$, the functors $F^\sharp$ and $F_\sharp$ 
induce morphisms of symmetric monoidal cocomplete categories:
\[       R \maps \svf(\fam(D^\op)) \to \svf(\fam(C^\op))  \] 
and
\[      L \maps \svf(\fam(C^\op)) \to \svf(\fam(D^\op))  \] 
with the former being right adjoint to the latter.  

Now recall that the category $\sig(D)$ is equivalent, as a monoidal
category, to the category $\End(\svf(\fam(D)^\op)).$
Thus we may identify $\sig(D)$ with this latter category, which is 
strictly monoidal.  A $D$-signature $S$ is then an endomorphism
\[      S \maps \svf(\fam(D^\op)) \to \svf(\fam(D^\op)), \]
and the composite
\[  R\circ S\circ L\maps \svf(\fam(C^\op)) \to \svf(\fam(C^\op)) \]  
is a $C$-signature.  This composition process extends to a functor
from $\sig(D)$ to $\sig(C)$, which one may check is equivalent to 
$F^\ast$.  

To make $F^\ast$ into a lax monoidal functor
it thus suffices to find a natural transformation 
$\Phi_{S,T} \maps F^\ast(S) \circ F^\ast(T) \to F^\ast(S \circ T)$
making the following diagram commute for any $D$-signatures $S,T,U$:
\be
\begin{diagram} [F^\ast(S) \circ F^\ast(T) \circ F^\ast(U)]
\node{F^\ast(S) \circ F^\ast(T) \circ F^\ast(U)} 
\arrow{e,t}{\Phi_{S,T}\, \circ 1}  
\arrow{s,l}{1 \circ \Phi_{T,U}} \node{(S \circ T) \circ F^\ast(U)} 
\arrow{s,r}{F^\ast_{S \circ T,U}}  \\
\node{F^\ast(S) \circ F^\ast(T \circ U)} \arrow{e,t}{\Phi_{S,T \circ U}}
\node{F^\ast(S \circ T \circ U)} 
\end{diagram}
\label{lax.monoidal} 
\ee
together with a morphism $\phi \maps 1_{\sig(D)} \to
F^\ast(1_{\sig(C)})$ making 
the following diagrams commute for any $D$-signature $S$:
\be
\begin{diagram} [F^\ast(1) \circ F^\ast(S)]
\node{1 \circ F^\ast(S)} \arrow{e,t}{1} \arrow{s,l}{\phi \circ 1}
\node{F^\ast(S)} \arrow{s,r}{1} \\
\node{F^\ast(1) \circ F^\ast(S)} \arrow{e,t}{\Phi_{1,S}} 
\node{F^\ast(1 \circ S)} 
\label{lax.monoidal.2} \\
\end{diagram}
\ee

\be
\begin{diagram} [F^\ast(1) \circ F^\ast(S)]
\node{F^\ast(S) \circ 1} \arrow{e,t}{1} \arrow{s,l}{1\circ \phi}
\node{F^\ast(S)} \arrow{s,r}{1} \\
\node{F^\ast(S)\circ F^\ast(1)} \arrow{e,t}{\Phi_{S,1}} 
\node{F^\ast(S\circ 1)} 
\label{lax.monoidal.3} \\
\end{diagram}
\ee
Since $R$ is the right adjoint of $L$, there is a natural
transformation $\epsilon \maps L\circ R \=> 1$, the counit of the adjunction.  
Since 
\[     F^\ast(S) \circ F^\ast(T) = R\circ S\circ L\circ R\circ T\circ L \]
while
\[     F^\ast(S \circ T) = R\circ S\circ T\circ L, \]
we may use $\epsilon$ to define 
\[      \Phi_{S,T} = 1_{R\circ S} \circ \epsilon \circ 1_{T\circ L}
\maps R\circ S\circ L\circ R\circ T\circ L \=> R\circ S\circ T\circ L
.\] 
The commutativity of (\ref{lax.monoidal}) is then easy to check.
Similarly, the unit $\iota \maps 1 \=> R \circ L$ of the adjunction
gives a morphism $\phi \maps 1 \to f^\ast(1) = R \circ L$.  The
commutativity of (\ref{lax.monoidal.2}) and (\ref{lax.monoidal.3}) then
follows from the triangle identities for an adjunction, which say that
\[ 
\begin{diagram} [L\circ R \circ L]
\node{R} \arrow{e,t}{\iota \circ 1} \node{R \circ L \circ R}
\arrow{e,t}{1 \circ \epsilon} \node{R} 
\end{diagram}
\]
and
\[ 
\begin{diagram} [L\circ R \circ L]
\node{L} \arrow{e,t}{1 \circ \iota} \node{L \circ R \circ L}
\arrow{e,t}{\epsilon \circ 1} \node{L} 
\end{diagram}
\]
are identity morphisms.  
\qed

Here we note another interesting wrinkle in the analogy between
set-theoretic linear algebra and categorified linear algebra.  A
function $f \maps S \to T$ from the finite set $S$ to the finite set $T$
induces a function $f^\ast \maps \End(k[T]) \to \End(k[S])$, using the
isomorphism $\End(k[S]) \iso k\langle FS \times S\rangle$.  However, in
contrast to Proposition \ref{lax.monoidal.functor}, this is not a monoid
homomorphism.

The same thing happens in the simpler context of matrix algebras.  For
any finite set $S$, the set $k\langle S \times S\rangle$ becomes a
monoid under matrix multiplication.  Similarly, for any category $C$,
$\svf(C^\op \times C)$ becomes a monoidal category, called the category
of {\it distributors} from $C$ to $C$.  Given a function $f \maps S \to
T$ between finite sets, the pullback $f^\ast \maps k\langle T \times
T\rangle \to k\langle S \times S\rangle$ is only a monoid homomorphism
when $f$ is one-to-one.  However, for any functor $F \maps C \to D$, the
pullback $F^\ast \maps \svf(D^\op \times D) \to \svf(C^\op \times C)$ is
a lax monoidal functor.  In fact, this follows from Proposition
\ref{lax.monoidal.functor}, using the fact that a distributor may be
regarded as a signature with only unary operations.

\subsection{The slice operad of an algebra} \label{slice.operads}

Given an $O$-algebra $A$, the slice operad $A^+$ is an operad whose
algebras are $O$-algebras over $A$, that is, equipped with a
homomorphism to $A$.   We give an explicit construction of the slice
operad and then prove it has this property.   

Recall that given a category $C$ and an object $A \in \svf(C)$, the
category $\elt(A)$ of {\rm elements} of $A$ has pairs $(x,y)$ with $x
\in C$ and $y \in A(x)$ as objects, and morphisms $f \maps
x \to x'$ with $A(f)(y) = y'$ as morphisms from $(x,y)$ to $(x',y')$.
Composition of morphisms is defined in the obvious manner.  
In this situation there is a functor $p \maps \elt(A) \to C$ with
$p(x,y) = x$ and $p(f) = f$.    

Now suppose that $O$ is a $C$-operad and $A$ is an $O$-algebra.  Then
$A$ is an object of $\svf(C)$, so as in the previous section we may form
the pullback  $p^\ast O$, which is an $\elt(A)$-operad.  Thus the
following makes sense: 

\begin{defn} \et For a $C$-operad $O$ and an $O$-algebra $A$, the {\rm
slice operad of $A$}, written $A^+$, is the sub-operad of $p^\ast O$ for
which an operation $g$ of $p^\ast O$ of profile $(a_1,\dots,a_k,a')$ is
included if and only if it satisfies $g(a_1,\dots,a_k) = a'$.  \end{defn}

\begin{prop} \label{slice.operad.algebras} \et Suppose $O$ is a
$C$-operad and $A$ is an $O$-algebra.   Then $A^+$-$\rm alg$ is 
equivalent to the category of $O$-algebras over $A$.  That is, an
$A^+$-algebra is an $O$-algebra $B$ equipped with an $O$-algebra
homomorphism $f_B \maps B \to A$, and a morphism between 
$A^+$-algebras is an $O$-algebra morphism $g \maps B \to B'$ for which
the following diagram commutes:

\[  
\begin{diagram}
\node{B} \arrow[2]{e,t}{g} \arrow{se,b}{f_B} 
\node[2]{B'} \arrow{sw,r}{f_{B'}}  \\
\node[2]{A}  
\end{diagram}
\]
\end{prop}

Proof - One may check this explicitly.
Alternatively, since the operations of $A^+$ are certain operations
of $O$, we obtain a forgetful functor from $O$-$\rm alg$ to $A^+$-$\rm
alg$.  This has a left adjoint $L \maps A^+$-${\rm alg} \to O$-$\rm alg$
sending the terminal object of $A^+$-$\rm alg$ to $A \in O$-$\rm alg$.  This
gives a functor from $A^+$-$\rm alg$ to the category of $O$-algebras equipped
with a homomorphism to $A$, which one may check is an equivalence.  \qed

\section{Opetopes and Opetopic Sets} \label{opetopes}

We now begin to address the crucial issue of categorification: the
process whereby, in passing from an $n$-categorical context to an
$(n+1)$-categorical context, laws are promoted to operations and these
new operations satisfy new laws of their own.  Our approach to this
issue relies heavily on operads.

In all that follows, we restrict attention to operads having a 
set of types, in the manner explained at the end
of Section \ref{typed.operads}.  Note that any such operad is determined by: 
\begin{enumerate}
\item its types
\item its operations
\item its reduction laws
\end{enumerate}
where by `reduction laws' we mean all equations stating that a given
composite of operations, possibly with their arguments permuted,
equals a given operation.  (Here we include unary and nullary
composites.)  Our approach to categorification relies on
a construction that yields for any operad $O$ a new operad $O^+$
having operations corresponding to the reduction laws of $O$. This
construction works roughly as follows.  In Section 
\ref{operad.operad}, we show that $S$-operads are themselves the
algebras of a certain operad.  This allows us to apply the slice
operad construction to $S$-operads, obtaining for each $S$-operad $O$
a new operad $O^+$ whose algebras are $S$-operads over $O$.  It turns
out that:
\begin{enumerate} 
\item The types of $O^+$ are the operations of $O$.  
\item The operations of $O^+$ are the reduction laws of $O$.   
\item The reduction laws of $O^+$ are the ways of combining reduction laws
of $O$ to give other reduction laws of $O$.   
\end{enumerate} 

We give numerous examples of this construction in Section
\ref{slice.operad.of.operad}.  In Section \ref{metatree} we introduce
the $n$-dimensional $O$-opetopes, which are the operations in the $n$th
iterated slice operad $O^{n+}$, and we describe a notation for them 
involving lists of labelled trees which we call `metatrees'.  We pay
special attention to the $I$-opetopes, or simply `opetopes',
because they serve as the basic shapes for cells in our approach 
to $n$-category theory.   In Section \ref{algebras} we
give a description of the algebras of $O^+$ for any $S$-operad $O$.
Finally, in Section \ref{opetopic}, we describe `opetopic sets'; 
$n$-categories are opetopic sets with certain properties.  

\subsection{The operad for operads} \label{operad.operad}

Given a small category $C$, we denote by $|C|$ the set of objects of
$C$.  We now show that for any set $S$, $S$-operads are the algebras
of a certain $|\prof(S)|$-operad. More precisely, recall from Section
\ref{typed.operads} that the category of $S$-operads, $\op(S)$, is the
category of monoid objects in $\sig(S)$.  Then we have:

\begin{thm} \label{op.c-op} \et For any set $S$, there is a
$|\prof(S)|$-operad whose category of algebras is equivalent to
$\op(S)$.  \end{thm}

Proof - We construct a $|\prof(S)|$-operad $X$ whose
category of algebras is equivalent to $\op(S)$.   The basic idea is that the
operations of $X$ are the ways of composing operations in
$S$-operads, while possibly permuting their arguments.

Note that any $S$-operad has an underlying $S$-signature, giving
us a functor
\[            R \maps \op(S) \to \sig(S) .\]
This functor has a left adjoint
\[            L \maps \sig(S) \to \op(S)  \]
assigning to each $S$-signature the free operad on that $S$-signature. 
Let $T_S$ denote the terminal $S$-operad,
and let $F$ equal $L(R(T_S))$, the free $S$-operad on the underlying
$S$-signature of $T_S$.  Note that the terminal $S$-operad has one
operation for each $S$-profile, so we may identify $S$-profiles with
operations of $T_S$.  We may think of $F$ as the $S$-operad freely generated
by all these operations.  

The operations of $F$ are in one-to-one correspondence with certain
labelled trees called {\it $T_S$-trees}.  A typical $T_S$-tree is shown
in Figure 10.  An $T_S$-tree is, first of all, a {\it combed tree}; it is
planar except at the very top, where we allow an arbitrary permutation
of the edges.  Second, each node is labelled with an operation of $T_S$, or
in other words, an $S$-profile.  A node labelled by the $S$-profile
$(x_1,\dots,x_k,x')$ must have $k$ edges coming into it from above.
Moreover, we require that it be possible to label every edge with an
element of $S$ in such a way that for any node labelled by the
$S$-profile $(x_1,\dots,x_k,x')$, the edges coming into that node from
above are labelled by the elements $x_1,\dots,x_k$ in that order from
left right, while the edge coming out of it from below is labelled by
the element $x'$.

\medskip

\begin{center}

\setlength{\unitlength}{0.000500in}%
\begingroup\makeatletter\ifx\SetFigFont\undefined
\def\x#1#2#3#4#5#6#7\relax{\def\x{#1#2#3#4#5#6}}%
\expandafter\x\fmtname xxxxxx\relax \def\y{splain}%
\ifx\x\y   
\gdef\SetFigFont#1#2#3{%
  \ifnum #1<17\tiny\else \ifnum #1<20\small\else
  \ifnum #1<24\normalsize\else \ifnum #1<29\large\else
  \ifnum #1<34\Large\else \ifnum #1<41\LARGE\else
     \huge\fi\fi\fi\fi\fi\fi
  \csname #3\endcsname}%
\else
\gdef\SetFigFont#1#2#3{\begingroup
  \count@#1\relax \ifnum 25<\count@\count@25\fi
  \def\x{\endgroup\@setsize\SetFigFont{#2pt}}%
  \expandafter\x
    \csname \romannumeral\the\count@ pt\expandafter\endcsname
    \csname @\romannumeral\the\count@ pt\endcsname
  \csname #3\endcsname}%
\fi
\fi\endgroup
\begin{picture}(5144,5444)(3279,-5183)
\thicklines

\put(5701,-3661){\circle*{150}}
\put(4201,-3061){\circle*{150}}
\put(5701,-3061){\circle*{150}}
\put(7201,-3061){\circle*{150}}
\put(7801,-1561){\circle*{150}}
\put(3601,-1561){\circle*{150}}
\put(4801,-1561){\circle*{150}}
\put(4201,-1561){\circle*{150}}
\put(6601,-1561){\circle*{150}}

\put(5701,-3661){\line( 0,-1){1500}}
\put(5701,-3661){\line(-5, 2){1500}}
\put(5701,-3661){\line( 5, 2){1500}}
\put(4201,-3061){\line( 0, 1){1500}}
\put(4201,-3061){\line(-2, 5){600}}
\put(4201,-3061){\line( 2, 5){600}}
\put(7201,-3061){\line( 2, 5){600}}
\put(3601,-1561){\line( 1, 3){300}}
\put(3601,-1561){\line(-1, 3){300}}
\put(4201,-1561){\line( 0, 1){900}}
\put(7201,-3061){\line(-2, 5){951.724}}
\put(5701,-3661){\line( 0, 1){600}}
\put(5701,-3061){\makebox(8.3333,12.5000){\SetFigFont{10}{12}{rm}.}}
\put(5701,-3061){\line( 0, 1){2400}}
\put(3901,-661){\line( 1, 3){300}}
\put(3301,-661){\line( 1, 3){300}}
\multiput(4201,-661)(-7.50000,7.50000){21}{\makebox(8.3333,12.5000){\SetFigFont{7}{8.4}{rm}.}}
\put(3901,-361){\line(-1, 1){300}}
\multiput(3301,239)(7.50000,-7.50000){21}{\makebox(8.3333,12.5000){\SetFigFont{7}{8.4}{rm}.}}
\put(5701,-661){\line( 0, 1){900}}
\put(6251,-661){\line( 4, 3){1200}}

\put(7801,-1561){\line(-1, 3){300}}
\put(7801,-1561){\line( 1, 3){300}}
\put(8101,-661){\line( 1, 3){300}}
\put(7501,-661){\line(-4, 3){516}}
\put(6301,239){\line( 4,-3){516}} 

\put(3256,-1681){\makebox(0,0)[lb]{\raisebox{0pt}[0pt][0pt]{$p_2$}}}
\put(5006,-1666){\makebox(0,0)[lb]{\raisebox{0pt}[0pt][0pt]{$p_4$}}}
\put(3886,-3241){\makebox(0,0)[lb]{\raisebox{0pt}[0pt][0pt]{$p_5$}}}
\put(5281,-3271){\makebox(0,0)[lb]{\raisebox{0pt}[0pt][0pt]{$p_1$}}}
\put(7366,-3286){\makebox(0,0)[lb]{\raisebox{0pt}[0pt][0pt]{$p_3$}}}
\put(5371,-3946){\makebox(0,0)[lb]{\raisebox{0pt}[0pt][0pt]{$p_6$}}}
\put(7966,-1726){\makebox(0,0)[lb]{\raisebox{0pt}[0pt][0pt]{$p_7$}}}
\put(3826,-1696){\makebox(0,0)[lb]{\raisebox{0pt}[0pt][0pt]{$p_8$}}}
\put(6251,-1711){\makebox(0,0)[lb]{\raisebox{0pt}[0pt][0pt]{$p_9$}}}

\end{picture}
\end{center}

\medskip
\centerline{10. A $T_S$-tree}
\medskip

\noindent In this graphical notation, we compose operations in $F$ by
combining trees essentially as in  Section \ref{untyped.operads}, and then
`combing' the resulting tree so that all the permutations of edges
occur at the very top.

Now let us turn to the $|\prof(S)|$-operad $X$.  The operations of $X$
are given as follows.  Suppose that $p_1,\dots, p_k$ are $S$-profiles.
Then $X$ has one operation $f$ of arity $(p_1,\dots,p_k)$ for each
operation $\overline f$ of $F$ that can be written as a composite of the
operations $p_1, \dots ,p_k$.  Given such an operation $f$ of $X$, we
define its output type to be the profile of $\overline f$.

Alternatively, we may describe the operations of $X$ using $T_S$-trees.
The operad $X$ has one operation of arity $(p_1,\dots,p_k)$ for
each $T_S$-tree with nodes labelled by the $S$-profiles $p_1,\dots,p_k$,
each $p_i$ labelling exactly one node.  This description makes it a bit
easier to visualize how each operation of $X$ is a way of composing
operations in an $S$-operad. For example, let $f$ be the operation of
$X$ of arity $(p_1,\dots,p_9)$ corresponding to the $T_S$-tree in
Figure 10.  Suppose that $O$ is any $S$-operad having operations $o_i$
with profiles $p_i$.  Then we can compose the $o_i$ and permute their
arguments, following the pattern given by the $T_S$-tree, to obtain the
operation
\[  (o_6 \cdot (o_5 \cdot (o_2,o_8,o_4), o_1, o_3 \cdot (o_9,o_7)))\sigma \]
where $\sigma$ is the permutation at the top of the $T_S$-tree, namely
\[    (1,2,3,4,5,6,7) \mapsto (3,1,2,4,6,5,7) .\]

In general, suppose that $f \in  X(p_1,\dots,p_k,p')$ and $O$ is a
$C$-operad.  Given operations $o_i$ of $O$ of type $p_i$,  we may
compose them and permute their arguments
in the manner described by the $T_S$-tree for $f$ to obtain an
operation of type $p'$, which we denote by $\alpha(f)(o_1,\dots,o_k)$.  
Thus we obtain a map 
\[    \alpha  \maps X(p_1,\dots,p_k,p') \to 
\hom(O(p_1) \times \cdots \times O(p_k),O(p'))   \]
where $O(p)$ denotes the set of operations of $O$ of type $p$.

Composition of operations of $X$ is defined as follows. Suppose $X$
has operations $f$ and $g_1, \dots, g_k$ of profiles for which the
composite $f \cdot (g_1, \dots, g_k)$ should be well-defined.  Let
$\overline g_i$ denote the operations of $F$ corresponding to the
operations $g_i$.  Then we define $f \cdot (g_1,\dots,g_k)$ by
\[     \overline{f \cdot (g_1,\dots, g_k)} = \alpha(f)(\overline
g_1,\dots, \overline g_k) .\]

We finish giving $X$ the structure of a $|\prof(S)|$-operad with the
help of Proposition \ref{c.op}.  First, the only morphisms in
$|\prof(S)|$ are identity morphisms, so for any $S$-profile $p$ we need
an operation $\iota(1_p) \in X(p,p)$.  We take this to be the unique
operation with that profile corresponding to the operation $p$ of $F$.
Second, for any operation $f \in X(p_1,\dots,p_k,p')$ and $\sigma \in
S_k$ we need an operation $f\sigma \in X(p_{\sigma(1)}, \dots,
p_{\sigma(k)},p')$.  We define $f\sigma$ to be the unique operation of
arity $(p_{\sigma(1)}, \dots, p_{\sigma(k)})$ corresponding to the
operation $\overline f$ of $F$.  One may then check that $X$ is a
$|\prof(S)|$-operad by verifying conditions a) - d) of Proposition
\ref{c.op}; we leave this to the reader.

Any $S$-operad $O$ becomes an $X$-algebra with the
help of Proposition \ref{c.op.algebras}.  We have already defined
the sets $O(p)$ for any $S$-profile $p$ and the action
\[    \alpha  \maps X(p_1,\dots,p_k,p') \to 
\hom(O(p_1) \times \cdots \times O(p_k),O(p'))  , \]
so one must only verify conditions a) - c).  We leave this to the 
reader as well.  Finally, it is straightforward to check that any
$X$-algebra is naturally a $S$-operad, and that a homomorphism of
$X$-algebras is the same as a homomorphism of $S$-operads.
\qed

In fact, there is also a $\prof(C)$-operad for $C$-operads for any
small category $C$.  This played an important role in an earlier
version of our approach \cite{BD2}, but for various reasons we now
prefer in what follows to work only with operads having a set, rather
than a category, of types.

\subsection{The slice operad of an operad}\label{slice.operad.of.operad}

\begin{defn} \label{slice.operad.def} \et  Given a $S$-operad $O$, let
the {\rm slice operad of $O$}, denoted $O^+$, be the $\elt(O)$-operad
whose algebras are $S$-operads over $O$, i.e., equipped with a
$C$-operad homomorphism to $O$.  \end{defn}

\noindent The existence of $O^+$ is guaranteed by Proposition
\ref{slice.operad.algebras} and Theorem \ref{op.c-op}.  The point is
that since $S$-operads are the algebras of a certain operad, we can
apply the slice operad construction to $S$-operads.   

Since $O^+$ is an $\elt(O)$-operad, it follows that the types of $O^+$
are the operations of $O$.  Also, by examining the proof of Theorem
\ref{op.c-op} one may check that the operations of $O^+$ are the
reduction laws of $O$, and the reduction laws of $O^+$ are the ways of
combining reduction laws of $O$ to obtain new reduction laws.  This
will become clearer in the next section.

To get a feel for this important construction, let us consider
some examples:

\begin{ex} \label{init.op} \et The initial untyped operad $I$ as the
operad for sets.  {\rm  Since $S$-operads form a category we may speak
of initial and terminal $S$-operads.   In the case $S = 1$, the initial
$S$-operad $I$ is the untyped operad whose only operation is the
identity.  In other words, $I$ is the untyped operad with only one unary
operation and no operations of higher arity.  Its algebras are simply
sets, so we say that $I$ is the operad for sets.
}\end{ex}

\begin{ex} \label{op.monoids} \et  $I^+$ as the operad for monoids.  
{\rm Note that $I^+$ is an $\elt(I)$-operad, but $\elt(I) = 1$, so
$I^+$ is an untyped operad.   By definition, it is the operad for
untyped operads over $I$.   An untyped operad admits a homomorphism to
$I$ only if all its operations are unary, in which case it has a
unique homomorphism to $I$.  An operad with only unary operations is
just a monoid, so $I^+$ is the operad for monoids.  The operad $I^+$
has $k!$ operations of arity $k$, corresponding to all the
elements of $S_k$, or in other words, the different
orderings in which one can multiply $k$

elements of a monoid.  The
symmetric group $S_k$ acts on these operations in an obvious way. 
 }\end{ex}

In the next example we consider an iterated slice operad.  Note that
in Definition \ref{slice.operad.def} above, $\elt(O)$ is a small
discrete category, or in other words just the {\it set} of operations
of $O$, since in applying Theorem \ref{slice.operad.algebras} we are
treating $O$ as a set-valued functor on the discrete category
$|\prof(S)|$.  Thus if $O$ is an operad with a set of types, so is $O^+$,
so we may iterate the slice operad construction. 

\begin{ex} \label{op.planar.ops} \et  $I^{++}$ as the operad for 
planar untyped operads. 
{\rm  By definition, $I^{++}$ is the $\elt(I^+)$-operad for
untyped operads over $I^+$.  An untyped operad $O$ admits a homomorphism
to $I^+$ if and only if $S_k$ acts freely on the set $O_k$ of
$k$-ary operations of $O$.   A homomorphism $f \maps O \to I^+$ is
then determined by the sets $P_k = f^{-1}(g) \subseteq
O_k$, where $g$ is the $k$-ary operation of $I^+$ corresponding to the
identity element of $S_k$ as in Example \ref{op.monoids}.  One can
check that the sets $P_k$ equipped with the composition operation of
$O$ form a planar untyped operad $P$, and conversely, any planar
untyped operad comes from an untyped operad over $I^+$ in this manner,
unique up to isomorphism.   Thus $I^{++}$ is the operad for planar
untyped operads.   
}\end{ex}

\begin{ex}\label{op.comm.monoids} \et The terminal untyped operad $T$ as the
operad for commutative monoids. 
{\rm   In the case $S = 1$, the terminal $S$-operad $T$ has one operation of 
each arity.  An algebra $A$ of $T$ is thus a commutative monoid, with
the unique $k$-ary operation of $T$ acting as the map
\ban          A^k &\to& A    \\
       (a_1, \dots, a_k) &\mapsto& a_1 \cdots a_k. \ean  
}\end{ex}

\begin{ex} \label{op.ops} \et  $T^+$ as the operad for 
untyped operads.  
{\rm  A $T^+$-algebra is an untyped operad over $T$.  Since $T$ is
terminal, a $T^+$-algebra is just an untyped operad, so $T^+$ is the
$\elt(T)$-operad for operads.   

More generally, for any set $S$ there is a terminal $S$-operad $T_S$,
having one operation of each profile.  Alternatively,
$T_S$ is the pullback of the operad $T$ along the
unique functor from $S$ to the terminal category $1$.
The slice operad $T_S^+$ is the operad for
$S$-operads.  In fact, $\elt(T_S)$ is isomorphic to $|\prof(S)|$, and
$T_S^+$ is the $|\prof(S)|$-operad for $S$-operads constructed in Theorem
\ref{op.c-op}.
}\end{ex}

At this point a comment is in order about why we base our approach on
operads rather than planar operads.  To bootstrap our way up to the
definition of $n$-categories, we want a simple sort of algebraic theory
that is powerful enough for theories of this sort to be themselves
models of a theory of this sort.  Theorem \ref{op.c-op} says that
operads have this property.  Planar operads are simpler than operads,
but planar operads are not sufficiently powerful: there
is, for example, no planar operad for planar untyped operads.  

More precisely, for any small category $C$ we define a {\it planar
$C$-operad} to be a monoid object in the category of endomorphisms of
the free 2-rig on $C$.  Taking $C = 1$ we recover the usual definition
of planar untyped operad.  Following Example \ref{op.planar.ops}, one
may check that a planar $C$-operad is the same as $C$-operad $O$
equipped with a `planar structure': a morphism to $I^+$ riding the
obvious monoidal functor from $\sig(C)$ to $\sig(1)$.  To  give the
operad for planar untyped operads a planar structure, one would need
such a morphism from $I^{++}$ to $I^+$.  One may check that no such
morphism exists. 

\subsection{Opetopes} \label{metatree}

Opetopes arise when we iterate the slice operad constuction:

\begin{defn} \et Given an $S$-operad $O$, we define $O^{0+}$ to be $O$, and
define $O^{(n+1)+} = (O^{n+})^+$ for $n \ge 1$.  \end{defn}

\begin{defn} \et Given an $S$-operad $O$, we define an {\rm $n$-dimensional
$O$-opetope} to be a type of $O^{n+}$.  We define an {\rm
$n$-dimensional opetope} to be a type of $I^{n+}$, where $I$ is
the initial untyped operad.
\end{defn}

Recall that in Theorem \ref{op.c-op} we constructed an operad for
$S$-operads, and in Example \ref{op.ops} we saw that this was just
$T_S^+$, the slice operad of the terminal $S$-operad.  The proof of
Theorem \ref{op.c-op} thus amounts to a description of the operations of
$T_S^+$ in terms of `$T_S$-trees': trees with nodes labelled by
$S$-profiles in a consistent way.  A $T_S$-tree is not quite enough to
specify a unique operation of $T_S^+$.  Rather, for any ordering
$p_1,\dots,p_k$ of the $S$-profiles labelling the nodes of an $T_S$-tree,
there is a unique operation of $T_S^+$ of arity $(p_1,\dots,p_k)$
corresponding to that $T_S$-tree.  We can keep track of this ordering by
labelling the nodes of the $T_S$-tree with additional distinct symbols
$A,B,C,\dots$, and drawing a second tree with one node having $n$ edges
coming into it from above labelled by these symbols in the desired
order.  This second tree must be planar; also, we use each symbol
exactly once as a label on this second tree.  An example is shown in
Figure 11.  Note that we use arbitrary symbols $A,B,C,\dots$ rather than the
$S$-profiles themselves to label the second tree, because the
$S$-profiles might not be distinct.

\vbox{
\begin{center}

\setlength{\unitlength}{0.000500in}%
\begingroup\makeatletter\ifx\SetFigFont\undefined
\def\x#1#2#3#4#5#6#7\relax{\def\x{#1#2#3#4#5#6}}%
\expandafter\x\fmtname xxxxxx\relax \def\y{splain}%
\ifx\x\y   
\gdef\SetFigFont#1#2#3{%
  \ifnum #1<17\tiny\else \ifnum #1<20\small\else
  \ifnum #1<24\normalsize\else \ifnum #1<29\large\else
  \ifnum #1<34\Large\else \ifnum #1<41\LARGE\else
     \huge\fi\fi\fi\fi\fi\fi
  \csname #3\endcsname}%
\else
\gdef\SetFigFont#1#2#3{\begingroup
  \count@#1\relax \ifnum 25<\count@\count@25\fi
  \def\x{\endgroup\@setsize\SetFigFont{#2pt}}%
  \expandafter\x
    \csname \romannumeral\the\count@ pt\expandafter\endcsname
    \csname @\romannumeral\the\count@ pt\endcsname
  \csname #3\endcsname}%
\fi
\fi\endgroup
\begin{picture}(8744,5444)(279,-5183)
\thicklines
\put(2401,-3661){\circle*{150}}
\put(3301,-2761){\circle*{150}}
\put(1501,-2761){\circle*{150}}
\put(901,-1561){\circle*{150}}
\put(7801,-3661){\circle*{150}}
\put(2101,-1561){\circle*{150}}
\put(2401,-3661){\line( 0,-1){1500}}
\put(2401,-3661){\line(-1, 1){900}}
\put(2401,-3661){\line( 1, 1){900}}
\put(1501,-2761){\line(-1, 2){600}}
\put(1501,-2761){\line( 1, 2){600}}
\put(901,-1561){\line(-1, 3){300}}
\put(3301,-2761){\line( 1, 2){600}}
\put(3301,-2761){\line( 0, 1){  0}}
\put(3301,-2761){\line(-1, 2){600}}
\put(601,-661){\line(-1, 3){300}}
\put(2701,-1561){\line( 1, 1){1800}}
\put(3901,-1591){\line(-1, 2){348}}
\put(3001,239){\line( 1,-2){450}}
\put(7801,-5161){\line( 0, 1){1500}}
\put(7801,-3661){\line(-1, 3){600}}
\put(7801,-3661){\line( 1, 3){600}}
\put(901,-1561){\line( 1, 2){900}}
\put(2101,-1561){\line(-2, 3){600}}
\put(901,239){\line( 2,-3){450}}
\put(7801,-1861){\line( 0,-1){1800}}
\put(7801,-3661){\line( 2, 3){1200}}
\put(7801,-3661){\line(-2, 3){1200}}
\put(601,-1786){\makebox(0,0)[lb]{\raisebox{0pt}[0pt][0pt]{$p_1$}}}
\put(1201,-3061){\makebox(0,0)[lb]{\raisebox{0pt}[0pt][0pt]{$p_4$}}}
\put(2026,-3811){\makebox(0,0)[lb]{\raisebox{0pt}[0pt][0pt]{$p_3$}}}
\put(3451,-2986){\makebox(0,0)[lb]{\raisebox{0pt}[0pt][0pt]{$B$}}}
\put(6451,-1786){\makebox(0,0)[lb]{\raisebox{0pt}[0pt][0pt]{$A$}}}
\put(7126,-1786){\makebox(0,0)[lb]{\raisebox{0pt}[0pt][0pt]{$B$}}}
\put(7726,-1786){\makebox(0,0)[lb]{\raisebox{0pt}[0pt][0pt]{$C$}}}
\put(8326,-1786){\makebox(0,0)[lb]{\raisebox{0pt}[0pt][0pt]{$D$}}}
\put(8926,-1786){\makebox(0,0)[lb]{\raisebox{0pt}[0pt][0pt]{$E$}}}
\put(1726,-1711){\makebox(0,0)[lb]{\raisebox{0pt}[0pt][0pt]{$p_5$}}}
\put(2256,-1636){\makebox(0,0)[lb]{\raisebox{0pt}[0pt][0pt]{$E$}}}
\put(1726,-2836){\makebox(0,0)[lb]{\raisebox{0pt}[0pt][0pt]{$D$}}}
\put(1126,-1636){\makebox(0,0)[lb]{\raisebox{0pt}[0pt][0pt]{$A$}}}
\put(2851,-2836){\makebox(0,0)[lb]{\raisebox{0pt}[0pt][0pt]{$p_2$}}}
\put(2626,-3736){\makebox(0,0)[lb]{\raisebox{0pt}[0pt][0pt]{$C$}}}
\end{picture}
\end{center}

\medskip
\centerline{11.  An operation of the operad for $S$-operads}
}
\medskip

It is easy to extend this notation to describe the operations of $O^+$
for any $S$-operad $O$.  Recalling the definition of slice operads
given in Section \ref{slice.operads}, it is clear that an operation of
$O^+$ can be specified as in Figure 12.  The first tree is an
arbitrary {\it $O$-tree}.  This is a combed tree with nodes labelled
by operations of $O$.  We require that a node labelled by a $k$-ary
operation have $k$ edges coming into it from above.  Moreover, we
require that it be possible to label every edge with an element of $S$
in such a way that for any node labelled by an operation with profile
$(x_1,\dots,x_k,x')$, the edges coming into that node from above are
labelled by the elements $x_1,\dots,x_k$ in that order, while the edge
coming out of it from below is labelled by the element $x'$.  As
before, we also label each node of this first tree with a distinct
symbol $A,B,C,$ etc..  Also as before, the second tree is planar and
has only one node, with $n$ edges coming into that node from above,
labelled by the same symbols $A,B,C,\dots$ in any order.  These
specify the order of the input types of the operation of $O^+$ we are
describing.

\medskip
\vbox{\begin{center}

\setlength{\unitlength}{0.000500in}%
\begingroup\makeatletter\ifx\SetFigFont\undefined
\def\x#1#2#3#4#5#6#7\relax{\def\x{#1#2#3#4#5#6}}%
\expandafter\x\fmtname xxxxxx\relax \def\y{splain}%
\ifx\x\y   
\gdef\SetFigFont#1#2#3{%
  \ifnum #1<17\tiny\else \ifnum #1<20\small\else
  \ifnum #1<24\normalsize\else \ifnum #1<29\large\else
  \ifnum #1<34\Large\else \ifnum #1<41\LARGE\else
     \huge\fi\fi\fi\fi\fi\fi
  \csname #3\endcsname}%
\else
\gdef\SetFigFont#1#2#3{\begingroup
  \count@#1\relax \ifnum 25<\count@\count@25\fi
  \def\x{\endgroup\@setsize\SetFigFont{#2pt}}%
  \expandafter\x
    \csname \romannumeral\the\count@ pt\expandafter\endcsname
    \csname @\romannumeral\the\count@ pt\endcsname
  \csname #3\endcsname}%
\fi
\fi\endgroup
\begin{picture}(8744,5144)(579,-5183)
\thicklines
\put(2401,-3661){\circle*{150}}
\put(3301,-2761){\circle*{150}}
\put(1501,-2761){\circle*{150}}
\put(901,-1561){\circle*{150}}
\put(7801,-3631){\circle*{150}}
\put(2101,-1561){\circle*{150}}
\put(2401,-2761){\circle*{150}}
\put(2401,-3661){\line( 0,-1){1500}}
\put(2401,-3661){\line(-1, 1){900}}
\put(2401,-3661){\line( 1, 1){900}}
\put(1501,-2761){\line(-1, 2){600}}
\put(1501,-2761){\line( 1, 2){600}}
\put(7801,-5161){\line( 0, 1){1500}}
\put(901,-1561){\line( 6, 5){1800}}
\put(3301,-2686){\line( 1, 4){652.941}}
\put(2101,-1561){\line(-1, 1){450}}
\put(1426,-886){\line(-1, 1){825}}
\put(2401,-3661){\line( 0, 1){900}}
\put(7801,-3661){\line(-1, 5){300}}
\put(7801,-3661){\line( 1, 5){300}}
\put(7801,-3661){\line(-3, 5){900}}
\put(7801,-3661){\line( 3, 5){900}}
\put(7801,-3661){\line(-1, 1){1500}}
\put(7801,-3661){\line( 1, 1){1500}}
\put(601,-1786){\makebox(0,0)[lb]{\raisebox{0pt}[0pt][0pt]{$f_1$}}}
\put(1201,-3061){\makebox(0,0)[lb]{\raisebox{0pt}[0pt][0pt]{$f_4$}}}
\put(2026,-3811){\makebox(0,0)[lb]{\raisebox{0pt}[0pt][0pt]{$f_3$}}}
\put(3451,-2986){\makebox(0,0)[lb]{\raisebox{0pt}[0pt][0pt]{$B$}}}
\put(1726,-1711){\makebox(0,0)[lb]{\raisebox{0pt}[0pt][0pt]{$f_5$}}}
\put(2256,-1636){\makebox(0,0)[lb]{\raisebox{0pt}[0pt][0pt]{$E$}}}
\put(1696,-2836){\makebox(0,0)[lb]{\raisebox{0pt}[0pt][0pt]{$D$}}}
\put(1126,-1636){\makebox(0,0)[lb]{\raisebox{0pt}[0pt][0pt]{$A$}}}
\put(2626,-3736){\makebox(0,0)[lb]{\raisebox{0pt}[0pt][0pt]{$C$}}}
\put(2101,-2911){\makebox(0,0)[lb]{\raisebox{0pt}[0pt][0pt]{$f_6$}}}
\put(2451,-2836){\makebox(0,0)[lb]{\raisebox{0pt}[0pt][0pt]{$F$}}}
\put(2926,-2761){\makebox(0,0)[lb]{\raisebox{0pt}[0pt][0pt]{$f_2$}}}
\put(6226,-2086){\makebox(0,0)[lb]{\raisebox{0pt}[0pt][0pt]{$A$}}}
\put(6826,-2086){\makebox(0,0)[lb]{\raisebox{0pt}[0pt][0pt]{$B$}}}
\put(7426,-2086){\makebox(0,0)[lb]{\raisebox{0pt}[0pt][0pt]{$C$}}}
\put(8026,-2086){\makebox(0,0)[lb]{\raisebox{0pt}[0pt][0pt]{$D$}}}
\put(8626,-2086){\makebox(0,0)[lb]{\raisebox{0pt}[0pt][0pt]{$E$}}}
\put(9226,-2086){\makebox(0,0)[lb]{\raisebox{0pt}[0pt][0pt]{$F$}}}
\end{picture}
\end{center}

\medskip
\centerline{12.  An operation of $O^+$}
\medskip
}

More generally, for any $n > 1$ one can specify any
$n$-dimensional $O$-opetope by means of an 
{\it $n$-dimensional $O$-metatree}, as in Figure 13.

\medskip
\vbox{
\begin{center}

\setlength{\unitlength}{0.00037500in}%
\begingroup\makeatletter\ifx\SetFigFont\undefined
\def\x#1#2#3#4#5#6#7\relax{\def\x{#1#2#3#4#5#6}}%
\expandafter\x\fmtname xxxxxx\relax \def\y{splain}%
\ifx\x\y   
\gdef\SetFigFont#1#2#3{%
  \ifnum #1<17\tiny\else \ifnum #1<20\small\else
  \ifnum #1<24\normalsize\else \ifnum #1<29\large\else
  \ifnum #1<34\Large\else \ifnum #1<41\LARGE\else
     \huge\fi\fi\fi\fi\fi\fi
  \csname #3\endcsname}%
\else
\gdef\SetFigFont#1#2#3{\begingroup
  \count@#1\relax \ifnum 25<\count@\count@25\fi
  \def\x{\endgroup\@setsize\SetFigFont{#2pt}}%
  \expandafter\x
    \csname \romannumeral\the\count@ pt\expandafter\endcsname
    \csname @\romannumeral\the\count@ pt\endcsname
  \csname #3\endcsname}%
\fi
\fi\endgroup
\begin{picture}(12944,5305)(579,-5183)
\thicklines
\put(2401,-3661){\circle*{150}}
\put(3301,-2761){\circle*{150}}
\put(1501,-2761){\circle*{150}}
\put(901,-1561){\circle*{150}}
\put(2101,-1561){\circle*{150}}
\put(2401,-2761){\circle*{150}}
\put(5101,-1861){\circle*{150}}
\put(5701,-2761){\circle*{150}}
\put(6601,-3661){\circle*{150}}
\put(10201,-3661){\circle*{150}}
\put(10801,-2461){\circle*{150}}
\put(12901,-3661){\circle*{150}}
\put(2401,-3661){\line( 0,-1){1500}}
\put(2401,-3661){\line(-1, 1){900}}
\put(2401,-3661){\line( 1, 1){900}}
\put(1501,-2761){\line(-1, 2){600}}
\put(1501,-2761){\line( 1, 2){600}}
\put(901,-1561){\line( 6, 5){1800}}
\put(3301,-2686){\line( 1, 4){652.941}}
\put(2101,-1561){\line(-1, 1){450}}
\put(1426,-886){\line(-1, 1){825}}
\put(2401,-3661){\line( 0, 1){900}}
\put(6601,-5161){\line( 0, 1){1500}}
\put(6601,-3661){\line( 1, 1){900}}
\put(6601,-3661){\line(-1, 1){900}}
\put(5701,-2761){\line(-2, 3){600}}
\put(5701,-2761){\line( 2, 3){600}}
\put(5701,-2761){\line( 0, 1){900}}
\put(5101,-1861){\line(-1, 4){300}}
\put(5101,-1861){\line( 0, 1){1200}}
\put(5101,-661){\line( 0, 1){600}}
\put(4801,-661){\line(-1, 4){150}}
\put(8101,-361){\line( 0, 1){300}}
\put(7501,-2761){\line( 1, 3){600}}
\put(8101,-961){\line( 0, 1){600}}
\put(5101,-1861){\line( 1, 3){300}}
\put(5401,-961){\line( 5, 3){1500}}
\put(5701,-1861){\line( 0, 1){975}}
\put(5701,-661){\line( 0, 1){600}}
\put(6301,-1861){\line( 0, 1){1350}}
\put(6301,-511){\line( 0,-1){ 75}}
\put(6301,-586){\line( 0, 1){ 75}}
\put(6301,-61){\line( 0,-1){225}}
\put(10201,-5161){\line( 0, 1){1500}}
\put(10201,-3661){\line( 1, 2){600}}
\put(10201,-3661){\line(-1, 2){1200}}
\put(10801,-2461){\line( 1, 2){600}}
\put(10801,-2536){\line(-1, 2){630}}
\put(12901,-5161){\line( 0, 1){1500}}
\put(12901,-3661){\line( 1, 2){600}}
\put(12901,-3661){\line(-1, 2){600}}
\put(9001,-1261){\line( 1, 1){1200}}

\put(10181,-1261){\line(-1, 1){525}}

\put(9001,-61){\line( 1,-1){525}}
\put(11401,-1261){\line( 0, 1){1200}}
\put(551,-1786){\makebox(0,0)[lb]{\raisebox{0pt}[0pt][0pt]{$f_1$}}}
\put(1151,-3061){\makebox(0,0)[lb]{\raisebox{0pt}[0pt][0pt]{$f_4$}}}
\put(1826,-3811){\makebox(0,0)[lb]{\raisebox{0pt}[0pt][0pt]{$f_3$}}}
\put(3451,-2986){\makebox(0,0)[lb]{\raisebox{0pt}[0pt][0pt]{$B$}}}
\put(1556,-1711){\makebox(0,0)[lb]{\raisebox{0pt}[0pt][0pt]{$f_5$}}}
\put(2326,-1636){\makebox(0,0)[lb]{\raisebox{0pt}[0pt][0pt]{$E$}}}
\put(1126,-1636){\makebox(0,0)[lb]{\raisebox{0pt}[0pt][0pt]{$A$}}}
\put(2626,-3736){\makebox(0,0)[lb]{\raisebox{0pt}[0pt][0pt]{$C$}}}
\put(1971,-2911){\makebox(0,0)[lb]{\raisebox{0pt}[0pt][0pt]{$f_6$}}}
\put(2816,-2761){\makebox(0,0)[lb]{\raisebox{0pt}[0pt][0pt]{$f_2$}}}
\put(2476,-2611){\makebox(0,0)[lb]{\raisebox{0pt}[0pt][0pt]{$F$}}}
\put(1651,-2761){\makebox(0,0)[lb]{\raisebox{0pt}[0pt][0pt]{$D$}}}
\put(4576, 14){\makebox(0,0)[lb]{\raisebox{0pt}[0pt][0pt]{$A$}}}
\put(5026, 14){\makebox(0,0)[lb]{\raisebox{0pt}[0pt][0pt]{$D$}}}
\put(5626, 14){\makebox(0,0)[lb]{\raisebox{0pt}[0pt][0pt]{$C$}}}
\put(6826, 14){\makebox(0,0)[lb]{\raisebox{0pt}[0pt][0pt]{$E$}}}
\put(6226, 14){\makebox(0,0)[lb]{\raisebox{0pt}[0pt][0pt]{$F$}}}
\put(8026, 14){\makebox(0,0)[lb]{\raisebox{0pt}[0pt][0pt]{$B$}}}
\put(4706,-2011){\makebox(0,0)[lb]{\raisebox{0pt}[0pt][0pt]{$G$}}}
\put(5306,-2986){\makebox(0,0)[lb]{\raisebox{0pt}[0pt][0pt]{$H$}}}
\put(6171,-3811){\makebox(0,0)[lb]{\raisebox{0pt}[0pt][0pt]{$I$}}}
\put(11026,-2611){\makebox(0,0)[lb]{\raisebox{0pt}[0pt][0pt]{$J$}}}
\put(10426,-3811){\makebox(0,0)[lb]{\raisebox{0pt}[0pt][0pt]{$K$}}}
\put(13426,-2386){\makebox(0,0)[lb]{\raisebox{0pt}[0pt][0pt]{$J$}}}
\put(12226,-2386){\makebox(0,0)[lb]{\raisebox{0pt}[0pt][0pt]{$K$}}}
\put(8926, 14){\makebox(0,0)[lb]{\raisebox{0pt}[0pt][0pt]{$G$}}}
\put(11326, 14){\makebox(0,0)[lb]{\raisebox{0pt}[0pt][0pt]{$H$}}}
\put(10126, 14){\makebox(0,0)[lb]{\raisebox{0pt}[0pt][0pt]{$I$}}}
\end{picture}
\end{center}

\medskip
\centerline{13.  A 3-dimensional $O$-metatree}
\medskip
}
This is a list of $n$ labelled trees, the last of which is a
planar tree with only one node, while the rest are combed trees.
The first tree is an arbitrary $O$-tree.  For $1 \le i < n$, every node
of the $i$th tree is labelled with a distinct symbol, and the same symbols
also label all the edges at the very top of the $(i + 1)$st tree, each
symbol labelling exactly one edge.  In addition, each edge of the
$(i + 1)$st tree must correspond to a subtree of the $i$th tree in such
a way that: 
\begin{enumerate}
\item The edge at the very top of the $(i + 1)$st tree labelled by a
given symbol corresponds to the subtree of the $i$th tree whose
one and only node is labelled by the same symbol.  
\item  The edge of the $(i + 1)$st tree coming out of a given note from
below corresponds to the subtree that is the union of the subtrees
corresponding to the edges coming into that node from above.
\item  The edge at the very bottom of the $(i+1)$st tree corresponds to
the whole $i$th tree.
\end{enumerate}
Special care must be taken when the node of the last tree has
no edges coming into it from above.  This can only occur when
all the previous trees are empty.  This sort of metatree describes
a nullary operation of $O^{(n-1)+}$ whose output type is an
identity operation $1_x$ of $O^{(n-2)+}$.  To specify which identity
operation, we need to label the edge coming out of the node of
the last tree from below with the operation $1_x$.  

We conclude this section with some examples which
begin to explain the role opetopes play in $n$-category theory.

\begin{ex} \label{I+.metatree}\et Metatree notation for
operations of $I^+$.  {\rm  Let $I$ be the initial untyped operad
as in Example \ref{init.op}.  
Since the only operation in $I$ is the unary operation 
$1$, a metatree for a typical operation of $I^+$ looks like
that in Figure 14.  As we expect from Example \ref{op.monoids},
$I^+$ has $n!$ operations of arity $n$.

\medskip
\vbox{
\begin{center}

\setlength{\unitlength}{0.000500in}%
\begingroup\makeatletter\ifx\SetFigFont\undefined
\def\x#1#2#3#4#5#6#7\relax{\def\x{#1#2#3#4#5#6}}%
\expandafter\x\fmtname xxxxxx\relax \def\y{splain}%
\ifx\x\y   
\gdef\SetFigFont#1#2#3{%
  \ifnum #1<17\tiny\else \ifnum #1<20\small\else
  \ifnum #1<24\normalsize\else \ifnum #1<29\large\else
  \ifnum #1<34\Large\else \ifnum #1<41\LARGE\else
     \huge\fi\fi\fi\fi\fi\fi
  \csname #3\endcsname}%
\else
\gdef\SetFigFont#1#2#3{\begingroup
  \count@#1\relax \ifnum 25<\count@\count@25\fi
  \def\x{\endgroup\@setsize\SetFigFont{#2pt}}%
  \expandafter\x
    \csname \romannumeral\the\count@ pt\expandafter\endcsname
    \csname @\romannumeral\the\count@ pt\endcsname
  \csname #3\endcsname}%
\fi
\fi\endgroup
\begin{picture}(4522,4544)(4201,-5183)
\thicklines
\put(7801,-3631){\circle*{150}}
\put(4501,-3361){\circle*{150}}
\put(4501,-4261){\circle*{150}} 
\put(4501,-2461){\circle*{150}} 
\put(4501,-1561){\circle*{150}} 
\put(7801,-5161){\line( 0, 1){1500}}
\put(7801,-3661){\line(-1, 5){300}}
\put(7801,-3661){\line( 1, 5){300}}
\put(7801,-3661){\line(-3, 5){900}}
\put(7801,-3661){\line( 3, 5){900}}
\put(4501,-5161){\line( 0, 1){4500}}
\put(4201,-4336){\makebox(0,0)[lb]{\raisebox{0pt}[0pt][0pt]{$1$}}}
\put(4201,-3436){\makebox(0,0)[lb]{\raisebox{0pt}[0pt][0pt]{$1$}}}
\put(4201,-2536){\makebox(0,0)[lb]{\raisebox{0pt}[0pt][0pt]{$1$}}}
\put(4201,-1636){\makebox(0,0)[lb]{\raisebox{0pt}[0pt][0pt]{$1$}}}
\put(4801,-1636){\makebox(0,0)[lb]{\raisebox{0pt}[0pt][0pt]{$A$}}}
\put(4801,-2536){\makebox(0,0)[lb]{\raisebox{0pt}[0pt][0pt]{$B$}}}
\put(4801,-3436){\makebox(0,0)[lb]{\raisebox{0pt}[0pt][0pt]{$C$}}}
\put(4801,-4336){\makebox(0,0)[lb]{\raisebox{0pt}[0pt][0pt]{$D$}}}
\put(6826,-2086){\makebox(0,0)[lb]{\raisebox{0pt}[0pt][0pt]{$A$}}}
\put(7426,-2086){\makebox(0,0)[lb]{\raisebox{0pt}[0pt][0pt]{$B$}}} 
\put(8026,-2086){\makebox(0,0)[lb]{\raisebox{0pt}[0pt][0pt]{$C$}}}
\put(8626,-2086){\makebox(0,0)[lb]{\raisebox{0pt}[0pt][0pt]{$D$}}}
\end{picture}
\end{center}

\medskip
\centerline{14.  An operation of $I^+$}
\medskip
}

The term `opetope' is explained by the fact that one can associate
to the $n$-dimensional opetopes
certain labelled $n$-dimensional combinatorial polytopes, or 
generalizations thereof.  In particular, the
operations of $I^+$ are the 2-dimensional opetopes, and 
the $k$-ary operations of $I^+$ correspond to polygons
with $k$ labelled `infaces' and one `outface'.  
For example, the $4$-ary operation in Figure 14 corresponds to the polygon
shown in Figure 15, with four labelled infaces and one outface.  

\begin{center}
\setlength{\unitlength}{0.000500in}%
\begingroup\makeatletter\ifx\SetFigFont\undefined
\def\x#1#2#3#4#5#6#7\relax{\def\x{#1#2#3#4#5#6}}%
\expandafter\x\fmtname xxxxxx\relax \def\y{splain}%
\ifx\x\y   
\gdef\SetFigFont#1#2#3{%
  \ifnum #1<17\tiny\else \ifnum #1<20\small\else
  \ifnum #1<24\normalsize\else \ifnum #1<29\large\else
  \ifnum #1<34\Large\else \ifnum #1<41\LARGE\else
     \huge\fi\fi\fi\fi\fi\fi
  \csname #3\endcsname}%
\else
\gdef\SetFigFont#1#2#3{\begingroup
  \count@#1\relax \ifnum 25<\count@\count@25\fi
  \def\x{\endgroup\@setsize\SetFigFont{#2pt}}%
  \expandafter\x
    \csname \romannumeral\the\count@ pt\expandafter\endcsname
    \csname @\romannumeral\the\count@ pt\endcsname
  \csname #3\endcsname}%
\fi
\fi\endgroup
\begin{picture}(3644,1919)(3279,-4658)
\thicklines
\put(3301,-4561){\line( 1, 2){600}}
\put(3901,-3361){\line( 2, 1){1200}}
\put(5101,-2761){\line( 2,-1){1200}}
\put(6301,-3361){\line( 1,-2){600}}
\put(3301,-4561){\line( 1, 0){3600}}

\multiput(5026,-4486)(9.37500,-4.68750){17}{\makebox(8.3333,12.5000){\SetFigFont{7}{8.4}{rm}.}}
\multiput(5026,-4636)(9.37500,4.68750){17}{\makebox(8.3333,12.5000){\SetFigFont{7}{8.4}{rm}.}}
\multiput(3526,-3886)(9.37500,4.68750){17}{\makebox(8.3333,12.5000){\SetFigFont{7}{8.4}{rm}.}}
\put(3676,-3811){\line( 0,-1){180}}
\put(4351,-3061){\line( 1, 0){150}}
\multiput(4501,-3061)(-4.68750,-9.37500){17}{\makebox(8.3333,12.5000){\SetFigFont{7}{8.4}{rm}.}}
\put(5551,-3061){\line( 1, 0){150}}
\multiput(5701,-3061)(-4.68750,9.37500){17}{\makebox(8.3333,12.5000){\SetFigFont{7}{8.4}{rm}.}}
\put(6601,-3961){\line( 0, 1){195}}
\multiput(6601,-3961)(-9.37500,4.68750){17}{\makebox(8.3333,12.5000){\SetFigFont{7}{8.4}{rm}.}}
\multiput(4951,-3886)(6.00000,-9.00000){26}{\makebox(8.3333,12.5000){\SetFigFont{7}{8.4}{rm}.}}
\multiput(5101,-4111)(6.00000,9.00000){26}{\makebox(8.3333,12.5000){\SetFigFont{7}{8.4}{rm}.}}
\put(5064,-4051){\line( 0, 1){390}}
\put(5146,-4051){\line( 0, 1){390}}
\put(3338,-3879){\makebox(0,0)[lb]{\raisebox{0pt}[0pt][0pt]{$A$}}}
\put(4306,-2978){\makebox(0,0)[lb]{\raisebox{0pt}[0pt][0pt]{$B$}}}
\put(5791,-2978){\makebox(0,0)[lb]{\raisebox{0pt}[0pt][0pt]{$C$}}}
\put(6736,-3924){\makebox(0,0)[lb]{\raisebox{0pt}[0pt][0pt]{$D$}}}
\end{picture}
\end{center}

\medskip
\centerline{15.  2-dimensional opetope represented as a polytope}
\medskip

The degenerate cases $k = 0$ and $k = 1$ are a bit of a nuisance
because one  cannot represent `unigons' and `bigons' as convex
geometrical polytopes.   Nonetheless, one can still draw them if one
allows curved edges, and these drawings are widely used in
2-categorical commutative diagrams. In fact, the bigon is the only
basic shape of $2$-cell in the traditional  globular approach to
$n$-category theory; to achieve the effect of 2-cells with other
shapes one resorts to pasting theorems \cite{Crans,Johnson,Power}.  In
the opetopic approach the the basic shapes of cells are the opetopes,
which may have any number of infaces but always exactly one outface.  
For example, we use a 2-cell shaped like the opetope in Figure 15 to
represent an operation  having the 1-cells $A$, $B$, $C$, and $D$ as
inputs and  the outface 1-cell as its output.  In particular, we use a
`universal'  2-cell of this sort --- as defined below in Section
\ref{coherent} --- to represent a process of composing the 1-cells
$A$, $B$, $C$, and $D$.   The outface is then called a `composite' of
these 1-cells.  
} \end{ex}

\begin{ex} \label{I++.metatree}\et Metatree notation for
operations of $I^{++}$.  {\rm  
A metatree for a typical operation of $I^{++}$ is shown in Figure
16.

\medskip
\vbox{
\begin{center}
\setlength{\unitlength}{0.000500in}%
\begingroup\makeatletter\ifx\SetFigFont\undefined
\def\x#1#2#3#4#5#6#7\relax{\def\x{#1#2#3#4#5#6}}%
\expandafter\x\fmtname xxxxxx\relax \def\y{splain}%
\ifx\x\y   
\gdef\SetFigFont#1#2#3{%
  \ifnum #1<17\tiny\else \ifnum #1<20\small\else
  \ifnum #1<24\normalsize\else \ifnum #1<29\large\else
  \ifnum #1<34\Large\else \ifnum #1<41\LARGE\else
     \huge\fi\fi\fi\fi\fi\fi
  \csname #3\endcsname}%
\else
\gdef\SetFigFont#1#2#3{\begingroup
  \count@#1\relax \ifnum 25<\count@\count@25\fi
  \def\x{\endgroup\@setsize\SetFigFont{#2pt}}%
  \expandafter\x
    \csname \romannumeral\the\count@ pt\expandafter\endcsname
    \csname @\romannumeral\the\count@ pt\endcsname
  \csname #3\endcsname}%
\fi
\fi\endgroup
\begin{picture}(6922,4768)(4201,-4958)
\thicklines
\put(8101,-3061){\circle*{150}}
\put(7501,-3961){\circle*{150}}
\put(10501,-3961){\circle*{150}}
\put(4501,-3961){\circle*{150}}
\put(4501,-3061){\circle*{150}}
\put(4501,-2161){\circle*{150}}
\put(7501,-3961){\line( 2, 3){600}}
\put(7501,-3961){\line(-2, 3){600}}
\put(8101,-3061){\line(-2, 3){600}}
\put(8101,-3061){\line( 2, 3){600}}
\put(6901,-3061){\line(-2, 3){600}}
\put(7501,-4036){\line( 0,-1){825}}
\put(10501,-4861){\line( 0, 1){900}}
\put(10501,-3961){\line( 2, 3){600}}
\put(10501,-3961){\line(-2, 3){600}}
\put(4501,-1261){\line( 0,-1){3675}}
\put(6301,-2161){\line( 2, 3){1200}}
\put(7501,-2161){\line(-2, 3){542.308}}
\put(6301,-361){\line( 2,-3){542.308}}
\put(8701,-2161){\line( 0, 1){1800}}
\put(4201,-3136){\makebox(0,0)[lb]{\raisebox{0pt}[0pt][0pt]{$1$}}}
\put(4201,-2236){\makebox(0,0)[lb]{\raisebox{0pt}[0pt][0pt]{$1$}}}
\put(4201,-4036){\makebox(0,0)[lb]{\raisebox{0pt}[0pt][0pt]{$1$}}}
\put(4801,-4036){\makebox(0,0)[lb]{\raisebox{0pt}[0pt][0pt]{$C$}}}
\put(4801,-3136){\makebox(0,0)[lb]{\raisebox{0pt}[0pt][0pt]{$B$}}}
\put(4801,-2236){\makebox(0,0)[lb]{\raisebox{0pt}[0pt][0pt]{$A$}}}
\put(6226,-286){\makebox(0,0)[lb]{\raisebox{0pt}[0pt][0pt]{$C$}}}
\put(7426,-286){\makebox(0,0)[lb]{\raisebox{0pt}[0pt][0pt]{$A$}}}
\put(8626,-286){\makebox(0,0)[lb]{\raisebox{0pt}[0pt][0pt]{$B$}}}
\put(8251,-3211){\makebox(0,0)[lb]{\raisebox{0pt}[0pt][0pt]{$D$}}}
\put(7106,-4111){\makebox(0,0)[lb]{\raisebox{0pt}[0pt][0pt]{$E$}}}
\put(9826,-2986){\makebox(0,0)[lb]{\raisebox{0pt}[0pt][0pt]{$E$}}}
\put(11026,-2986){\makebox(0,0)[lb]{\raisebox{0pt}[0pt][0pt]{$D$}}}

\end{picture}
\end{center}

\medskip
\centerline{16.  An operation of $I^{++}$}
\medskip
}

\noindent The operations of $I^{++}$ are the 3-dimensional opetopes,
and we can associate to them certain 3-dimensional combinatorial
polytopes or generalizations thereof.   For example, the operation of
Figure 16 corresponds to the polytope shown in Figure 17, having two
triangular `infaces' labelled $D$ and $E$ on top, and having the
quadrilateral  on the bottom as `outface'.   Note that while this is a
combinatorial polytope, it cannot be realized as a convex geometrical
polytope.  As in the 2-dimensional case, there are also `degenerate'
3-dimensional opetopes that cannot be realized as combinatorial
polytopes in the strict sense. Also note that Figure 17 does not
record all the information needed to uniquely specify an operation of
$I^{++}$, because it does not keep track of the permutations in the
metatree of Figure 16.   Because of these problems we find it better
to describe opetopes using metatrees.  Nonetheless, the polytopes may
help the reader relate our approach to other work on $n$-categories.

In the opetopic approach to $n$-categories, we use a universal 3-cell
shaped like that in Figure 17 to represent the process of composing
the 2-cells $D$ and $E$ in the indicated manner to obtain a 2-cell
shaped like the outface.   More generally, an $n$-dimensional opetope
always has some number of $(n-1)$-dimensional  opetopes as infaces,
pasted together in a manner described by a tree, together with a
single $(n-1)$-dimensional opetope as outface. A universal $n$-cell of
this shape then describes a process of composing $(n-1)$-cells shaped
like the infaces to obtain an $(n-1)$-cell shaped like the outface.  

\begin{center}
\setlength{\unitlength}{0.000500in}%
\begingroup\makeatletter\ifx\SetFigFont\undefined
\def\x#1#2#3#4#5#6#7\relax{\def\x{#1#2#3#4#5#6}}%
\expandafter\x\fmtname xxxxxx\relax \def\y{splain}%
\ifx\x\y   
\gdef\SetFigFont#1#2#3{%
  \ifnum #1<17\tiny\else \ifnum #1<20\small\else
  \ifnum #1<24\normalsize\else \ifnum #1<29\large\else
  \ifnum #1<34\Large\else \ifnum #1<41\LARGE\else
     \huge\fi\fi\fi\fi\fi\fi
  \csname #3\endcsname}%
\else
\gdef\SetFigFont#1#2#3{\begingroup
  \count@#1\relax \ifnum 25<\count@\count@25\fi
  \def\x{\endgroup\@setsize\SetFigFont{#2pt}}%
  \expandafter\x
    \csname \romannumeral\the\count@ pt\expandafter\endcsname
    \csname @\romannumeral\the\count@ pt\endcsname
  \csname #3\endcsname}%
\fi
\fi\endgroup
\begin{picture}(3644,2068)(7779,-4658)
\thicklines
\put(7801,-4561){\line( 1, 0){3600}}
\put(7801,-4561){\line( 1, 2){900}}
\put(8701,-2761){\line( 1, 0){1800}}
\put(10501,-2761){\line( 1,-2){900}}
\put(8701,-2761){\line( 3,-2){2700}}

\multiput(9526,-4486)(9.37500,-4.68750){17}{\makebox(8.3333,12.5000){\SetFigFont{7}{8.4}{rm}.}}
\multiput(9526,-4636)(9.37500,4.68750){17}{\makebox(8.3333,12.5000){\SetFigFont{7}{8.4}{rm}.}}

\multiput(9526,-2836)(9.37500,4.68750){17}{\makebox(8.3333,12.5000){\SetFigFont{7}{8.4}{rm}.}}
\multiput(9526,-2686)(9.37500,-4.68750){17}{\makebox(8.3333,12.5000){\SetFigFont{7}{8.4}{rm}.}}
\multiput(8176,-3661)(7.50000,7.50000){21}{\makebox(8.3333,12.5000){\SetFigFont{7}{8.4}{rm}.}}
\put(8326,-3736){\line( 0, 1){225}}
\put(10951,-3661){\line( 0, 1){225}}
\multiput(10951,-3661)(-7.50000,7.50000){21}{\makebox(8.3333,12.5000){\SetFigFont{7}{8.4}{rm}.}}
\put(9901,-3136){\line( 0,-1){225}}
\put(9901,-3361){\line( 1, 0){225}}
\put(9976,-3361){\line( 1, 1){300}}
\put(9901,-3286){\line( 1, 1){300}}
\put(8701,-3661){\line( 0,-1){375}}
\put(8821,-3661){\line( 0,-1){375}}
\multiput(8596,-3961)(8.14988,-6.79157){22}{\makebox(8.3333,12.5000){\SetFigFont{7}{8.4}{rm}.}}
\multiput(8761,-4111)(7.50000,7.50000){21}{\makebox(8.3333,12.5000){\SetFigFont{7}{8.4}{rm}.}}
\multiput(10216,-3661)(5.36029,-8.93383){21}{\makebox(8.3333,12.5000){\SetFigFont{7}{8.4}{rm}.}}
\multiput(10321,-3841)(-10.21672,2.55418){20}{\makebox(8.3333,12.5000){\SetFigFont{7}{8.4}{rm}.}}
\put(9301,-3886){\line( 1,-2){150}}
\put(9451,-4186){\line( 1, 2){150}}
\put(9376,-3661){\line( 0,-1){375}}
\put(9451,-3661){\line( 0,-1){525}}
\put(9526,-3661){\line( 0,-1){375}}
\put(8026,-3586){\makebox(0,0)[lb]{\raisebox{0pt}[0pt][0pt]{$A$}}}
\put(9526,-2611){\makebox(0,0)[lb]{\raisebox{0pt}[0pt][0pt]{$B$}}}
\put(11101,-3586){\makebox(0,0)[lb]{\raisebox{0pt}[0pt][0pt]{$C$}}}
\put(10276,-3361){\makebox(0,0)[lb]{\raisebox{0pt}[0pt][0pt]{$D$}}}
\put(8476,-3886){\makebox(0,0)[lb]{\raisebox{0pt}[0pt][0pt]{$E$}}}

\end{picture}
\end{center}

\medskip
\centerline{17.  A 3-dimensional opetope represented as a polytope}
\medskip

\noindent 

} \end{ex}

\subsection{Algebras of slice operads} \label{algebras}

The following examples lead up to a concrete description, for
any $S$-operad $O$, of the algebras of $O^+$.  

\begin{ex}\label{op.pointed.sets} \et The free operad 
on one nullary operation, $K$, as the operad for pointed sets.
{\rm Let $K$ be the untyped operad with one nullary operation $k$,
one unary operation $1$ (the identity operation), and no other
operations.  A $K$-algebra is simply a pointed set.  
} \end{ex}

\begin{ex} \label{op.monoid.actions} \et $K^+$ as the operad for
monoid actions.  {\rm Since $K$ has two operations, $K^+$
has two types, $k$ and $1$.   The operations of $K^+$ include the
three operations shown in metatree notation in Figure 18:
a nullary operation with output type $1$, a binary operation with
profile $(1,1,1)$, and a binary operation with profile $(k,1,k)$.   All 
the operations of $K^+$ are generated from these three 
by composition.    A $K^+$-algebra $A$ thus consists of a set
$A(1)$ and a set $A(k)$ together with a special element $i \in 
A(1)$, a map $m \maps A(1) \times A(1) \to A(1)$ and a map
$a \maps A(k) \times A(1) \to A(1)$ satisfying certain laws.  
One may check that these laws say precisely that $A(1)$ is a monoid with
a right action on $A(k)$. 
\medskip

\begin{center}
\setlength{\unitlength}{0.000500in}%
\begingroup\makeatletter\ifx\SetFigFont\undefined
\def\x#1#2#3#4#5#6#7\relax{\def\x{#1#2#3#4#5#6}}%
\expandafter\x\fmtname xxxxxx\relax \def\y{splain}%
\ifx\x\y   
\gdef\SetFigFont#1#2#3{%
  \ifnum #1<17\tiny\else \ifnum #1<20\small\else
  \ifnum #1<24\normalsize\else \ifnum #1<29\large\else
  \ifnum #1<34\Large\else \ifnum #1<41\LARGE\else
     \huge\fi\fi\fi\fi\fi\fi
  \csname #3\endcsname}%
\else
\gdef\SetFigFont#1#2#3{\begingroup
  \count@#1\relax \ifnum 25<\count@\count@25\fi
  \def\x{\endgroup\@setsize\SetFigFont{#2pt}}%
  \expandafter\x
    \csname \romannumeral\the\count@ pt\expandafter\endcsname
    \csname @\romannumeral\the\count@ pt\endcsname
  \csname #3\endcsname}%
\fi
\fi\endgroup
\begin{picture}(3622,5512)(3901,-5783)
\thicklines
\put(7201,-361){\circle*{150}}
\put(7201,-2761){\circle*{150}}
\put(7201,-5161){\circle*{150}}
\put(4201,-5161){\circle*{150}}
\put(4201,-4561){\circle*{150}}
\put(4201,-2761){\circle*{150}}
\put(4201,-2161){\circle*{150}}
\put(6901,-2161){\line( 1,-2){300}}
\put(7201,-2761){\line( 1, 2){300}}
\put(7201,-2761){\line( 0,-1){600}}
\put(7201,-361){\line( 0,-1){600}}
\put(7101,-1301){\makebox(0,0)[lb]{\raisebox{0pt}[0pt][0pt]{$1$}}}
\put(7201,-5161){\line( 0,-1){600}}
\put(6901,-4561){\line( 1,-2){300}}
\put(7201,-5161){\line( 1, 2){300}}
\put(4201,-3361){\line( 0, 1){1800}}
\put(4201,-5761){\line( 0, 1){1200}}
\put(3901,-2236){\makebox(0,0)[lb]{\raisebox{0pt}[0pt][0pt]{$1$}}}
\put(3901,-2836){\makebox(0,0)[lb]{\raisebox{0pt}[0pt][0pt]{$1$}}}
\put(3901,-4636){\makebox(0,0)[lb]{\raisebox{0pt}[0pt][0pt]{$k$}}}
\put(3901,-5161){\makebox(0,0)[lb]{\raisebox{0pt}[0pt][0pt]{$1$}}}
\put(6826,-2011){\makebox(0,0)[lb]{\raisebox{0pt}[0pt][0pt]{$A$}}}
\put(7426,-2011){\makebox(0,0)[lb]{\raisebox{0pt}[0pt][0pt]{$B$}}}
\put(6826,-4411){\makebox(0,0)[lb]{\raisebox{0pt}[0pt][0pt]{$A$}}}
\put(7426,-4411){\makebox(0,0)[lb]{\raisebox{0pt}[0pt][0pt]{$B$}}}
\end{picture}
\end{center}

\centerline{18.  Three operations of $K^+$}
\medskip
} \end{ex}

\begin{ex} \label{op.binary} \et The free operad on one
unary operation, $F_1$, as the operad for functions.
{\rm Let $F_1$ be the operad with two types, say $x$ and $x'$, and
three operations: a unary operation $f$ with profile $(x,x')$,
and the two identity operations, which we call $1_{x}$ and
$1_{x'}$.  An $F_1$-algebra is simply a function.
} \end{ex}

\begin{ex} \label{op.biaction}\et  $F_1^+$ as the operad for 
monoid bi-actions.
{\rm  Since $F_1$ has three operations, $F_1^+$ has three types:
$f$, $1_x$, and $1_{x'}$.  Following Example \ref{op.monoid.actions},
one may check that an $F_1^+$-algebra $A$ consists of two monoids
$A(1_x)$ and $A(1_{x'})$ together with a set $A(f)$ equipped with an
an action of $A(1_x)^\op \times A(1_{x'})$.  
} \end{ex}

\begin{ex} \label{op.kary} \et The free operad on one $k$-ary
operation, $F_k$, as the operad for $k$-ary multi-functions. {\rm
Generalizing from the previous examples, we let $F_k$ be the operad
with $k+1$ types, say $x_1,\dots,x_k,x'$, and one $k$-ary operation
$f$ with profile  $(x_1,\dots,x_k,x')$, together with the operations
required by the definition of an operad: the $k+1$ identity
operations, which we call $1_{x_1}, \dots, 1_{x_k}, 1_{x'}$, and the 
$k$-ary operations obtained from $f$ by the action of the permutation
group $S_k$.  An $F_k$-algebra $A$ is a collection of sets
$A_1,\dots,A_k,A'$ and a function   from $A_1 \times \cdots \times
A_k$ to $A'$.   We call this a {\it $k$-ary multi-function}.
} \end{ex}

\begin{ex} \label{op.multiaction} \et $F_k^+$ as the operad for 
$(k,1)$ monoid multi-actions.
{\rm An $F_k^+$-algebra consists of $k+1$ monoids $A_1,
\dots, A_k, A'$ and a set equipped with an action of 
\[   A_1^\op \times  \cdots \times A_k^\op \times A'.\]
We call this a {\it $(k,1)$ multi-action} of the monoids
in question, since it can be thought of as $k$ right actions
and one left action, all of which commute.  
} \end{ex}

Since every $S$-operad $O$ may be presented as the quotient of a free
operad on some set of operations, Example \ref{op.multiaction}
suggests the following general picture of $O^+$-algebras.  Given an 
$S$-operad $O$, let us say that an operation $f$ of $O^{n+}$ is {\it
degenerate} if $n = 0$ and $f$ is an identity operation, or if $n > 0$
and $f$ is either an identity operation, a nullary operation, or an
operation with one or more degenerate operations as input types.  For
example, all the operations of $I^{n+}$ are degenerate in this sense.  

\begin{thm} \label{o+.algebras}
\et For any $S$-operad $O$, an $O^+$-algebra $A$ consists of:
\begin{enumerate}  {\rm
\item {\it for each type $x$ of $O$, a monoid $A(x)$}
\item {\it for each nondegenerate operation $g$ of $O$ with profile $(x_1,\dots,
x_k,x')$, a set $A(g)$ equipped with a $(k,1)$ multi-action of the monoids
$A(x_1), \dots, A(x_k), A(x')$}
\item {\it for each nondegenerate reduction law of $O$ --- that is, for
each nondegenerate operation $G$ of $O^+$ with profile
$(g_1,\dots,g_k,g')$ --- a morphism
\[  A(G) \maps G(A(g_1), \dots ,A(g_k)) \to A(g')\]
of multi-actions}
\item {\it for each nondegenerate way of combining reduction
laws of $O$ to obtain another reduction law --- that is, for each
nondegenerate operation $\G$ of $O^{++}$ with profile
$(G_1,\dots,G_k,G')$ --- an equation}
\[  \G(A(G_1),\dots,A(G_k)) = A(G'). \]
} \end{enumerate}
\end{thm}

Proof - First, points 3 and 4 require a bit of clarification.  An
operation $G$ of $O^+$ with profile $(g_1,\dots,g_k,g')$  corresponds
to an $O$-metatree, and this metatree gives a recipe for tensoring the
multi-actions on $A(g_1),\dots, A(g_k)$ in a tree-like pattern,
obtaining a set we denote by $G(A(g_1),\dots,A(g_k))$, equipped with
a multi-action of the same monoids that act on $A(g')$.   Similarly,
an operation $G$ of $O^{++}$ with profile $(G_1,\dots,G_k,G')$
corresponds to a metatree that specifies how to compose the morphisms
$A(G_1), \dots, A(G_k)$ in a tree-like pattern, obtaining a morphism
with the same source and target as $A(G')$, which we denote by
$\G(A(G_1), \dots, A(G_k))$.     

Next, suppose $A$ is an $O^+$-algebra.  We have seen that $A$ consists
of: 
\begin{alphalist}
\item  for each type $g$ of $O^+$, a set $A(g)$
\item for each operation $G$ of $O^+$ with profile
$(g_1,\dots,g_k,g')$, a function
\[  A(G) \maps A(g_1) \times \dots \times A(g_k) \to A(g')\]
\item for each reduction
law of $O^+$ --- that is, for each operation $\G$ of $O^{++}$ with profile
$(G_1,\dots,G_k,G')$ --- an equation  
\[  \G(A(G_1),\dots,A(G_k)) = A(G') \]
\end{alphalist}
where again we use metatree notation to compose the functions
$A(G_1),\dots,A(G_n)$ in the tree-like pattern specified by the
operation $\G$ of $O^{++}$.   In what follows we show how (a)-(c)
give 1-4; by examining our argument one can check that the converse
holds as well.  

Recall first that the types of $O^+$ are the operations of $O$.  These are
either identity operations or nondegenerate operations.  Item (a) applied
to any identity operation $1_x$ of $O$ gives a set which we denote as
$A(x)$.  Item (a) applied to any nondegenerate operation $g$ of $O$ gives
a set $A(g)$.  

Recall next that the operations of $O^+$ are the reduction laws of
$O$.  Any operation of $O^+$ is either an identity operation, a
nullary  operation, an operation with an identity operation of $O$ as
an input type, or a nondegenerate operation.  We consider
these cases in turn.

Item (b) applied to any identity operation $1_g$ of $O^+$ gives a
function from $A(g)$ to itself.  However, (c) applied to the nullary
operation of $O^{++}$ with $1_g$ as output type implies that this
function is the identity.  

There is one nullary operation of $O^+$ with output type $1_x$ for
each type $x$ of $O$.  Item (b) applied to this operation equips the
set $A(x)$ with a distinguished element.  

There are many operations of $O^+$ having an identity operation of $O$
as an input type, but they are all composites of nondegenerate operations
with operations of the following three kinds, so by (c) it suffices
to consider only these three kinds.  First, there are
identity operations $1_{1_x}$ of $O$, which we have already treated. 
Second, there is the binary operation of $O^+$ with profile
$(1_x,1_x,1_x)$.  By (b) it follows that $A(x)$ is equipped with a
binary product, and (c) then implies that $A(x)$ is a monoid with this
product and its distinguished element.   Third, there are the
operations of composing an operation $g$ of $O$ with profile
$(x_1,\dots,x_k,x')$ with the identity operations $1_{x_1}, \dots,
1_{x_k}, 1_{x'}$.   By (b) and (c) it follows that $A(g)$ is equipped
with a $(k,1)$ multi-action of the monoids $A(x_1),\dots,
A(x_k),A(x')$.

Item (b) applied to any nondegenerate operation $G$ of $O^+$ 
with profile $(g_1,\dots,g_k,g')$ gives a function
\[  A(G) \maps A(g_1) \times \dots \times A(g_k) \to A(g')\]
and (c) implies that this function defines a morphism of 
multi-actions
\[  A(G) \maps G(A(g_1), \dots ,A(g_k)) \to A(g'). \]

Recall finally that the operations of $O^{++}$, or reduction laws of
$O^+$, are ways of combining reduction laws of $O$ to give other
reduction laws of $O$.  Applying (c) to an operation $G$ of $O^{++}$ 
with profile $(G_1,\dots,G_k,G')$ we obtain an equation
\[  \G(A(G_1),\dots,A(G_k)) = A(G') .\]
One can check that the equations coming from nondegenerate operations $\G$ 
imply those coming from degenerate operations. 
\qed 

\subsection{Opetopic sets} \label{opetopic}

In topology it is common to take simplices as the basic shapes 
for cells.   There is a category with simplices as objects and face
and degeneracy maps as morphisms.  Presheaves on this category ---
i.e., set-valued functors on the opposite category --- are called
`simplicial sets'.  In our approach to $n$-category theory  we take
opetopes as the basic shapes for cells.  Opetopes form a category, and
presheaves on this category are called `opetopic sets'.   

Here, however, we give a recursive definition of opetopic sets that does
not rely on the category of opetopes.   For this it is convenient to
introduce some notation.  

\begin{defn} \et Given a set $S$, a {\it set over $S$} is
a set $Y$ equipped with a function to $S$.  Given an $S$-operad
$O$ and a set $Y$ over $S$, we define $O_Y$ to be the pullback operad
$F^\ast O$, where $F$ is the function from $Y$ to $S$.  \end{defn}

\noindent We then define opetopic sets as follows:

\begin{defn} \et Given an $S$-operad $O$, an {\rm $O$-opetopic set}
$X$ is defined recursively as a set $X(0)$ over $S$ together with a
$(O_{X(0)})^+$-opetopic set.     \end{defn}

If we work out the implications of this definition, we see that if
$O$ is an $S$-operad, an $O$-opetopic set $X$ consists of an
set $X(n)$ over $S(n)$ for each integer $n \ge 0$, where 
\[      S(0) = S, \qquad S(n+1) = \elt(O(n)_{X(n)}), \]
and $O(n)$ is the $S(n)$-operad given by
\[      O(0) = O, \qquad O(n+1) = (O(n)_{X(n)})^+ .\]
Note also that
\[    S(n) = \type(O(n)). \]

\begin{defn}\et  Let $O$ be an $S$-operad and $X$ an $O$-opetopic set. 
We define an {\rm $n$-dimensional cell} (or {\rm $n$-cell}) of $X$ to
be an element of $X(n)$.  We define an {\rm $n$-dimensional frame} in
$X$ to be an element of $S(n)$.  For $n \ge 1$, we define an {\rm
$n$-dimensional opening} in $X$ to be an operation of $O(n-1)$. 
\end{defn}

Since $X(n)$ is a set over $S(n)$, there is a map from $n$-dimensional
cells to $n$-dimensional frames, and for any cell of $X$ we may speak
of the frame {\it of} that cell.  Also, for $n \ge 1$, the tautologous
morphism from the pullback of an operad to the operad itself gives a
map from operations of $O(n-1)_{X(n-1)}$, which are $n$-dimensional
frames, to operations of $O(n-1)$, which are $n$-dimensional openings. 
Thus for $n \ge 1$ we may speak of any frame $s$ of $X$ as being {\it in}
some opening $o$, and given any cell $x$ with frame $s$, we also say
that $x$ is {\it in} $o$.  

Let $o$ be an $n$-dimensional opening in $X$.  
We define an {\it $o$-cell} to be a cell in $o$.   The frame of $x$ is
an operation of $O(n)_{X(n)}$, and has profile $(a_1,\dots,a_k,b)$ for
some $(n-1)$-dimensional cells $a_1,\dots,a_k,b$.   It is convenient
to use the following schematic picture of $x$: 
\[ 
\begin{diagram}[(a_1,\dots,a_k)] 
\node{(a_1,\dots,a_k)}\arrow{x}\node{b} 
\end{diagram} \] 
We call $a_1,\dots,a_k$ the {\it infaces} of $x$, and $b$ the {\it
outface} of $x$.     

Similarly, we define an {\it $o$-frame} to be a frame in $o$, 
and depict an $o$-frame with profile $(a_1,\dots,a_k,b)$
as follows:
\[ 
\begin{diagram}[(a_1,\dots,a_k)]
\node{(a_1,\dots,a_k)}\arrow{e,t}{?} \node{b}
\end{diagram} 
\]

An `$o$-niche' is like an $o$-frame with the outface missing.  
Suppose that the opening $o$ has profile
$(s_1,\dots,s_k,t)$.  We define an {\it $o$-niche} to be a tuple
$(a_1,\dots,a_k)$ of $(n-1)$-dimensional cells with $a_i$ having $s_i$
as its frame.   We depict this $o$-niche as follows:
\[ 
\begin{diagram}[(a_1,\dots,a_k)]
\node{(a_1,\dots,a_k)}\arrow{e,t}{?} \node{?}
\end{diagram} 
\]
The concept of niche serves as our
substitute for the concept of a horn in a simplicial set.

Similarly, a `punctured $o$-niche' is like an $o$-frame with the outface
and one inface missing.   We define a {\it punctured $o$-niche} to be
a tuple $(a_1,\dots,a_{j-1},a_{j+1},\dots,a_k)$ of cells with $a_i$
having $s_i$ as its frame, and depict this as:
\[ 
\begin{diagram}[(a_{1},a_{j+1},\dots,a_{k})]
\node{(a_{1},\dots,a_{j-1},?,a_{j+1},\dots,a_{k})}\arrow{e,t}{?} \node{?}
\end{diagram} 
\]

In the case where one of these configurations ($o$-frame, $o$-niche,
or punctured $o$-niche) can be extended to an actual $o$-cell, the
$o$-cell is called an {\it occupant} of the configuration.  Occupants
of the same frame (resp.\ niche) are called {\it frame-competitors}
(resp.\ {\it niche-competitors}).

To make $O$-opetopic sets into a category we need to define morphisms
between them.  Roughly speaking, a morphism $\phi \maps X \to X'$
between $O$-opetopic sets is a function sending each cell $x \in X(n)$
to a cell $\phi(x) \in X'(n)$ of the same shape, such that $\phi$ of
any face of $x$ is the corresponding face of $\phi(x)$.  To make this
precise requires a bit of technical work.

We begin with some remarks on the functoriality of the slice operad
construction.  Suppose $O$ is an $S$-operad, $O'$ is an $S'$-operad, and
$F \maps S \to S'$ is a function.  By Proposition
\ref{lax.monoidal.functor} we obtain a lax monoidal functor $F^\ast
\maps \sig(S') \to \sig(S)$.   As in Section \ref{monoid.objects} this
allows us to speak of morphisms from $S'$-operads to $S$-operads, but we
can also define morphisms going the other way.  Namely, we define an
{\it operad morphism} $f \maps O \to O'$ {\it riding $F$} to be an
operad homomorphism $f \maps O \to F^\ast(O')$.   

Given such an operad morphism there is an obvious function from from
$|\elt(O)|$ to $|\elt(O')|$, which we call $F^+$.  We also obtain a
operad morphism $f^+ \maps O^+ \to O'^+$ riding this function. To see
this it is easiest to use metatree notation:  an operation of $O^+$ is
given by a 1-dimensional $O$-metatree, and using $f \maps O \to O'$
one can convert this to a 1-dimensional $O'$-metatree, which specifies an
operation of $O'$ and thus of ${F^+}^\ast(O')$.  One can then
check this defines an operad morphism $f^+ \maps O^+ \to O'^+$.  

Now suppose that $Y$ is a set over $S$ and $Y'$ is a set over $S'$. 
We define a {\it function $\phi \maps Y \to Y'$ over} $F \maps S
\to S'$ to be a function making the following diagram commute:
\[   
\begin{diagram}[Y']
\node{Y} \arrow{e,t}{\phi} \arrow{s} \node{Y'} \arrow{s} \\
\node{S} \arrow{e,t}{F} \node{S'} 
\end{diagram} 
\]  
Given an operad morphism $f \maps O \to O'$ riding $F$, 
there is an obvious operad morphism from $O_Y$ to 
$O'_{Y'}$ riding $\phi$, which we call $f_\phi$.   

Finally, suppose that $X$ is an $O$-opetopic set and $X'$ is an
$O'$-opetopic set.   Suppose that $f \maps O \to O'$ is an operad
morphism riding $F \maps S \to S'$.   We define an {\it opetopic map}
$\phi \maps  X \to X'$ {\it riding} $f$ to consist of,  for each $n
\ge 0$, a function
\[       \phi_n \maps X(n) \to X'(n)  \]
over the function 
\[       F_n \maps S(n) \to S'(n)   \] 
given as follows.  We set $F_0 = F$, and define $F_n$ for higher $n$
recursively, along with a sequence of operad morphisms  
\[      f_n \maps O(n) \to O'(n),  \] 
starting with $f_0 = f$.   To do so, we let
\[   f_{n+1} = ((f_n)_{\phi_n})^+ \] 
and note that this operad morphism gives a map from $S(n+1)$ to
$S'(n+1)$, which we take as $F_{n+1}$.   Unrolling this recursive
construction one sees that, fixing $f$ and $F$, the morphism $\phi
\maps X \to X'$ is completely determined by the functions $\phi_n$
sending $n$-cells of $X$ to $n$-cells of $X'$.

\begin{defn} \et Given an $S$-operad $O$, we define the {\rm category of
$O$-opetopic sets} to be that with $O$-opetopic sets as objects
and opetopic morphisms riding the identity function as morphisms.  
\end{defn}

\noindent In fact, this category is equivalent to the category of presheaves on
a certain category of $O$-opetopes.  To save space we shall not prove
this here, but only seek to make it plausible by showing that every
$n$-cell of an $O$-opetopic set $X$ has some $n$-dimensional opetope
as its `shape'.  This is trivial in the case $n = 0$, so we assume $n
\ge 1$.

Recall that every $n$-dimensional cell of $X$ is in some opening,
which is an operation of $O(n-1)$.  On the other
hand, each $n$-dimensional opetope is an operation of $O^{(n-1)+}$.
Thus to associate an $n$-dimensional opetope to each $n$-cell
of $X$, we construct, for all $n \ge 0$, an operad morphism
\[      p_n \maps  O(n) \to O^{n+} . \]
Since $O(0) = O$, we take $p_n$ to be the identity when
$n = 0$.  Given $p_n$, to 
define $p_{n+1}$ we first form the composite
\[  
\begin{diagram}[O(n)_{X(n)}]
\node{O(n)_{X(n)}} \arrow{e} \node{O(n)} \arrow{e,t}{p_n} \node{O^{n+}} 
\end{diagram}
\]
where the first arrow is the tautologous morphism 
from a pullback of an operad to the operad itself.   Taking the `+' of
this composite, we then obtain $p_{n+1}$.  

\section{$n$-Categories} \label{n-categories}

In Section \ref{coherent} we define `$n$-coherent $O$-algebras'.   The
basic idea is that for any operad $O$, an $n$-coherent $O$-algebra is
an $n$ times categorified analog of an $O$-algebra.  For example, just
an $I$-algebra is a set, an $n$-coherent $I$-algebra is an
$n$-category.    Other examples are also interesting: just as an
$I^+$-algebra is a monoid, an $n$-coherent $I^+$-algebra is a
`monoidal $n$-category', and just as $T$-operad is a commutative
monoid, an $n$-coherent $T$-algebra is a `stable $n$-category'. 
Stable $n$-categories play an important role in the program sketched
in HDA0, and also in the foundations of $n$-category theory itself,
since the $(n+1)$-category of all $n$-categories will be a
stable $(n+1)$-category.  

In Section \ref{n-functor} we define `$k$-ary virtual $n$-functors' to
be $n$-coherent $F_k$-algebras, where $F_k$ is the free operad on one
$k$-ary operation.  This concept allows us to reinterpret and  clarify
some of the previous material.  For example, in Theorem 
\ref{n-coherent.o-alg} we use them to give a recursive
characterization of $n$-coherent $O$-algebras that is often more
useful than the original definition.  We also use them in Propositions
\ref{balanced} and \ref{universal} to characterize the concepts of
`balanced' punctured niche and `universal' niche-occupant, introduced
in the previous section. Finally, in Section \ref{microcosm} we give a
rather general precise statement of the `microcosm principle'.  

\subsection{$n$-Coherent $O$-algebras} \label{coherent}

In what follows we fix a nonnegative integer $n$ and define
the notion of `$n$-coherent $O$-algebra', which will be an
$O$-opetopic set with certain properties.  To do so, we need
the notions of `balanced punctured niche' and `universal
niche-occupant', which we define in a recursively interlocking way.

As the definitions are a bit complicated, let us first explain them
in a heuristic way.  We shall see in Section \ref{n-functor} 
that in an $n$-coherent $O$-algebra, any $m$-dimensional punctured
niche 
\[ 
\begin{diagram}[(a_{1},a_{j+1},\dots,a_{k})]
\node{(a_{1},\dots,a_{j-1},?,a_{j+1},\dots,a_{k})}\arrow{e,t}{?} \node{?}
\end{diagram}
\]
determines a `virtual $(n-m)$-functor'.  In Proposition \ref{balanced} 
we show that the punctured niche is balanced if and only if this virtual
$(n-m)$-functor is an `equivalence'.  On the other hand, for a
niche-occupant   
\[
\begin{diagram}[(c_1,\dots,c_k)]
\node{(c_1,\dots,c_k)}  \arrow{e,t}{u} \node{d}
\end{diagram}
\]
to be `universal' means roughly that any other occupant of the same niche 
factors through the given one --- at least `up to equivalence'.   We
make this precise in Proposition \ref{universal}.

The definitions are as follows:

\begin{defn} \et For an $m$-dimensional opening $o$, a punctured $o$-niche:
\[ 
\begin{diagram}[(a_{1},a_{j+1},\dots,a_{k})]
\node{(a_{1},\dots,a_{j-1},?,a_{j+1},\dots,a_{k})}\arrow{e,t}{?} \node{?}
\end{diagram}
\]
is said to be {\rm balanced} if and only if $m>n+1$ or:
\begin{enumerate}
\item any extension
\[ 
\begin{diagram}[(a_{1},a_{j+1},\dots,a_{m})]
\node{(a_{1},\dots,a_{j-1},?,a_{j+1},\dots,a_{k})}\arrow{e,t}{?} \node{b}
\end{diagram} 
\]
extends further to:
\[ 
\begin{diagram}[(a_{1},a_{j+1},\dots,a_{k})]
\node{(a_{1},\dots,a_{j-1},a_j,a_{j+1},\dots,a_{k})}\arrow{e,t}{u} \node{b}
\end{diagram} 
\]
with $u$ universal in its niche, and 
\item  for any occupant
\[ 
\begin{diagram}[(a_{1},a_{j+1},\dots,a_{k})]
\node{(a_{1},\dots,a_{j-1},a_j,a_{j+1},\dots,a_{k})}\arrow{e,t}{u} \node{b}
\end{diagram} 
\]
universal in its niche, and frame-competitor $a'_{j}$ of $a_{j}$, the
$(m+1)$-dimensional punctured niches:
\[
\begin{diagram}
[((c_{1},\dots,c_{m}) \to d,\; d \to d')]
\node{(a'_{j} \mapright{?} a_{j}, \;
(a_{1},\dots,a_{j-1},a_{j},a_{j+1},\dots,a_{k}) \mapright{u} b)}  
\arrow{s,r}{?}  \\
\node{(a_{1},\dots,a_{j-1},a'_{j},a_{j+1},\dots, a_{k}) \mapright{?} b}
\end{diagram}
\]
and
\[
\begin{diagram}
[((c_{1},\dots,c_{m}) \to d,\; d \to d')]
\node{((a_{1},\dots,a_{j-1},a_{j},a_{j+1},\dots,a_{k}) \mapright{u} b, \;
a'_{j} \mapright{?} a_{j})} 
\arrow{s,r}{?}  \\
\node{(a_{1},\dots,a_{j-1},a'_{j},a_{j+1},\dots, a_{k}) \mapright{?} b}
\end{diagram}
\]
are balanced.  
\end{enumerate}
\end{defn}

\begin{defn}\et An $m$-dimensional niche-occupant:
\[
\begin{diagram}[(c_1,\dots,c_k)]
\node{(c_1,\dots,c_k)}  \arrow{e,t}{u} \node{d}
\end{diagram}
\]
is said to be {\rm universal} if and only if $m>n$ and $u$ is its own unique
niche-competitor, or $m\le n$ and for any frame-competitor $d'$ of $d$, the
$(m+1)$-dimensional punctured niches:
\[
\begin{diagram}[((c_{1},\dots,c_{k}) \to d,\; d \to d')]
\node{((c_{1},\dots,c_{k}) \mapright{u} d,\; d \mapright{?} d')} 
\arrow{s,r}{?} \\
\node{(c_{1},\dots,c_{k}) \mapright{?} d'}  
\end{diagram}
\]
and
\[
\begin{diagram}[((c_{1},\dots,c_{k}) \to d,\; d \to d')]
\node{(d \mapright{?} d',\;(c_{1},\dots,c_{k}) \mapright{u} d)} 
\arrow{s,r}{?} \\
\node{(c_{1},\dots,c_{k}) \mapright{?} d'}  
\end{diagram}
\]
are balanced.  
\end{defn}

\begin{defn}\et  Given a universal $o$-cell: 
\[ 
\begin{diagram}[(a_1,\dots,a_k)]
\node{(a_1,\dots,a_k)}\arrow{e,t}{u} \node{b}
\end{diagram} 
\]
we call $b$ a {\rm composite} of $(a_{1},\dots,a_{k})$, or
{\rm $o$-composite} if we need to be more specific.  
\end{defn}

\begin{defn} \et An {\rm $n$-coherent $O$-algebra} is an $O$-opetopic
set such that 1) every niche has a universal occupant, and 2) composites
of universal cells are universal.  \end{defn}

\noindent The dependence on $n$ in this definition is implicit in how
the definition of `universal' depends on $n$.   Note that in an
$n$-coherent $O$-algebra, for $m > n$ every $m$-dimensional niche has
a unique occupant, which is automatically universal,  and for $m >
n+1$ every $m$-dimensional punctured niche is balanced. One can also
check that for $m > n + 1$ every $m$-dimensional frame has a unique
occupant.  This is analogous to how a Kan complex represents an
$n$-groupoid if, for $m > n+1$, any configuration in which all the
faces of $m$-simplex are filled in by $(m-1)$-dimensional cells in a
consistent way can be uniquely extended to a $m$-dimensional cell.   

A 0-coherent $O$-algebra is essentially the same thing as an
$O$-algebra.  Given a 0-coherent $O$-algebra $A$, the types of $O$ are
the 0-dimensional frames of $A$, so for any type $s$ there is a set
$\tilde A(s)$ of 0-cells of $A$ having $s$ as frame.  
For any operation $f$ of $O$ with profile $(s_1,\dots,s_k,s')$, and
any 0-cells $a_i \in \tilde A(s_i)$, the 1-dimensional niche
\[ 
\begin{diagram}[(a_1,\dots,a_k)]
\node{(a_1,\dots,a_k)}\arrow{e,t}{?} \node{?}
\end{diagram} 
\]
has a unique occupant
\[ 
\begin{diagram}[(a_1,\dots,a_k)]
\node{(a_1,\dots,a_k)}\arrow{e,t}{u} \node{a'.}
\end{diagram} 
\]
Thus one can check that there is an $O$-algebra $\tilde A$ with
the sets $\tilde A(s)$ given as above, and with
the operation $f$ acting by 
\[     f(a_1,\dots,a_k) = a' .\]
In fact, one can check that this construction gives an equivalence
between the category of $O$-algebras and the category of 0-coherent
$O$-algebras in which morphisms are defined as follows:

\begin{defn} \et Let $O$ be an $S$-operad and let $A,A'$ be
$n$-coherent $O$-algebras.   Then a morphism of $O$-opetopic sets $f
\maps A \to A'$ is called an {\rm $n$-coherent $O$-algebra morphism}
if it preserves universality of niche-occupants.   \end{defn}

We study $n$-coherent $O$-algebras for higher $n$
in the following two sections.   In Theorem \ref{n-coherent.o-alg} we
recursively describe $n$-coherent $O$-algebras in terms of
$(n-1)$-coherent $O$-algebras.  In Theorem \ref{1-coherent.o-alg} we
use this to give a concrete description of 1-coherent $O$-algebras.

The simplest sort of operad algebra is an $I$-algebra, which 
by Example \ref{init.op} is just a set.  Similarly, the
simplest sort of $n$-coherent $O$-algebra is an $n$-category:

\begin{defn} \et An {\rm $n$-category} is an $n$-coherent $I$-algebra.
An {\rm $n$-functor} is a morphism of $n$-coherent $I$-algebras.
\end{defn}  

\begin{ex} \et \label{1.cat}
1-categories as categories. 
{\rm A 1-coherent $I$-algebra $C$ 
has a set $C(0)$ of 0-cells, and given 0-dimensional cells $c$ and $c'$
we may denote the set of occupants of the frame
\[ 
\begin{diagram}[c']
\node{c}\arrow{e,t}{?} \node{c'}
\end{diagram} 
\]
as $\hom(c,c')$.  Given a 0-cell $c$ the 2-dimensional niche
\[
\begin{diagram}[c \mapright{?} c']
\node{} 
\arrow{s,r}{?} \\
\node{c \mapright{?} c}  
\end{diagram}
\]
has a unique occupant 
\[
\begin{diagram}[c \mapright{?} c']
\node{} 
\arrow{s,r}{u} \\
\node{c \mapright{1_c} c}  
\end{diagram}
\]
so we have $1_c \in \hom(c,c)$.  Similarly, given 0-cells $c,c',c''$,
the 2-dimensional niche 
\[
\begin{diagram}[c \mapright{f} c',\; c' \mapright{g} c'']
\node{(c \mapright{f} c',\; c' \mapright{g} c'')} 
\arrow{s,r}{?} \\
\node{c \mapright{?} c''}  
\end{diagram}
\]
has a unique occupant 
\[
\begin{diagram}[c \mapright{f} c',\; c' \mapright{g} c'']
\node{c \mapright{f} c',\; c' \mapright{g} c''} 
\arrow{s,r}{u} \\
\node{c \mapright{fg} c''}  
\end{diagram}
\]
so given $f \in \hom(c,c'), g \in \hom(c',c'')$ we get $fg \in
\hom(c,c'')$.  By examining the 3-dimensional cells of $C$ one can
check that these operations give a category $\tilde C$ with
$C(0)$ as its set of objects and the sets $\hom(c,c')$ as 
hom-sets.  One can also check that this construction gives an
equivalence between the category with 1-categories as objects and
1-functors as morphisms, and the category with small categories as
objects and functors as morphisms.  
}
\end{ex}

In Examples \ref{op.monoids} and \ref{op.comm.monoids} we saw that an
$I^+$-algebra is a monoid, and a $T$-algebra is a commutative monoid. 
By analogy we make the following definitions:

\begin{defn} \et A {\rm monoidal $n$-category} is an $n$-coherent
$I^+$-algebra.  \end{defn}  

\begin{defn} \et A {\rm stable $n$-category} is an $n$-coherent $T$-algebra. 
\end{defn}  

Since there are unique operad homomorphisms from $I$ to $I^+$ and from
$I^+$ to $T$, the following result lets us extract an $n$-category
from any monoidal $n$-category, and a monoidal $n$-category from any
stable $n$-category.

\begin{prop} \label{pullback.coherent.algebra} \et  Suppose $O$ is an
$S$-operad, $O'$ is an $S'$-operad, $F \maps S \to S'$ is a
function, and $f \maps O \to O'$ is an operad morphism riding $F$. 
Suppose $X'$ is an $O'$-opetopic set and $X = f^\ast X'$ is the pullback
$O$-opetopic set.  Then a punctured niche in $X$ is balanced if and
only if the corresponding punctured niche in $X'$ is balanced, and a
niche-occupant in $X$ is universal if and only the corresponding
niche-occupant in $X'$ is universal.  Thus $X$ is an
$n$-coherent $O$-algebra if $X'$ is an $n$-coherent $O'$-algebra.  
\end{prop}

Proof - The proof is a straightforward verification once we have
clarified the notion of `pullback' used here.
Suppose that $O$ is an $S$-operad and $O'$ is an $S'$-operad.  Given an
operad morphism $f \maps O \to O'$ riding a function $F \maps S \to
S'$, the pullback $X = f^\ast X'$ of an $O'$-opetopic set
$X'$, which is an $O$-opetopic set.   The set $X(0)$ over $S$ is
defined to be the pullback of the set $X'(0)$ over $S$, and the
underlying $(O_{X(0)})^+$-opetopic set of $X$ is defined (recursively)
to be the pullback of the underlying $(O_{X'(0)})^+$-opetopic set of
$X'$.   \qed

In a future paper we plan to discuss the stable $(n+1)$-category of
$n$-categories, $n$Cat.  This is needed for most of the interesting
applications of $n$-category theory.  The 1-cells in $n$Cat are 
`$k$-ary virtual functors'.  We study a version of these in the
following section, defined in a way that is convenient now but not
necessarily best in the long run.  

\subsection{$k$-ary virtual $n$-functors} \label{n-functor}

As we saw in Example \ref{op.kary}, the free operad on one $k$-ary
operation, $F_k$, is the operad for $k$-ary multi-functions. 
By analogy we make the following definition:

\begin{defn} \et A {\rm $k$-ary virtual $n$-functor} is an
$n$-coherent $F_k$-algebra.  We omit the term `$k$-ary' if $k = 1$,
and the reference to $n$ if $n = 1$.
\end{defn}

Suppose that $A$ is $k$-ary virtual $n$-functor.  
Recall that $F_k$ has one operation $f$ of type $(x_1,\dots,x_k,x')$,
together with $k+1$ identity operations $1_{x_1}, \dots, 1_{x_k},$ and
$1_{x'}$.  Thus there are $k+1$ operad morphisms from $I$ to $F_k$, and 
by Proposition \ref{pullback.coherent.algebra}, the pullback of
$A$ along any one of these is an $n$-category.  Calling these
$n$-categories $C_1,\dots,C_k$ and $C'$, respectively, we say that $A$ is a
$k$-ary virtual $n$-functor {\it from} $C_1 \times \dots \times C_k$
{\it to} $C'$, and write
\[          A \maps C_1 \times \dots \times C_k \rightharpoonup C'.\]

\begin{ex} \label{anafunctor} \et Virtual functors as saturated anafunctors.  
{\rm A virtual functor is essentially the same as what Makkai
\cite{Makkai} calls a `saturated anafunctor', which may be
viewed as a special sort of
distributor.   A {\it distributor} $A$ from the category $C$ to the
category $D$ is a functor $A \maps C^\op \times D \to \Set$, and $A$
is a {\it saturated anafunctor} if for every object $c \in C$, the functor
$A(c,\cdot)$ is naturally isomorphic to $\hom(d,\cdot)$ for some object
$d \in D$.  Thus, in keeping with the philosophy of this paper, a saturated
anafunctor does not specify a unique object $d \in D$ for each
object $c \in C$.  Instead, it specifies a universal property, which
automatically determines an object $d \in D$ up to a specified
isomorphism.  

Suppose that $A \maps C \rightharpoonup D$ is a virtual functor.
Then we obtain 1-categories $C$ and $D$, which by Example \ref{1.cat}
we may think of as categories.  Given objects $c \in C, d \in D$, 
we denote the set of occupants of the $f$-frame
\[ 
\begin{diagram}[c']
\node{c} \arrow{e,t}{?} \node{d}
\end{diagram} 
\]
by $\tilde A(c,d)$.  Since 1-cells in a 1-coherent
$O$-algebra have unique composites, any morphism $f \maps c \to
c'$ in $C$ gives a function
\[     \tilde A(c',d) \to \tilde A(c,d) \]
for each $d \in D$, and any morphism $f \maps d \to d'$ gives a
function 
\[     \tilde A(c,d) \to \tilde A(c,d') \]
for each $c \in C$.  Thus $\tilde A$ can be thought of as 
a distributor from $C$ to $D$.  Because every $f$-niche
\[ 
\begin{diagram}[c']
\node{c} \arrow{e,t}{?} \node{?}
\end{diagram} 
\]
has a universal occupant, $\tilde A$ is a saturated anafunctor.   
Conversely, every saturated anafunctor can be thought of as a
virtual functor.
}\end{ex}

\begin{ex}\et 0-ary virtual functors as representable presheaves.
{\rm Generalizing the previous example, one can show that $k$-ary 
virtual functors are essentially the same as `$k$-ary saturated anafunctors'.
The case $k = 0$ is particularly interesting.  A $0$-ary virtual functor
with codomain $C$ is just a functor $P \maps C \to \Set$ that is
naturally isomorphic to $\hom(c,\cdot)$ for some $c \in C$.  This is 
also called a `representable presheaf' on $C^\op$.

Recall from Example \ref{op.pointed.sets} that we gave the operad
$F_0$ another name, $K$.  This stands for `constant', since a
$K$-algebra is just a pointed set.  Here we see that a 1-coherent
$K$-algebra, or in other words a representable presheaf, is a
categorified version of a pointed set: it is a category equipped,
not quite with a distinguished object, but with a universal property
that determines an object up to natural isomorphism.  

Generalizing, we call an $n$-coherent $K$-algebra a {\it representable
$n$-prestack}.   It follows from Theorem \ref{n-coherent.o-alg} that
an $n$-prestack $P$ may be regarded as a special sort of {\it
$n$-prestack}, which we define as an $(n-1)$-coherent
$(K_{P(0)})^+$-algebra.  We expect that prestacks are to the `stacks'
sought by Grothendieck \cite{Gro} as presheaves are to sheaves.  
}\end{ex}

As noted earlier, the concept of balanced punctured niche is closely
related to the concept of `equivalence'.  We can now begin to make
this more precise:

\begin{defn} \et A virtual $n$-functor $A \maps C \rightharpoonup C'$ is
an {\rm $n$-equivalence}, or simply an {\rm equivalence}, if the
punctured $f$-niche  
\[ 
\begin{diagram}[?]
\node{?} \arrow{e,t}{?} \node{?}
\end{diagram} 
\]
is balanced.  
\end{defn}

A functor is an equivalence if and only if it is essentially
surjective and fully faithful.  The same is true for virtual
$n$-functors.  Note the similarity of the following two definitions to
the two clauses in the definition of `balanced':

\begin{defn} \et A virtual $n$-functor $A \maps C \rightharpoonup C'$ is
{\rm essentially surjective} if any extension 
\[ 
\begin{diagram}[c']
\node{?} \arrow{e,t}{?} \node{c'}
\end{diagram} 
\]
of the punctured $f$-niche extends further to
\[ 
\begin{diagram}[c']
\node{c} \arrow{e,t}{u} \node{c'}
\end{diagram} 
\]
with $u$ universal in its niche. 
\end{defn}

\begin{defn} \et A virtual $n$-functor $A \maps C \rightharpoonup C'$ is
{\rm fully faithful} if for any universal occupant 
\[ 
\begin{diagram}[c']
\node{c} \arrow{e,t}{u} \node{c'}
\end{diagram} 
\]
of the punctured $f$-niche, and any niche-competitor $b$ of
$c$, the punctured niches
\[
\begin{diagram}[(b \mapright{?} c, \; c \mapright{u} c')]
\node{(b \mapright{?} c, \; c \mapright{u} c')}  
\arrow{s,r}{?}  \\
\node{b \mapright{?} c'}
\end{diagram}
\]
and
\[
\begin{diagram}[(b \mapright{?} c, \; c \mapright{u} c')]
\node{(c \mapright{u} c', \; b \mapright{?} c)}  
\arrow{s,r}{?}  \\
\node{b \mapright{?} c'}
\end{diagram}
\]
are balanced.  
\end{defn}

\begin{prop} \et A virtual $n$-functor is an equivalence if and
only if it is essentially surjective and fully faithful.  \end{prop}

Proof - A straightforward consequence of the definitions.  \qed

To further explain the relation between balanced punctured niches and
equivalences, we need the following characterization of $n$-coherent
$O$-algebras.   Recall that if $O$ is an $S$-operad, an $O$-opetopic
set $A$ consists of a set $A(0)$ over $S$ together with
an $(O_{A(0)})^+$-opetopic set.   

\begin{thm} \label{n-coherent.o-alg} \et Suppose that $O$ is an
$S$-operad.  For any $n \ge 1$, an $O$-opetopic set $A$ is an
$n$-coherent $O$-algebra if and only if:  
{\rm 
\begin{enumerate}
\item {\it The underlying $(O_{A(0)})^+$-opetopic set of $A$ 
is an $(n-1)$-coherent $(O_{A(0)})^+$-algebra. }
\item {\it For any $k$-ary operation of $O$, the pullback of $A$
along the resulting operad morphism from $F_k$ to $O$ is a $k$-ary
virtual $n$-functor.}
\item {\it Composites of universal 1-cells in $A$ are universal.}  
\end{enumerate} }
\noindent 
\end{thm}

Proof - We denote the underlying $(O_{A(0)})^+$-opetopic set of $A$ by
$A^-$.  Suppose that $A$ is an $n$-coherent $O$-algebra: in other
words, every niche of $A$ has a universal occupant, and composites of
universal niche-occupants are universal.   One can check using the
formalism developed in Section \ref{opetopic} that for $m \ge 1$, the
$m$-dimensional frames (resp.\ cells) of $A$ correspond to the
$(m-1)$-dimensional frame (resp.\ cells) of $A^-$.  The same is also
true for openings, niches and punctured niches when $m \ge 2$.  Also, the
definitions of `balanced' and `universal' are set up so that an
$m$-dimensional punctured niche of $A$ is balanced if and only if the
corresponding punctured niche of $A^-$ is balanced, and an
$m$-dimensional niche-occupant of $A$ is universal if and only if the
corresponding niche-occupant of $A^-$ is universal.  Thus 1 holds. 
Proposition \ref{pullback.coherent.algebra} implies 2, and 3 is
immediate.  

Conversely, suppose that 1, 2, and 3 hold.  By 1, for $m \ge 2$
every $m$-dimensional niche of $A$ has a universal occupant, and composites
of $m$-dimensional universal niche-occupants are universal.  The
former also holds for $m = 1$ by 2, and the latter holds for $m = 1$
by 3.  \qed

Let $O$ be an $S$-operad and let $A$ be an $n$-coherent
$O$-algebra.  Given $m \le n$, we now describe how:
\begin{enumerate}
\item Every $m$-dimensional frame in $A$ determines an
$(n-m)$-category.  
\item For $m \ge 1$, every $m$-dimensional opening in $A$ determines a
$k$-ary virtual \hfill\break $(n-m+1)$-functor.
\item For $m \ge 1$, every $m$-dimensional punctured niche in $A$ determines
a virtual \hfill\break $(n-m+1)$-functor.
\item For $m \ge 1$, every $m$-dimensional niche in $A$ determines a
representable \hfill\break $(n-m+1)$-prestack.
\end{enumerate}  
For example, when $A$ is an $n$-category there is a unique
1-dimensional frame for any pair of 0-cells $a,b$ in $A$, and we
denote the corresponding $(n-1)$-category by $\hom(a,b)$.  

As in Section \ref{opetopic}, let 
\[      S(0) = S, \qquad S(i+1) = \elt(O(i)_{A(i)}), \]
where $O(i)$ is the $S(i)$-operad given by
\[      O(0) = O, \qquad O(i+1) = (O(i)_{A(i)})^+ .\]
Also let $A^{0-} = A$, and let $A^{(i+1)-}$ be the underlying 
$O(i+1)$-opetopic set of the $O(i)$-opetopic set $A^{i-}$.  
By Theorem \ref{n-coherent.o-alg}, $A^{i-}$ is an $(n-i)$-coherent
$O(i)$-algebra if $i \le n$.  By remarks in the proof of the theorem,
the $m$-dimensional cells (resp.\ frames) of $A$ correspond to the
$(m-i)$-dimensional cells (resp.\ frames) of $A^{i-}$ if $i \le m$, 
and the same is true for openings, niches, and punctured niches if
$i < m$.

Using this `level-shifting' trick, to deal with 1-4 above it suffices
to explain how:
\begin{enumerate}
\item Every $0$-dimensional frame in $A$ determines an
$n$-category.  
\item Every $1$-dimensional opening in $A$ determines a $k$-ary
virtual $n$-functor.
\item Every $1$-dimensional punctured niche in $A$ determines
a virtual $n$-functor.
\item Every $1$-dimensional niche in $A$ determines a representable
$n$-prestack.
\end{enumerate}

For 1, note that a $0$-dimensional frame in $A$ is just an element
$s$ of the set $S$ of types of $O$.  This determines a unique operad
morphism from $I$ to $O$ riding the function $F \maps 1 \to S$
that sends the one element of $1$ to $s$.  The pullback of $A$ under
this morphism is the desired $n$-category. 

For 2, recall that a $1$-dimensional opening in $A$ is simply
an operation of $O$.  As noted in Theorem \ref{n-coherent.o-alg}, 
any $k$-ary operation $o$ of $O$ determines an operad morphism
from $F_k$ to $O$, and the pullback of $A$ under this morphism is a
$k$-ary virtual $n$-functor, say 
\[       G \maps C_1 \times \cdots \times C_k \rightharpoonup C' \]

For 3 and 4, note that if we fix an operation $o$ of $O$,
an $o$-niche then consists of a choice of one 0-cell from each of the
$n$-categories $C_i$, while a punctured $o$-niche consists of a
choice of 0-cells from all but one of the $C_i$.   Thus it
suffices to explain how to extract a $(k - \ell)$-ary virtual
$n$-functor from $G$ by choosing 0-cells in $\ell$ of the
$n$-categories $C_i$.   By induction it suffices to consider the case
$\ell = 1$, so supposing without loss of generality that we have chosen a
0-cell $c_k \in C_k$, let us construct a $(k-1)$-ary virtual $n$-functor 
\[    H \maps C_1 \times \cdots \times C_{k-1} \rightharpoonup C' .\]

By Theorem \ref{n-coherent.o-alg}, $G$ gives an $(n-1)$-coherent
$((F_k)_{G(0)})^+$-algebra $G^-$.  Concretely, $G(0)$ is the
$(k+1)$-tuple of disjoint sets $(C_1(0), \dots, C_k(0),C'(0))$,
where each $C_i(0)$ is the set of 0-cells of the corresponding
$n$-category $C_i$, and $C'(0)$ is the set of 0-cells of the $n$-category
$C'$.   To construct $H$, we first construct an $(n-1)$-coherent
$((F_{k-1})_{H(0)})^+$-algebra $H^-$, where $H(0)$ is the $k$-tuple of
disjoint sets $(C_1(0), \dots, C_{k-1}(0),$ $ C'(0))$.  Note that there
is a unique operad morphism
\[       f \maps (F_{k-1})_{H(0)} \to (F_k)_{G(0)} \]
sending each operation with profile $(c_1,\dots,c_{k-1},c')$ to 
the unique operation with profile $(c_1,\dots,c_k,c')$.  This gives
an operad morphism 
\[       f^+ \maps ((F_{k-1})_{H(0)})^+ \to ((F_k)_{G(0)})^+ ,\]
and we define 
\[       H^- = (f^+)^\ast G^- .\]
Together with $H(0)$, $H^-$ defines an $F_{k-1}$-opetopic set $H$.
To see that $H$ is an $n$-coherent $F_{k-1}$-algebra, it suffices
to check that 1-dimensional niches have universal occupants, and that
composites of universal 1-cells are universal.  These follow from the
corresponding properties for $G$.  Thus $H$ is a $(k-1)$-ary
virtual $n$-functor as desired.  

Now we can finish clarifying the relationship between balanced
punctured niches and equivalences:

\begin{prop} \label{balanced} \et Suppose that $A$ is an $n$-coherent
$O$-algebra.  Then an $m$-dimensional punctured niche in $A$ is
balanced if and only if the $(n-m)$-functor it defines is an
$(n-m)$-equivalence.   
\end{prop}

Proof - Suppose that an $m$-dimensional punctured $o$-niche $p$ in $A$
defines the $(n-m)$-functor $G$.  Then one can check that 
$p$ is balanced if and only if the punctured $f$-niche of $G$ is
balanced, that is, if and only if $G$ is an $(n-m)$-equivalence.  
 \qed

We conclude by explaining the sense in which a given niche-occupant is
universal if and only if any of its niche-competitors factors
through it, up to equivalence.  Recall that associated to 
any $m$-dimensional $o$-frame 
\[ 
\begin{diagram}[(a_{1},\dots,a_{k})]
\node{(a_{1},\dots,a_{k})}\arrow{e,t}{?} \node{b\;,}
\end{diagram} 
\]
in $A$ there is an $(n-m)$-category.  We denote this by
$\hom_o(a_1,\dots,a_k,b)$, though 
when $o$ is an identity operation, we may follow more traditional
practice and omit it.  

Suppose that in above situation $b'$ is a frame-competitor of $b$.
Then there is an $(n-m)$-category $\hom(b,b')$.  Given any
0-cell
$x \in \hom_o(a_1,\dots,a_k,b)$, there are two virtual $(n-m)$-functors
\[    x_1^\ast,\; x_2^\ast \maps \hom(b,b') \to \hom_o(a_1,\dots,a_k,b') \]
either one of 
which we may think of as `composition with $x$'.  The first is the virtual
$(n-m)$-functor determined by the $(m+1)$-dimensional punctured
niche in $A$,
\[
\begin{diagram}[((c_{1},\dots,c_{k}) \to d,\; d \to d')]
\node{((a_{1},\dots,a_{k}) \mapright{u} b,\; b \mapright{?} b')} 
\arrow{s,r}{?} \\
\node{(a_{1},\dots,a_{k}) \mapright{?} b'}  
\end{diagram}
\]
The second is the one determined by the punctured niche
\[
\begin{diagram}[((c_{1},\dots,c_{k}) \to d,\; d \to d')]
\node{(b \mapright{?} b', \;(a_{1},\dots,a_{k}) \mapright{u} b)} 
\arrow{s,r}{?} \\
\node{(a_{1},\dots,a_{k}) \mapright{?} b'}  
\end{diagram}
\]
We now show that $x$ is universal in its niche if and only if both
these are equivalences --- i.e., heuristically speaking, all the
niche-competitors of $x$ factor through $x$, up to equivalence.   

\begin{prop} \label{universal} \et Suppose that $A$ is an $n$-coherent
$O$-algebra.  Let 
\[ 
\begin{diagram}[(a_{1},a_{j+1},\dots,a_{k})]
\node{(a_{1},\dots,a_{k})}\arrow{e,t}{x} \node{b}
\end{diagram} 
\]
be an occupant of an $m$-dimensional $o$-niche 
\[ 
\begin{diagram}[(a_{1},a_{j+1},\dots,a_{k})]
\node{(a_{1},\dots,a_{k})}\arrow{e,t}{?} \node{?}
\end{diagram} 
\]
Then $x$ is universal if and only if for any frame-competitor $b'$ of 
$b$, the virtual $(n-m)$-functors $x_1^\ast$ and $x_2^\ast$
are $(n-m)$-equivalences.  \end{prop}

Proof - By definition $x$ is universal if and only the 
punctured niches corresponding to $x_1^\ast$ and $x_2^\ast$
above are balanced, or equivalently, by Proposition \ref{balanced},
if $x_1^\ast$ and $x_2^\ast$ are equivalences. \qed

It is a bit annoying to have two virtual $(n-m)$-functors with 
an equal claim to being `composition with $x$', but it is
not very surprising in the present context.  In fact we conjecture
that $x_1^\ast$ is an equivalence if and only if $x_2^\ast$ is.  

\subsection{The microcosm principle}\label{microcosm}

In Section \ref{monoid.objects} we gave a rough statement of the
microcosm principle as follows: {\it certain algebraic structures can
be defined in any category equipped with a categorified version of the
same structure.}   To make this more precise one needs to work with
some particular class of algebraic structures.  Since our approach to
$n$-categories is especially suited to studying operad algebras, we
work with these.   

Recall that for any $S$-operad $O$, a $1$-coherent $O$-algebra can be
thought of as a categorified analog of an $O$-algebra.  Here we show
the following version of the microcosm principle: {\it $O$-algebra
objects can be defined in any $1$-coherent $O$-algebra}.  For example,
monoid objects can be defined in any monoidal 1-category, and
commutative monoid objects can be defined in any stable 1-category.
Another example is the fact that we may define morphisms `riding'
virtual functors.  These are simply $F_1$-algebra objects in
1-coherent $F_1$-algebras.

More generally, we show that {\it $n$-coherent $O$-algebra objects can
be defined in any $(n+1)$-coherent $O$-algebra}.   For example,
`monoidal $n$-category objects' can be defined in any monoidal
$(n+1)$-category, and `stable $n$-category objects' can be defined in
any stable $(n+1)$-category.  

\begin{prop} \et Let $O$ be an $S$-operad.  There exists a terminal
$n$-coherent $O$-algebra $\tau$, that is, one such that for any
$n$-coherent $O$-algebra $A$, there is a unique $n$-coherent
$O$-algebra morphism $f \maps A \to \tau$.  \end{prop}

Proof - Let $\tau$ be a terminal $O$-opetopic set, that is, one having
only one cell occupying each frame, and thus one cell for each
$O$-opetope.   We prove that $\tau$ is an $n$-coherent 
$O$-algebra by showing inductively `from the top down' that every
niche-occupant in $\tau$ is universal, so that every niche has
a universal occupant and composites of universal cells are universal.  
It then follows that $\tau$ is universal as an $n$-coherent
$O$-algebra, since is already terminal as an $O$-opetopic set.  

We claim that every occupant of an $m$-dimensional niche is universal,
and every $m$-dimensional punctured niche is balanced.  By the
definition of $n$-coherent $O$-algebra, both of these are true if $m >
n + 1$.  Supposing they are true for a given $m$, let us show they
hold for $m - 1$.  Given an $(m-1)$-dimensional punctured niche,
condition 1 in the definition of `balanced' holds because every frame
has an occupant, while condition 2 holds by our inductive hypothesis. 
Similarly, every $(m-1)$-dimensional niche-occupant is universal by
our inductive hypothesis.  \qed

\begin{defn}\et Let $O$ be an $S$-operad, let $A$ be an $(n+1)$-coherent
$O$-algebra, and let $\tau$ be the terminal $(n+1)$-coherent $O$-algebra.
Then we define an {\rm $n$-coherent $O$-algebra object in} $A$ to be
a morphism of $O$-opetopic sets $a \maps \tau \to A$.  If $n = 0$,
we call this simply an {\rm $O$-algebra object in} $A$.  
\end{defn}

Since $\tau$ has one cell for each $O$-opetope, we see that 
an $n$-coherent $O$-algebra object in $A$ gives:
\begin{enumerate}
\item a $0$-cell of $A$ for each type of $O$
\item a $1$-cell of $A$ for each operation of $O$
\item a $2$-cell of $A$ for each reduction law of $O$
\item a $3$-cell of $A$ for each way of combining reduction laws of
$O$ to obtain another reduction law
\end{enumerate}
and so on, satisfying certain conditions.  We can work out what this
amounts to quite explicitly in the case $n = 1$.  First we give a
`nuts-and-bolts' description of 1-coherent $O$-algebras:

\begin{thm} \label{1-coherent.o-alg} \et A $1$-coherent $O$-algebra
$A$ consists of:   
{\rm 
\begin{enumerate}
\item {\it for each type $x$ of $O$, a category $A(x)$}
\item {\it for each nondegenerate operation $g$ of $O$ with profile
$(x_1,\dots, x_k,x')$, a $k$-ary virtual functor} 
\[    A(g) \maps A(x_1) \times \cdots \times A(x_k) \rightharpoonup (x') \]
\item {\it for each nondegenerate reduction law of $O$ --- that is, for
each nondegenerate operation $G$ of $O^+$ with profile
$(g_1,\dots, g_k,g')$ --- a natural isomorphism}
\[  A(G) \maps G(A(g_1), \dots ,A(g_k)) \to A(g') \]
\item {\it for each nondegenerate way of combining reduction laws of $O$
to obtain another reduction law --- that is, for each operation $\G$
of $O^{++}$ with profile $(G_1,\dots,G_k,G')$ --- an equation}
\[  \G(A(G_1),\dots, A(G_k)) = A(G') \]
\end{enumerate}
}
\end{thm}

Proof - Note from Example \ref{anafunctor} that a $k$-ary virtual
functor $F \maps C_1 \times \cdots \times C_k \rightharpoonup C'$ is a
special sort of set-valued functor on $C_1^\op \times \dots \times
C_k^\op \times C'$, so the concept of `natural isomorphism' between $k$-ary
virtual functors makes sense.  Also, much as in Theorem
\ref{o+.algebras},  metatree notation makes it clear how to compose
the $k_i$-ary virtual functors $A(g_i)$ in a tree-like pattern
specified by the operation $G$ of $O^+$ to obtain
$G(A(g_1),\dots,A(g_k))$, and how to compose the natural isomorphisms
$A(G_i)$ in a tree-like pattern specified by the operation $\G$ of
$O^{++}$ to obtain a natural isomorphism $\G(A(G_1),\dots,A(G_k))$.   

By item 1 of Theorem \ref{n-coherent.o-alg}, the 1-coherent
$O$-algebra $A$ has an underlying 0-coherent $(O_{A(0)})^+$-algebra,
which  we may think of as simply an $(O_{A(0)})^+$-algebra.  Theorem
\ref{o+.algebras} implies that such an algebra consists of 1-4 as
above, but with a $k$-ary distributor $A(g)$ for each nondegenerate
$k$-ary operation $g$ of $O$, and a natural transformation $A(G)$
between $k$-ary distributors for each nondegenerate operation $G$ of
$O^+$.  Item 2 of Theorem \ref{n-coherent.o-alg} implies that the
$A(g)$ are $k$-ary virtual functors, and item 3 of that theorem
implies that the $A(G)$ are natural isomorphisms.  Conversely, one
can show that 1-4 as above give a 1-coherent $O$-algebra.  \qed

In particular, we see that monoidal 1-categories and stable
1-categories are almost the same as monoidal categories and symmetric
monoidal categories, respectively, though there is a bit of work
required to translate between our concepts and the traditional ones.  

To describe $O$-algebra objects in the language of Theorem
\ref{1-coherent.o-alg}, it is convenient to define a {\it $k$-ary morphism
$b \maps c_1 \times \cdots \times c_k \to c'$ riding} the $k$-ary
virtual functor $B \maps C_1 \times \cdots \times C_k \to C'$ to be an
$F_k$-algebra object $b$ in the 1-coherent $F_k$-algebra $H$.
Concretely, this amounts to a choice of objects $c_i \in C_i$ and $c'
\in C'$, together with a $0$-cell $b$ in
$\hom_f(c_1,\dots,c_k,c')$.

\begin{thm} \label{o-alg.object} \et Given a 1-coherent $O$-algebra
$A$, an $O$-algebra object $a$ in $A$ consists of: 
{\rm
\begin{enumerate}
\item {\it for each type $x$ of $O$, an object $a(x)$ in the category $A(x)$}
\item {\it for each nondegenerate operation $g$ of $O$ with profile
$(x_1,\dots,x_k,x')$, a $k$-ary morphism $a(g) \maps a(x_1) \times \cdots
\times a(x_k) \to a(x')$ riding the $k$-ary virtual functor $A(g)$}
\item {\it for each nondegenerate reduction law $G$ of $O$ with profile
$(g_1,\dots,g_k,g')$, an equation}
\[  A(G)(a(g_1),\dots, a(g_k)) = a(g') \]
\end{enumerate} 
}
\end{thm}

Proof - This is straightforward except that item 3 may need some
clarification.  Given $\ell_i$-ary morphisms $a(g_i)$ riding the
$\ell_i$-ary virtual functors $A(g_i)$, and
a reduction law $G$ of $O$ with profile $(g_1,\dots,g_k,g')$, one
obtains an $\ell'$-ary morphism riding the $\ell'$-ary virtual functor 
$G(A(g_1),\dots,A(g_k))$.  Applying the natural isomorphism $A(G)$
to this we obtain an $\ell'$-ary morphism riding $A(g')$, which we
call $A(G)(a(g_1),\dots,a(g_k))$.  In item 3 we require this to equal
$a(g')$.  \qed

\section{Conclusions}

In addition to our approach to weak $n$-categories, there are a number
of others.  We have already mentioned Street's original simplicial
approach \cite{Street}.  After a sketch of our definition appeared
\cite{BD2}, Makkai has begun studying it, and a version adapted to the
globular approach has been developed by Makkai, Hermida, and Power,
but the details of this have not yet been published.  Independently,
Tamsamani \cite{Tamsamani} developed an approach using multisimplicial
sets: simplicial objects in the category of simplicial objects in the
category of simplical objects in the category of $\dots$ sets.  More
recently, Batanin has developed a globular approach to weak
$\omega$-categories, and thus in particular weak $n$-categories, using
the notion of an `$\omega$-operad'.  We expect that as time goes by
even more definitions will be proposed.

The question thus arises of when two definitions of weak $n$-category
may be considered `equivalent'.  This question was already raised, and
a solution proposed, in Grothendieck's 600-page letter to Quillen
\cite{Gro}.  Suppose that for all $n$ we have two different
definitions of weak $n$-category, say `$n$-${\rm category}_1$' and
`$n$-${\rm category}_2$'.  Then we should try to construct the
$(n+1)$-${\rm category}_1$ of all $n$-${\rm categories}_1$ and the
$(n+1)$-${\rm category}_1$ of all $n$-${\rm categories}_2$ and see if
these are equivalent as objects of the $(n+2)$-${\rm category}_1$ of
all $(n+1)$-${\rm categories}_1$.  If so, we may say the two
definitions are equivalent as seen from the viewpoint of the
first definition.   Of course, there are some `size' issues
involved here, but they should not be a serious problem.  More
importantly, there is some freedom of choice involved in constructing
the two $(n+1)$-${\rm categories}_1$ in question.  Also, we would be in
an embarrassing position if we got a different answer for the question
with the roles of the two definitions reversed.  Nonetheless, it
should be interesting to compare different definitions of weak
$n$-category in this way.

A second solution is suggested by homotopy theory, where many
superficially different approaches turn out to be fundamentally
equivalent.  Different approaches use objects from different `model
categories' to represent homotopy types: compactly generated
topological spaces, CW complexes, Kan complexes, and so on
\cite{Baues,Brown}.  These categories are not equivalent, but each one is
equipped with a class of morphisms playing the role of
homotopy equivalences.  Given a category $C$ equipped with a specified
class of morphisms called `equivalences', under mild assumptions one
can `localize' $C$ with respect to this class, which amounts to
adjoining inverses for these morphisms \cite{GZ}.  The resulting
category is called the `homotopy category' of $C$.  Two categories
with specified equivalences may be considered the same for the
purposes of homotopy theory if their homotopy categories are
equivalent.  All the model categories above are the same in this
sense.

It is natural to adopt the same attitude in $n$-category theory. 
(Indeed, this attitude is also implicit in Grothendieck's letter to 
Quillen, which was in part inspired by the latter's work on model
categories \cite{Quillen}.)  Thus we propose the following homotopy
category of $n$-categories.  We define an $n$-functor $F \maps C \to
D$ to be an {\it equivalence} if:
\begin{enumerate}
\item Every 0-cell in $C$ is connected to a 0-cell in
the image of $F$ by a universal 1-cell.
\item For any 0-cells $c,c'$ in $C$, the restriction of $F$ to
the $(n-1)$-category $\hom(c,c')$ is an equivalence.
\end{enumerate}
\noindent where to ground this recursive definition we define equivalences
between $0$-categories to be bijections, using the identification
of $0$-categories with sets.   Condition 1 above says that  $F$ is
`essentially surjective', while condition 2 says that $F$ is `fully
faithful'.   We then define the following category:

\begin{defn} \et The {\rm homotopy category of $n$-categories} is the
localization of the category of $n$-categories and $n$-functors with
respect to the equivalences. \end{defn}

\noindent We regard any other definition of $n$-category as
fundamentally `the same' as ours if it gives an equivalent homotopy
category of $n$-categories.  

\subsection*{Acknowledgements}

We wish to thank Michael Batanin, Sjoerd Crans, John Power, Michael
Makkai, Bruce Smith, James Stasheff, Ross Street, Zouhair Tamsamani,
Todd Trimble, and Dominic Verity for discussions and correspondence on
the subject of $n$-categories.

\end{document}